\documentclass[oneside, a4paper,11pt,reqno]{amsart}
\makeatletter
\def\timenow{\@tempcnta\time
  \@tempcntb\@tempcnta
  \divide\@tempcntb60
  \ifnum10>\@tempcntb0\fi\number\@tempcntb
  \multiply\@tempcntb60
  \advance\@tempcnta-\@tempcntb
  :\ifnum10>\@tempcnta0\fi\number\@tempcnta}
\makeatother

\usepackage[ddmmyyyy,hhmmss]{datetime}
\usepackage{color}

\usepackage{euscript}

\usepackage[applemac]{inputenc}
\usepackage{amssymb,amsmath,amsthm,dsfont}
\usepackage{graphics,graphicx}
\usepackage[T1]{fontenc}
\usepackage{color}
\usepackage{epsfig,psfrag}

\newtheorem{theo}{Theorem}[section]
\newtheorem{prop}[theo]{Proposition}
\newtheorem{lemma}[theo]{Lemma}

\newtheorem{cor}[theo]{Corollary}

\newtheorem{rem}[theo]{Remark}
\newtheorem{fact}[theo]{Fact}

\newcommand{\egloi}{\stackrel{\mathcal{L}}{=}}
\newcommand{\cvloi}{\stackrel{\mathcal{L}}{\rightarrow}}

\newcommand{\mysection}{\setcounter{equation}{0} \section}

\def\un{\mathds{1}}
\def\a{\alpha}

\def\d{\delta}
\def\e{\varepsilon}

\def\E{\mathbb{E}}
\def\B{\mathcal{E}}

\def\o{\omega}

\def\N{\mathbb{N}}
\def\P{\mathbb{P}}

\def\R{\mathbb{R}}

\def\tt{\theta^-}

\def\Z{\mathbb{Z}}
\def\Vi{{V}^{(i)}}
\def\tV{\tilde{V}^{(i)}}

\def\mV{\mathcal{V}}
\def\tA2{\tilde{A}_2\big({\tilde L_2}\big)}
\def \tts2{\tilde \tau_2(  h_t/2 )}

\def\mI{{\bf I}}
\def\Io{{\mathcal{I}}}
\def\Ip{{\mathcal{I}}^+}
\def\Im{{\mathcal{I}}^-}

\def\tt{\tilde{\tau}}

\def\mt{\tilde{m}}
\def\Mt{\tilde{M}}
\def\mf{{m}}
\def\Mf{{M}}

\def\Vi{V^{(i)}}

\def\g{\zeta}

\def\tm{\tilde m_i}
\def\tL{\tilde L_j}
\def\Lm{\tilde L_j^-}

\def\tN{\tilde N}

\def\lo{\mathcal L}

\def\bU{ {\bf U} }
\def\k{\kappa}
\def\dd{\textnormal{d}}
\def\wk{W_{\k}}
\def\eo{\E^{W_{\kappa}}}

\def\lneg{\lo_X^{*-}(+\infty)}
\def\lx{\lo_X}

\title[Localization and number of visited valleys for a transient diffusion in random environment]
{Localization and number of visited valleys for a transient diffusion in random environment} 

\author{Pierre Andreoletti}
\address{Laboratoire MAPMO - Fédération Denis Poisson, Universit\'e d'Orl\'eans, B\^atiment de math\'ematiques - Rue de Chartres
B.P. 6759 - 45067 Orl\'eans cedex 2,
France}
\email{pierre.andreoletti@univ-orleans.fr}

\author{Alexis Devulder}
\address{Universit\'e de Versailles Saint-Quentin-en-Yvelines, Laboratoire de Math\'ematiques de
Versailles, CNRS UMR 8100, B\`at. Fermat,
45 avenue des Etats-Unis,
78035 Versailles Cedex, France.}
\email{devulder@math.uvsq.fr }

\keywords{Diffusion in a random potential, Localization, Aging, Renewal theorem.\\
This research was supported by the french ANR projects MEMEMO 2 (2010 BLAN 0125 02) and Randymeca}

\date{\today}

\begin{document}

\pagestyle{plain}

\begin{abstract}
\underline{}We consider a transient diffusion in a $(-\kappa/2)$-drifted Brownian potential $W_{\kappa}$ with $0<\kappa<1$.
 We prove its localization at time $t$ in the neighborhood of some random points depending only on the environment, which are
 the positive $h_t$-minima of the environment, for $h_t$ a bit smaller than $\log t$.
 We also prove an Aging phenomenon for the diffusion,
 a renewal theorem for the hitting time of the farthest visited valley, and provide a  central limit theorem for the number of valleys  visited up to time $t$.

The proof relies on a
decomposition of the trajectory of $\wk$ in the neighborhood of
$h_t$-minima, with the help of results of Faggionato \cite{Faggionato},
and on a precise analysis of exponential functionals of $W_{\kappa}$ and of
$\wk$ Doob-conditioned to stay positive.
\end{abstract}

\maketitle



\mysection{Introduction and notation}\label{intro}

\subsection{Presentation of the model}
We are interested in a diffusion $(X(t),\ t\geq 0)$ in a random c\`adl\`ag potential $(V(x),\ x\in\R)$.
It is defined informally by $X(0)=0$ and
$$\text{d}X(t)=d\beta(t)-\frac{1}{2}V'(X(t))\text{d}t,$$
where $(\beta(t),\ t\geq 0)$ is a Brownian motion independent of $V$. More rigorously, $X$ is
a diffusion process, starting from $0$, and whose conditional generator given $V$ is
$$
\frac{1}{2}e^{V(x)}\frac{\text{d}}{\text{d} x}\left(e^{-V(x)}\frac{\text{d}}{\text{d} x}\right).
$$
These diffusions in random potentials are considered as continuous time analogues of random walks in random environment (RWRE)
(see e.g. P. R\'ev\'esz \cite{Revesz}, B.D. Hughes \cite{Hug}, Z. Shi \cite{Shi1} and O. Zeitouni \cite{Zeitouni} for reviews on RWRE).

The study of such a process starts with a choice for $V$. A classic one, originally introduced by S. Schumacher  \cite{Schumacher} and T. Brox \cite{Brox}, is to choose $V$ as a L\'evy process. In fact only a few papers deal with the discontinuous case, see for example P. Carmona \cite{Carmona} or A. Singh \cite{Singh, Singh1}. Most of the results concern diffusions in a continuous L\'evy potential $V$, that is,
$$
V(x)=\wk(x):=W(x)-\frac{\k}{2}x, \qquad x\in\R,
$$
where $\kappa \in \R$ and $(W(x),\ x\in\R)$ is a two sided Brownian motion.
We denote by $P$ the probability measure associated to $\wk(.)$.
The probability conditionally on the potential $W_{\kappa}$ is denoted by $\P^{W_{\kappa}}$
and is called the {\it quenched probability}.
 We also define the {\it annealed probability} as
$$
    \P(.)
:=
    \int \P^{W_{\kappa}}(.) P(\wk\in \dd\o).
$$
We denote respectively by $\E^{W_{\kappa}}$, $\E$, and $E$ the expectations with regard to $\P^{W_{\kappa}}$, $\P$ and $P$.

In the case $\kappa=0$, the diffusion $X$ is a.s. recurrent. More precisely, T. Brox \cite{Brox} shows that it is  sub-diffusive with asymptotic behavior in $(\log t)^2$, and that it is localized, at time $t$,  in the neighborhood of a random point $b_{\log t}$ depending only on $t$ and $W$,
similarly as Sinai's walk (see Ya. G. Sinai \cite{Sinai}).
More precisely, this result can be written:
\begin{theo} (Brox \cite{Brox})  Assume $\kappa=0$. Then, for all $\e>0$,
\begin{align}
\lim_{t \rightarrow + \infty} \P\left[ X(t) \in [b_{\log t}-\e (\log t)^2, b_{\log t}+\e (\log t)^2] \right]=1. \label{Brox}
\end{align}
\end{theo}
\noindent The limit law of $b_{\log t}/(\log t)^2$ and therefore of $X(t)/(\log t)^2$ was made explicit independently by H. Kesten \cite{Kesten2} and A. O. Golosov \cite{Golosov1}. For recent results for this recurrent case, see for example P. Andreoletti et al. \cite{AndeolettiDiel}, and R. Diel \cite{Diel}.

In the case $\kappa\neq 0$, the diffusion $X$ is a.s. transient, with a wide range of limiting behaviors, depending on the value of $\kappa$. It was first studied by K. Kawazu and H. Tanaka. Let us denote by $H(r)$ the hitting time of $r\in\R$ by $X$:
$$H(r):=\inf\{s>0,\ X(s)=r\}.$$
Kawazu et al. \cite{KawazuTanaka} proved in particular that
under the annealed probability $\P$,   $H(r)/r^{1/\k}$ converges in law to a stable distribution when $0<\k<1$,
whereas $H(r)/(r\log r)$ converges in probability to $4$ when $\k=1$,
and $H(r)/r$ converges almost surely to $4/(\k-1)$ if $k>1$
 (see also Y. Hu et al. \cite{HuShiYor}, and H. Tanaka \cite{Tanaka2}).
More recently we mention the results for large and moderate deviations, by M. Taleb (\cite{Taleb0} and \cite{Talet}), A. Devulder \cite{Devulder2}
 and G. Faraud \cite{Faraud}.


In this paper we study the case $0<\k<1$. We follow a different approach from Y. Hu et al. \cite{HuShiYor} and K. Kawazu et al. \cite{KawazuTanaka}. Indeed we focus on a quenched study, which has attracted
much interest for transient RWRE in the last few years, see for example the works of N. Enriquez et al. \cite{ESZ4}, \cite{ESZ1}, \cite{ESZ2}, \cite{ESZ3}, D. Dolgopyat et al.  \cite{DoGo}, and J. Peterson et al.
\cite{PetersonSamo2010}, \cite{PetersonSamo}, \cite{PetersonZeitouni}.
Heuristically, the diffusion $X$ goes to locations where the potential is low, hence it goes to $+\infty$, but it is slowed by "valleys" of the potential, which trap the diffusion for some time. The
 diffusion even spends most of its time in these valleys.
We will prove this more in details in the present paper.


\subsection{Main results}
The goals of this paper are to localize the diffusion $X$, when $ 0 <\kappa<1 $, in some valleys of the potential $\wk$,
to understand the differences with Brox's result given by \eqref{Brox}, and to prove an Aging phenomenon,
corresponding to results obtained by Enriquez et al. in their
papers \cite{ESZ1}, \cite{ESZ2} and  \cite{ESZ3} for transient zero-speed RWRE.
We moreover obtain a central limit theorem for the number of valleys visited up to time $t$. We also prove some intermediate results, which we think will be useful for obtaining new results about the maximum local time of $X$, as explained later in  this introduction.



Let $t \mapsto \phi(t)$ be a positive increasing function, such that
$\phi(t)=o(\log t)$ and $\log\log t=o(\phi(t))$ as $t\to+\infty$,
where $f(t)=o(g(t))$ means $\lim_{t\to+\infty}f(t)/g(t)=0$. We prove the following aging phenomenon:

\begin{prop}\label{ThAging}
 Assume $0 <\kappa <1$. For all $\alpha>1$, we have
\begin{align*}
    \lim_{t \rightarrow + \infty}
    \P\Big( \left|X(\alpha t)-X(t) \right| \leq    \phi(t) \Big)
=
    \frac{\sin(\kappa \pi)}{\pi} \int_{0}^{1/\alpha}u^{\kappa-1}(1-u)^{-\kappa}\textnormal{d}u.
\end{align*}
\end{prop}

\noindent

More generally, {\it aging} usually denotes dynamical out-of-equilibrium physical phenomenons, which appear in some disordered systems.
It refers to the existence of a limit for a given two-time correlation function of the system as both times diverge but keep a fixed ratio between them. This subject has received a considerable attention in physics.
For a physical or a mathematical point of view on aging, see e.g. respectively Bouchaud et al. \cite{BouchaudAging}
and Zindy \cite{Zindy}, and references therein.

Proposition \ref{ThAging} is actually a consequence of Theorem \ref{ththm}.
Before stating it, we first recall the notion of {\it $h$-extrema},
which was first introduced by J. Neveu et al. \cite{NP}, and studied in the case of drifted Brownian motion by
A. Faggionato \cite{Faggionato}. For $h>0$, we say that $x\in\R$ is an {\it $h$-minimum}
for a given function $f$, $\R\to \R$, if there exist $u<x<v$ such that
$f(x)=\inf_{y\in[u,v]} f(y)$, $f(u)\geq f(x)+h$ and $f(v)\geq f(x)+h$.
Moreover, $x$ is an {\it $h$-maximum} for $f$ iff $x$ is an $h$-minimum for $-f$.
Finally, $x$ is an {\it $h$-extremum} for $f$ iff it is an $h$-maximum or an $h$-minimum for $f$.

Since we want to study the diffusion $X$ until time $t>0$,
we are more especially interested in the $h_t$-extrema of $\wk$, where
$$
h_t:=\log t-\phi(t).
$$
It is known (see \cite{Faggionato}) that almost surely, the $h_t$-extrema of $\wk$ form a sequence indexed by
$\Z$, unbounded from below and above, and that the $h_t$-minima and $h_t$-maxima alternate.
We denote respectively by $(\mf_j,\ j \in \Z)$ and $(\Mf_j,\ j \in \Z)$ the increasing sequences of $h_t$-minima
and of $h_t$-maxima of $\wk$,
such that $\mf_{0}\leq 0<\mf_1$ and $\mf_j<\Mf_j<\mf_{j+1}$ for every  $j\in\Z$.
These $h_t$-minima $m_i$, $i\in\Z$, can be considered as the bottoms of some valleys of the potential $\wk$,
of height at least $h_t$, that will be defined more precisely in Section \ref{environments}.
We also introduce
$$
N_t:=\max\Big\{k\in\N,\ \sup_{ 0 \leq s \leq t}X(s) \geq \mf_{k}\Big\},
$$
so that $m_{N_t}$ is the largest $h_t$-minimum visited by $X$ until time $t$
if $N_t>0$ (and $\lim_{t\to+\infty}\P(N_t>0)=1$ as proved later, see Proposition \ref{propTRT2} or Lemma \ref{LemmaMajorationNt}).
The main result of this paper concerns the localization of the diffusion. It is stated as follows:
\begin{theo} \label{ththm}
Assume $ 0<\kappa<1$. There exists a constant $ \mathcal{C}_1 >0$, such that
$$
    \lim_{t \rightarrow + \infty }
    \P\Big( |X(t)- \mf_{N_{t}}| \leq \mathcal{C}_1 \phi(t) \Big)
=
    1.
$$
\end{theo}

We first recall that $X(t)$ is asymptotically of order $t^{\kappa}$ (see Kawazu et al. \cite{KawazuTanaka}),
and that the typical distance between two $h_t$-minima of $\wk$ is asymptotically of order $e^{\k h_t}=t^{\k}e^{-\k \phi(t)}$
(see Faggionato \cite{Faggionato} Prop. 1, partly recalled in our Fact \ref{Fact_Faggio} below).
So, the size $2 \mathcal{C}_1\phi(t)$ of the intervals in which $X$ is localized in Theorem \ref{ththm}
is very small, since it is $o(\log t)$ and can be as small as, for example,  $(\log \log t)^{1+\e}$, $\e>0$.
Notice that it depends on the minimum height $h_t$ of our valleys.
We could not say however if it the best interval size that can be obtained.\\
\indent
The main difference with the result of Brox \eqref{Brox} is the appearance of the (random) integer $N_t$, which is the number of typical valleys of height $h_t$ visited before time $t$. In the recurrent case of Brox, the diffusion $X$ is, with a large probability, localized near the bottom of a unique valley of the potential, whereas in our transient case, the diffusion is localized near the bottom of one among several valleys of the potential. This, and the absence of scaling for the potential in the case $0<\k<1$, contrarily to the case $\k=0$, makes the study much more involved technically.


We also prove a renewal theorem for hitting time of the bottom $m_{N_t}$ of the last
valley
visited by $X$ before $t$:


\begin{prop}\label{propRenewalJoint}
Assume $0<\k<1$.
We have the following convergence in law under the annealed probability $\P$,
$$
    \bigg(\frac{H(\mf_{N_t})}{t}, \frac{H(\mf_{N_t+1})}{t}\bigg)
\cvloi_{t\to+\infty}
    \frac{\k\sin(\pi\k)}{\pi}(y-x)^{-\k-1}x^{\k-1}\un_{[0,1]}(x)\un_{[1,\infty)}(y)\dd x \dd y.
$$
\end{prop}
This unables us to get
the following results, which are useful for the proofs of Proposition \ref{ThAging}
and Theorem \ref{ththm}:

\begin{cor} \label{propTRT}
Assume $ 0<\kappa<1$ and  let\,  $0\leq r <s\leq 1$ and $v\geq 0$.
Then,
\begin{eqnarray}
    \lim_{t \rightarrow + \infty}
    \P\bigg(1-s\leq  \frac{H(\mf_{N_t})}{t} \leq 1-r\bigg)
& = &
    \frac{\sin(\pi \kappa )}{ \pi}
    \int_{1-s}^{1-r} x^{ \kappa-1} (1-x)^{-\kappa} \textnormal{d}x
\label{HNt},
\\
    \lim_{t \rightarrow + \infty} \P\bigg(\frac{H(\mf_{N_t+1})}{t}\geq 1+ v\bigg)
& = &
    \frac{\sin(\pi \kappa )}{ \pi}
    \int_{v}^{+ \infty} (1+x)^{-1} x^{-\kappa}\dd x
\label{Nt}.
\end{eqnarray}
Moreover, the total time spent in the last valley of height at least $h_t$ visited before time $t$ renormalized by $t$,
that is
$[H(\mf_{N_t+1})-H(\mf_{N_t})]/t$, converges in law under $\P$ to a r.v. with density
$
\sin(\pi\k)\pi^{-1}x^{-\k-1}[(1-(1-x)^\k)\un_{[0,1]}(x)+\un_{(1,+\infty)}(x)]
$.
\end{cor}

\noindent
Let $\lfloor x \rfloor$ denote the integer part of $x$, for any $x\in\R$. We introduce
$0<\delta<1$ and
\begin{equation}\label{eqDefnt}
    n_t
:=
    \lfloor e^{\kappa \phi(t)(1+ \delta)} \rfloor.
\end{equation}
We will see in Section 4 and 5 that
Proposition \ref{propRenewalJoint} is a consequence of the fact that for any integer $1\leq k \leq n_t$,
the hitting time $H(m_{k})$ can approximated by a sum of i.i.d. random variables
having  the law of a r.v. $\bU$. This r.v. $\bU$ is an approximation of the time the diffusion $X$ spends in a typical valley of height at least $h_t$ before escaping this valley.

These results are in accordance with those obtained by Enriquez et al. in their three papers \cite{ESZ1}, \cite{ESZ2} and \cite{ESZ3} for transient RWRE. Compared to their study, we have the advantage of being able to use some powerful stochastic tools. However, some other technical difficulties
appear in continuous time: for example local time and excursions are more complicated to deal with in continuous time than in discrete time.
The present paper is self contained, in particular we provide in this same paper the technical study of the Laplace transform of the first exit time $\bU$ of a typical valley. The study of the environment mainly requires continuous arguments of stochastic calculus,
starting by a
decomposition of the trajectory of $\wk$ near its $h_t$-minima, which mainly comes from results of A. Faggionato \cite{Faggionato}.

The number $N_t$ of valleys having height at least $h_t$, visited before time $t$ by the diffusion $X$,  goes to $+ \infty$ as $t\to+\infty$.
However, we prove that $\P(N_t\leq n_t)\to_{t\to+\infty}1$, which explains why we study the potential $\wk(x)$ only for $x\leq m_{n_t}$,
and the hitting times $H(m_k)$ only for $k\leq n_t$.
More precisely, we prove the following central limit theorem for $N_t$, with renormalization $e^{\kappa \phi(t)}$:

\begin{prop} \label{propTRT2}
Assume $ 0<\kappa<1$.
Then ${N_t} e^{-\kappa \phi(t)}\to_{t\to+\infty}\mathcal{N}$ in law under the annealed law $\P$;
the law of $\mathcal{N}$  is determined by its Laplace transform:
\begin{equation}
\label{eqTransfoLaplaceLimiteNbValleesVisitees}
    \forall u>0,
\qquad
    \E\left(e^{-u\mathcal{N}}\right)
=
    \sum_{j=0}^{+ \infty}
    \frac{1}{\Gamma(\kappa j+1) }
    \left(\frac{ -u}{C_{\kappa}}\right)^j
\end{equation}
where $C_{\kappa}>0$ is explicitly known (see Proposition \ref{proplap}).
This r.v. $\mathcal{N}$ has then a Mittag-Leffler distribution of order $\k$.
\end{prop}


Moreover we expect that the results of this paper will be useful to study other properties for the diffusion.
In particular,
let $(\lo_X(t,x),\ t\geq 0,\ x\in\R)$ be a bicontinuous version of the local time of $X$. It is known that the maximum local time of $X$ at time $t$, that is $\lo_X^*(t):=\max_{x\in \R} \lo_X(t,x)$, satisfies $\limsup_{t\to+\infty}\lo_X^*(t)/t=+\infty$ a.s. in the cases $\k=0$
(see Z. Shi \cite{Shi} and R. Diel \cite{Diel}) and even in the transient case $0<\k<1$ (see A. Devulder \cite{Devulder}). Hence the maximum local time of $X$ exhibits very interesting properties, that contrast with those of the maximum local time of RWRE at time $t$, which is naturally bounded by $t/2$.

We especially think that the better understanding of the localization of $X$ and some intermediate results provided in this paper will be useful to prove new results about $\lo_X^*$.
Indeed, in a work in progress
with G. Vechambre \cite{AndeolettiDevulderVechambre},
we use the methods and results of the present paper to study the local time of the diffusion.
In particular we expect to obtain the limit law of the maximal local time after suitable renormalization in the case $0<\kappa<1$.
Note that the local time plays a crucial role for estimation problems for random walk in random environment, recently studied e.g. in \cite{Pierre7}, \cite{AndDiel2} and \cite{Comets_etal2}.
In particular in \cite{Comets_etal2}, the limit law of the local time process in the neighborhood of the minima, obtained previously by \cite{GanPerShi}, is used.
In the same way
a better understanding of the local time of $X$ may be useful
to study estimation problems for diffusions in a random potential.


\subsection{Sketch of the proof and organization of the paper}\label{SubSectSketchOfTheProof}
We now give a general idea of the proof and provide, at the same time, the organization of the paper.
This subsection contains some non rigorous heuristics which will be made rigorous and explained in details
in the following sections.

First, in Section \ref{environments}, we build some valleys
$(\tilde L_i^- , \tilde m_i, \tilde L_i)$
of the potential,
of height at least $h_t$ (see Figure \ref{fig1} page \pageref{fig1}).
The $i$-th valley is the potential $\wk$, restricted to some interval $[\tilde L_i^-, \tilde L_i]$;
the minimum of $\wk$ in this interval is attained at a location called $\tilde m_i$.
The height of this valley is at least $h_t$, and more precisely,
$\wk(\tilde L_i)-\wk(\tilde m_i)\geq h_t$,
$\wk(\tilde L_i^-)-\wk(\tilde m_i)\geq (1+\k+2\delta)h_t$,
These valleys, when recentered at $\tilde m_i$, are i.i.d.
We prove in Lemma \ref{CVs} that with probability nearly $1$, $\tilde m_i=m_i$ for $1\leq i \leq n_t$,
that is, the $\tilde m_i$, $1\leq i \leq n_t$ are in fact the $n_t$ first positive $h_t$-minima.
We also provide in this Section \ref{environments}
 different tools to study the law of $\wk$ near $\tilde m_i$ or $m_i$,
and in particular drifted Brownian motions Doob-conditioned to stay positive.

Then in Section \ref{SectionQuasiIndep}, we prove (see Lemma \ref{LemmaProbaQuitterVallesaDroite}) that
with probability nearly $1$, after hitting the bottom $\tilde m_i$ of
the each valley $(\tilde L_i^- , \tilde m_i, \tilde L_i)$
the diffusion $X$ leaves this valley on the right, that is, on
$\tilde L_i$.
Moreover, we prove (see Lemma \ref{LemmaProbaRetourEnmi})
that with probability nearly $1$,
the diffusion $X$ visits successively the valleys
$(\tilde L_1^- , \tilde m_1, \tilde L_1)$, $(\tilde L_2^- , \tilde m_2, \tilde L_2)$, $\dots$,
$(\tilde L_i^- , \tilde m_i, \tilde L_i)$, $\dots$,
and does not come back to previously visited valleys.
Moreover, we show in Lemma \ref{lemtps} that  the time spent outside these valleys up to time $t$ is negligible compared to $t$.
That is, the hitting time $H(\tilde m_n)$ of the bottom $\tilde m_n$
of each valley number $n\leq n_t$ can be approximated as
\begin{equation}\label{eqApproxHmnparSommeUi}
    H(\tilde m_n)
\approx
    U_1+U_2+\dots+ U_{n-1},
\end{equation}
where $U_i$ can be seen as the time spent in valley number $i$.
We also prove in Lemma \ref{propind}
that these r.v. $U_i$ are "nearly" independent under the annealed and quenched probabilities $\P$ and $\P^{\wk}$.
Moreover, we prove that they have (under the annealed probability $\P$) the same law as some random variable $\bU$, defined by
\begin{equation*}
    \bU
:=
    \int_{\tilde L_2^-}^{\tilde L_2}
    e^{-[W_{\kappa}(u)-W_{\kappa}(\tilde m_2)]}
    \lo_B\big[\tau^B\big(\tilde{A}_2\big(\tilde L_2\big)\big),\tilde{A}_2(u)\big] \textnormal{d}u,
\end{equation*}
where $\tilde{A}_2(z):=\int_{\tilde m_2}^{z} e^{W_{\kappa}(x)-W_{\kappa}(\tilde m_2)}\textnormal{d}x$, $z\in\R$,
and $\mathcal{L}_B$ and $\tau^B$ denote respectively the local time and hitting time of some Brownian motion $B$ independent of $\wk$.
This r.v. $\bU$, which depends on $t$,  can be seen as the time spent by $X$ in a typical valley, which $X$ leaves on its right.

So, this reduces the study of $H(\tilde m_n)$ to the study of a sum $U_1+U_2+\dots+ U_{n-1}$
of independent random variables having the same law as $\bU$ under $\P$.

The main goal of Section \ref{SectionTimeStandardValley} is to study the asymptotics of the Laplace transform of $\bU$.
More precisely, we prove in Proposition \ref{proplap} that for all $\lambda>0$,
\begin{equation}\label{eqApproxLaplaceU}
     \E\big(e^{-\lambda \bU/t}\big)
=
    1-C_{\kappa}\lambda^{\kappa}e^{ -\kappa \phi(t)}+o(e^{ -\kappa \phi(t)})
\end{equation}
as $t\to+\infty$. To this aim, we first approximate $\bU$ by a more tractable expression,
which is done in Proposition \ref{lemU1} and Lemma \ref{propIV}.
First, by scaling, $\bU$ is equal in law to
\begin{equation}\label{eqApproxUSktetch}
    \int_{\tilde L_2^-}^{\tilde L_2}
    e^{-[W_{\kappa}(u)-W_{\kappa}(\tilde m_2)]}
    \tilde{A}_2\big(\tilde L_2\big)
    \lo_B\big[\tau^B(1),\tilde{A}_2(u)/\tilde{A}_2\big(\tilde L_2\big)\big] \textnormal{d}u.
\end{equation}
Loosely speaking, for $u$ "close" to $\tilde m_2$,
$\tilde{A}_2(u)/\tilde{A}_2(\tilde L_2)\approx 0$, and so (using our Lemma \ref{LemmaContinuiteEn0}) we have
$
    \lo_B[\tau^B(1),\tilde{A}_2(u)/\tilde{A}_2(\tilde L_2)]
\approx
    \lo_B[\tau^B(1),0]
=:
    {\bf e_1}
$,
which is by the first Ray Knight theorem an exponential variable of mean $2$, and is independent of $\wk$.

One the other hand, $\tilde m_2$ is the minimum of the potential $\wk$ in
$[\tilde L_2^-,  \tilde L_2]$
so for $u$ "far" from $\tilde m_2$  we have
$
    e^{-[W_{\kappa}(u)-W_{\kappa}(\mf_2)]}
\approx
    0
$.
And moreover,
$
    \tilde{A}_2(\tilde L_2)
\approx
    \int_{\tilde \tau_2(h_t/2)}^{\tilde L_2}
    e^{W_{\kappa}(x)-W_{\kappa}(\tilde m_2)}\textnormal{d}x
$
where $\tilde \tau_2(h)$ is the first hitting time of $h>0$
by $W_{\kappa}-W_{\kappa}(\tilde m_2)$ after $\tilde m_2$, because
$e^{W_{\kappa}(x)-W_{\kappa}(\tilde m_2)}$ is negligible compared to
$\tilde{A}_2(\tilde L_2)$ for $\tilde m_2\leq x \leq \tilde \tau_2(h_t/2)\leq \tilde \tau_2(h_t)<\tilde L_2$.
All this leads to the approximation
\begin{eqnarray*}
    \bU
& \approx &
    \left(
    \int_{\tilde L_2^-}^{\tilde m_2}
    +
    \int_{\tilde m_2}^{\tilde \tau_2(h_t/2)}
    e^{-[W_{\kappa}(u)-W_{\kappa}(\mf_2)]}
    \textnormal{d}u
    \right)
    \left(
    \int_{\tilde \tau_2(h_t/2)}^{\tilde \tau_2(h_t)}
    +
    \int_{\tilde \tau_2(h_t)}^{\tilde L_2}
    e^{W_{\kappa}(x)-W_{\kappa}(\tilde m_2)}\textnormal{d}x
    \right)
    {\bf e_1}
\\
& \approx: &
    (\Im_1+\Im_2)
    (\Ip_1+\Ip_2)
    {\bf e_1},
\end{eqnarray*}
which is the product of $5$ independent random variables, the first $4$ depending only on the potential $\wk$,
and ${\bf e_1}$ being independent of $\wk$.
This approximation, the asymptotics of the Laplace transforms of $\Im_1$, $\Im_2$, $\Ip_1$ and $\Ip_2$
provided by Lemma \ref{lem4.5}, and some technical calculations
help us to prove \eqref{eqApproxLaplaceU} as claimed in our Proposition \ref{proplap}.

Section \ref{SectionPreuveLocEtAging} is devoted to the proofs of the main results of this paper.
First, using the asymptotics  \eqref{eqApproxLaplaceU} of the Laplace transform of $\bU$ as well as
\eqref{eqApproxHmnparSommeUi}, we prove the renewal results
Proposition \ref{propRenewalJoint} and Corollary \ref{propTRT}, using the same kind of techniques as in
Enriquez et al. \cite{ESZ3}, inspired by Feller \cite{Feller}.
We also show Proposition \ref{propTRT2} about the number of visited valleys with a similar method.
We prove before, in Lemma \ref {LemmaMajorationNt},
that with probability nearly $1$, the number $N_t$ of valleys visited by $X$ up to time $t$  is less than $n_t$,
which explains why we  only consider the first $n_t$ valleys.

We then turn to the proof of the localization, that is, our Theorem \ref{ththm}.
To this aim, using the previous renewal results, we prove that with probability nearly $1$, at time $t$,
the diffusion $X$  has already spent a quite large amount of time in the last valley visited,
that is, on $[\tilde L_{N_t}^- ,  \tilde L_{N_t}]$, and that $X(t)$ still belongs to this interval.
This allows us to prove that, knowing $N(t)=j$,
the quenched law of $X(t)$ is nearly equal to the invariant probability $\tilde \mu_j$
of a diffusion $Y_j$ in the potential
$\wk$, starting inside $[\tilde L_j^-,\tilde L_j]$ and reflected at $\tilde L_j^-$ and $\tilde L_j$.
This is a kind of convergence to the invariant probability measure, which we prove as in Brox \cite{Brox}, by
using a coupling between $Y_j$ starting from $Y_j(0)$ distributed as its invariant measure $\tilde \mu_j$, and
$X(.+H(\tilde m_j))$. Since $\tilde \mu_j$ is proportional to
$
    \exp(-[\wk(x)-\wk(\tilde m_j)) \un_{[\Lm,\tL]}(x)\dd x
$,
it is highly concentrated on a small neighborhood of $\tilde m_j$,
which leads to the localization of $X$ at time $t$ in this small neighborhood of $\tilde m_j$.

We then prove the Aging, that is, Proposition \ref{ThAging}.
For this, we apply the localization (Theorem \ref{ththm}), first at time $t$ with function $\phi$,
and second at time $\alpha t$ but with another function $\phi_\alpha$ defined by $\log(\alpha t)-\phi_\alpha(\alpha t) = \log t-\phi(t)$, so that the r.v. $m_i$ are the same in both cases.


Our hypothesis $0<\k<1$ is used many times in this paper.
We recall that the typical distance between two $h_t$-minima of $\wk$, or valleys of depth at least $h_t$, is asymptotically of order $e^{\k h_t}=t^{\k}e^{-\k \phi(t)}$. Moreover, the time spent by the diffusion $X$ in such valleys
is approximatively proportional to the exponential of the depth of such valleys.
So $X$ is trapped a long time, of order at least $t e^{-\phi(t)}$, and then is slowed by such valleys if (and only if) it hits some of them.\\
\indent
In particular,  the first $h_t$-minimum, or valley,  appears at a distance
of order $e^{\k h_t}=t^{\k}e^{-\k \phi(t)}$ which is much smaller than $t$ if $0<\k<1$, and
much larger than $t$ if $\k>1$.
Heuristically, if $0<\k<1$, there are many $h_t$-minima in $[0, t^{\k+\e}]$, $0<\e<1-\k$, which  trap and slow the diffusion $X$,
which explains why $X$ is zero-speed. However if $\k>1$, loosely speaking,
the first positive $h_t$-minimum is very far from the origin, at a distance much large than $t$, and then
at time $t$ the diffusion $X$ has not yet reached it; $X$ has then
not been trapped nor slowed by such deep valleys, and in this case $X$ has a positive speed.

In all the paper, $0<\kappa<1$ is fixed, and $C_+$ and $c_+$ (resp. $C_-$ and $c_-$) denote
positive constants that may increase (resp. decrease) from line to line and may only depend on our fixed constant $\kappa$. Moreover some events are denoted by $\B_{\text{i}}^{\text{j.k}}$ for some $i$, $j$, $k$; for example
$\B_3^{\ref{propIV}}$ is the event number $3$ introduced in the proof of Lemma \ref{propIV}.


\mysection{Standard valleys and path decomposition of the potential}\label{environments}




\subsection{$3$-dimensional drifted Bessel processes}\label{SubSectDefDriftedBessel}
In this subsection, we introduce $3$-dimensional drifted Bessel processes
as drifted Brownian motions conditioned to stay positive.
These processes are helpful to describe the law of the potential $\wk$ near the $h_t$-minima $m_i$
and then to estimate relevant quantities depending on this potential $\wk$, mainly with formulas
\eqref{eqLienEntreRetWzPositif} and \eqref{eqLienEntreRetWzNul}
below.


\noindent
For any process $(U(t),\ t\geq 0)$ and any $a\in\R$, we denote the hitting time of $a$ by $U$ as
$$
\tau^{U}(a)
:=
\inf\{t>0,\  U_{}(t)=a\},
$$
with the convention $\inf\emptyset=+\infty$.
We denote by $(\lo_U(t,x),\ t\geq 0,\ x\in\R)$ the bicontinuous version of the local time of $U$ when it exists, which is the case for $X$ and for Brownian motions.
We also denote by $U^a$ the process $U$ starting from $a$, with the notation $U=U^0$.
We sometimes write $P^a(U\in .):=P(U^a\in.)$.
In particular, for $x\in\R$, $\zeta\neq 0$, $W_\zeta^x$ is a $(-\zeta/2)$-drifted Brownian motion starting from $x$.

Let $\zeta\neq 0$.
We recall the definition of a
{\it $(-\zeta/2)$-drifted Brownian motion $W_{\zeta}$ Doob-conditioned to stay positive}
(see \cite{Faggionato}, 5. p. 1783, or \cite{Bertoin}, Chapter VII.3 and references therein for more details).
We consider the $\sigma$-fields $\mathcal{F}_t$ defined on $C([0,\infty),\R)$ by
$\mathcal{F}_t:=\sigma(Y(s),\ 0\leq s \leq t)$, $t\geq 0$, and $\mathcal{F}_\infty:=\sigma(Y(s),\ s\geq 0)$,
for a generic element $(Y(s), s\geq 0)$ of the path space  $C([0,\infty),\R)$.
Following (\cite{Faggionato}, p. 1783), for $z>0$, the probability measure $P_z^{\zeta/2, \uparrow}$ is defined on $C([0,+\infty), \R)$ by
\begin{equation}\label{EqDefLoiDoobContidionnement}
    P_z^{\zeta/2, \uparrow}(\Lambda)
:=
    \frac{1}{1-e^{\zeta z}}
    E\Big[\big[1-\exp(\zeta W_{\zeta}^z(t))\big],\, W_{\zeta}^z\in\Lambda,\, t<\tau^{W_{\zeta}^z}(0)\Big],
\qquad
    \Lambda\in\mathcal{F}_t,\ t\geq 0.
\end{equation}
This induces a unique probability measure
$P_z^{\zeta/2, \uparrow}$ on $(C([0,\infty),\R), \mathcal{F}_\infty)$.
Moreover, $P_z^{\zeta/2, \uparrow}$ converges weakly as $z\to 0^+$, in the space of Skorokhod $D([0,\infty),\R)$
(see \cite{Bertoin} VII.3 Prop. 14 and comments below)
and in $C([0,\infty),\R)$ (see \cite{Faggionato} p. 1784) to a probability measure
on $C([0,\infty),\R)$ denoted by $P_0^{\zeta/2, \uparrow}$.
The canonical process, which we denote by $(R(s),\ s\geq 0)$,
takes values in $\R_+$. It
is a Feller process for the family $(P_z^{\zeta/2, \uparrow}, \ z\geq 0)$,
and then is strong Markov.
Its infinitesimal generator is given for every $x>0$ by
\begin{eqnarray}
    \label{infinitesimalGenerator}
    \frac{1}{2}\frac{d^2}{dx^2}+\frac{\zeta}{2} \coth\left(\frac{\zeta}{2} x\right)\frac{d}{dx}.
\end{eqnarray}
This infinitesimal generator is given by  Lemma 6 of \cite{Faggionato}, which is true in the case of positive or negative drift $\zeta/2$.

This process $R$ can be thought of as a
$(-\zeta/2)$-drifted Brownian motion $W_{\zeta}$ Doob-conditioned to stay positive, with the terminology of
\cite{Bertoin}, which is called Doob conditioned to reach $+\infty$ before $0$ in \cite{Faggionato}.
We notice in particular that, by \eqref{infinitesimalGenerator}, or by (\cite{Faggionato}, eq. (5.4)),
or directly by \eqref{EqDefLoiDoobContidionnement} combined with Girsanov theorem,  the law of $R$ is the same if $\zeta$ is replaced by $-\zeta$.
That is, $W_\zeta$ Doob-conditioned to stay positive has the same law as $W_{-\zeta}$ Doob-conditioned to stay positive.
This is the case in particular in (\cite{Faggionato}, Thm. 2) and then also in (\cite{Faggionato}, eq. (1.1)), where the sign should be a $+$ in every case.

This process $R$ is also shown in Rogers et al.  (\cite{RogPit}, Thm. 3 and eq. (13)) to be equal in law to the euclidian norm of a $3$-dimensional drifted
Brownian motion, with drift $\zeta/2$ in some direction given by a unit vector of $\R^3$.
We do not use this result in this paper, but for this reason,
in the rest of the paper, we call for $z\geq 0$ the process $R$ with law $P_z^{\zeta/2, \uparrow}$  a
$3$-dimensional $|\zeta/2|$-drifted Bessel process
starting from $z$.
As in \cite{RogPit}, its law is denoted by $\text{BES}^z(3, |\zeta/2|)$,
and by $\text{BES}(3, |\zeta/2|)=\text{BES}^0(3, |\zeta/2|)$ if it starts from $z=0$.

In the rest of the paper, it is often useful to consider a process
$(R(s),\ s\geq 0)$ with law $\textnormal{BES}^z(3,\k/2)$ for some $z\geq 0$. We have by the previous remark,
when $R$ starts from $z>0$,
\begin{equation}\label{eqLienEntreRetWzPositif}
    P^z(R\in\Lambda)
=
    P_z^{\kappa/2, \uparrow}(\Lambda)
=
    P_z^{-\kappa/2, \uparrow}(\Lambda)
=
    P(W_{-\k}^z\in\Lambda|\tau^{W_{-\k}^z}(0)=\infty),
\qquad
    z>0,\ \Lambda\in\mathcal{F}_\infty,
\end{equation}
where the last equality is noticed in \cite{Faggionato} just before its Lemma 6
for $\Lambda\in \mathcal{F}_t$, $t\geq 0$ since $-\kappa/2<0$ and then $W_{-\k}$ has a positive drift $\k/2$,
and so this is true for all $\Lambda\in \mathcal{F}_\infty$.
As a consequence, when $R$ starts from $0$, that is, when the law of $(R(s), \ s\geq 0)$ is $\textnormal{BES}(3,\k/2)$,
we have for all $\Lambda\in \mathcal{F}_\infty$ such that $P(R\in \partial\Lambda)\neq 0$,
\begin{equation}\label{eqLienEntreRetWzNul}
    P(R\in \Lambda)
=
    P_0^{\kappa/2, \uparrow}(\Lambda)
=
    P_0^{-\kappa/2, \uparrow}(\Lambda)
=
    \lim_{z\downarrow 0}
    P_z^{-\kappa/2, \uparrow}(\Lambda)
=
    \lim_{z\downarrow 0} P(W_{-\k}^z\in\Lambda|\tau^{W_{-\k}^z}(0)=\infty).
\end{equation}



\subsection{Path decomposition of the potential $\wk$ in the neighborhood of the $h_t$-minima $\mf_i$ }
\label{SubSectWilliamsDecomposition}
The point of view of $h$-extrema has been used recently in some studies of random walks or diffusions in random environment in the recurrent case, see e.g. Cheliotis \cite{Cheliotis_Favorite},
Bovier et al. \cite{Bovier_Faggionato} and  Devulder \cite{Devulder_Persistence}.
We now recall some results for $h$-extrema of $\wk$.
Let
$$
V^{(i)}(x):=W_{\kappa}(x)-W_{\kappa}(\mf_i), \qquad x\in\R,\  i\in\N^*,
$$
which is the potential $\wk$ translated so that it is $0$ at the local minimum $\mf_i$. We also define
\begin{equation*}
    \tau_i^-(h) := \sup \{s < \mf_i,\  V^{(i)}(x)=h\},
\qquad
    \tau_i(h):=  \inf \{s > \mf_i,\  V^{(i)}(x)=h\},
\qquad
     h>0.
\end{equation*}

The following result has been proved by Faggionato \cite{Faggionato}
for {\bf (i)} and {\bf (ii)}.

\begin{fact}\label{Fact_Williams} (path decomposition  of $\wk$ around the $h_t$-minima $m_i$)\\
{\bf (i)} The truncated trajectories
$\big(V^{(i)}(\mf_i-s),\ 0 \leq s \leq \mf_i-\tau_i^-(h_t)\big)$,
$\big(V^{(i)}(\mf_i+s),\  0\leq s \leq \tau_i(h_t)-\mf_i\big)$, $i\geq 1$
are independent.
\\
{\bf (ii)} Let $(R(s),\ s\geq 0)$ be a process with law $BES(3,\k/2)$.
All the truncated trajectories
$\big(V^{(i)}(\mf_i-s),\ 0 \leq s \leq \mf_i-\tau_i^-(h_t)\big)$ for $i\geq 2$
and
$\big(V^{(j)}(\mf_j+s),\  0\leq s \leq \tau_j(h_t)-\mf_j\big)$ for $j\geq 1$
are equal in law to
$\big(R(s),\ 0 \leq s \leq \tau^R(h_t) \big)$.
\\
{\bf (iii)} For $i\geq 1$, the truncated trajectory
$\big(V^{(i)}(s+\tau_i(h_t)), \ s \geq 0\big)$
is independent of $\big(\wk(s),\ s\leq \tau_i(h_t)\big)$
and is equal in law to $\big(\wk^{h_t}(s),\ s\geq 0\big)$, that is, to a $(-\k/2)$-drifted Brownian motion starting from $h_t$.
\end{fact}


We point out that for reasons linked to renewal theory, the first trajectory in {\bf (ii)} for $i=1$ has a different law, which we will not use in this paper.

\noindent {\bf Proof:}
Notice that $M_{i-1}<\tau_i^-(h_t)<m_i<\tau_i(h_t)<M_i$, $i\in\Z$.
Moreover the $h_t$-extrema of $\wk$ are the r.v. $m_i$ and $M_i$, $i\in\Z$.
So  {\bf(i)} follows from the independence of
the truncated trajectories between consecutive $h_t$-extrema
proved by
Faggionato (\cite{Faggionato}, Theorem 1; notice in particular the comment about independence just before its equation (2.26)).
Result {\bf (ii)} is proved in Faggionato (\cite{Faggionato}, Theorem 2 p. 1785),
since, as explained in the paragraph after \eqref{infinitesimalGenerator} of the present paper,
a Brownian motion starting at $0$ with drift $\k/2$ Doob conditioned to reach $+\infty$ before $0$
has the same law as a Brownian motion starting at $0$ with drift $-\k/2$ Doob conditioned to reach $+\infty$ before $0$, and the same law as $BES(3,\k/2)$.

Finally, let $i\geq 1$. For every $x\geq 0$,
$\tau_i(h_t)\leq x$ if and only if  the function
$s\mapsto \wk(s)\un_{(-\infty, x]}(s)+\wk(x)\un_{(x,+\infty)}(s)$
has at least $i$ $h_t$-minima in $(0,x]$.
Consequently,
$\tau_i(h_t)$ is a stopping time for the $\sigma$-field
$\sigma(\wk(s), s\in(-\infty, x])$, $x\geq0$.
Hence the strong Markov property gives {\bf (iii)}.
\hfill$\Box$

We now introduce
\begin{align}
    \tau_1^*(h)
:=
    \inf\{u\geq 0,\ \wk(u)-\inf_{[0,u]}\wk\geq h\},
\qquad
    h\geq 0.
\label{tauetoile}
\end{align}
We gather here some other results proved by Faggionato \cite{Faggionato}
that will be useful in the following:

\begin{fact}\label{Fact_Faggio} (Faggionato \cite{Faggionato}).
The random variables $\wk(M_i)-\wk(m_{i+1})-h_t$ and $\wk(M_i)-\wk(m_i)-h_t$,
$i\geq 1$, are
independent. They are exponential variables of mean respectively
$
    (2/\k)\sinh(\k h_t/2)e^{\k h_t/2}
$
and
$
    (2/\k)\sinh(\k h_t/2)e^{-\k h_t/2}
$.
Moreover for $i\geq 1$,
\begin{align}\label{EqLaplaceEll-}
    E \big(e^{- \alpha (m_{i+1}-M_i)} \big)
& =
    \frac{e^{-\kappa h_t/2 }\sqrt{2\alpha+ \kappa^2/4}}
    {\sqrt{2\alpha+ \kappa^2/4}\cosh( h_t\sqrt{2\alpha+ \kappa^2/4})
    -(\kappa/2) \sinh( h_t\sqrt{2\alpha+ \kappa^2/4})},
\quad \a>0,
\\
        E \left(e^{-\a (\tau_i(h)-m_i)} \right)
& =
    \frac{2\sqrt{2\a +\k^2/4}\sinh(\k h/2)}{\k \sinh(h\sqrt{2\a +\k^2/4})},
\qquad
    0<h\leq h_t,
\qquad
    \alpha>-\k^2/8,
\label{eqFaggioLaplaceHitting}
\end{align}
\begin{equation}\label{eqProbaM0avantm1Fact}
    P(0< \Mf_0 <\mf_{1})
=
    P(0\in[m_0,M_0))
\sim_{t\to+\infty}
    \k h_t e^{-\k h_t}.
\end{equation}
Also, $-\inf_{[0,\tau_1^*(h)]}\wk$ is, for $h>0$,
exponentially distributed
with mean $2\k^{-1} \sinh(\k h/2)e^{\k h/2}$.
So for large $h$ for all $x>0$,
\begin{equation}\label{InegBeta*}
    P\big(\inf\nolimits_{[0,\tau_1^*(h)]}\wk \geq -x\big)
\leq
    e^{-\k h} x.
\end{equation}
Finally,
\begin{equation}\label{InegProbaTau1Etolie}
    E\big(\tau_1^*(h)\big)
\leq
    C_+ e^{\kappa h},\qquad h>0,
\end{equation}
\end{fact}

\noindent{\bf Proof:}
All this is proved in Faggionato \cite{Faggionato} with $\mu=\k/2$, so we just explain where these results are stated in this paper
\cite{Faggionato}.

Thanks to (\cite{Faggionato} Thm. 1 and the remark before its (2.26)),
$(\wk(\Mf_i+x)-\wk(\Mf_i),\ 0\leq x\leq \mf_{i+1}-\Mf_i)$
and
$(\wk(m_i+x)-\wk(m_i),\ 0\leq x\leq M_i-m_i)$, $i\geq 1$,
are independent, and their laws are respectively $P_{-}^{\k/2}$ and $P_{+}^{\k/2}$,
which are defined in (\cite{Faggionato} p. 1769), with $h=h_t$ in our case.
In particular, the lengths $m_{i+1}-M_i$ and $M_i-m_i$
of these truncated and translated trajectories, called $h_t$-slopes in \cite{Faggionato},
and their excess heights $\wk(M_i)-\wk(m_{i+1})-h_t$ and $\wk(M_i)-\wk(m_i)-h_t$ are
denoted respectively by $\ell_-$, $\ell_+$, $\zeta_-$ and $\zeta_+$
in (\cite{Faggionato}, Prop. 1, formulas (2.11), and (2.14) with $\lambda=0$, see also (2.1)),  which provides their
law, and in particular this proves the present fact up to \eqref{EqLaplaceEll-}.

For $0<h\leq h_t$, $m_i$ is a $h$-minimum, so for $i\geq 1$, by (\cite{Faggionato}, Thm. 1),
$(\wk(x+m_i)-\wk(m_i),\ 0\leq x\leq \tau_i(h)-m_i)$
has the law defined as $P_+^{\k/2}$ for this $h$ in (\cite{Faggionato}, eq. (2.9)).
Hence $\tau_i(h)-m_i$ has the same law as the r.v. called $\tau-\sigma$ in (\cite{Faggionato} eq. (2.9), see also (2.2)).
Its Laplace transform is given in (\cite{Faggionato} eq. (2.3) of Lem. 1). This gives \eqref{eqFaggioLaplaceHitting}.

Now, $0< \Mf_0 <\mf_{1}$ if and only if $0\in[m_0,M_0)$, that is,
the translated trajectory between the two consecutive $h_t$-extrema surrounding $0$
belongs to the set of slopes starting at an $h_t$-minimum and ending at an $h_t$-maximum,
denoted by $\mathcal{W}_+$ in (\cite{Faggionato} after eq. (2.10)). The probability
of this event is given in
(\cite{Faggionato}, Thm. 1 and eq. (2.25)  in the case $\mathcal{W}_+$, $h=h_t$),
which leads to \eqref{eqProbaM0avantm1Fact}.

We turn to \eqref{InegProbaTau1Etolie}.
For $h>0$, $\tau_1^*(h)$ is denoted by $\tau$
in (\cite{Faggionato}, eq (2.2)).
Let $\ell_-$ and $\ell_+$ be as in (\cite{Faggionato}, Prop. 1),
that is,
$\ell_+\egloi \tau-\sigma+\sigma'\geq \tau-\sigma$ by (\cite{Faggionato}, eq. (2.9))
and $\ell_-\egloi \tau'-\sigma'+\sigma\geq \sigma$ by (\cite{Faggionato}, eq. (2.10)),
where the values of $\sigma'\geq 0$ and $\tau'-\sigma'\geq 0$ are not important here. Applying
(\cite{Faggionato}, eq. (2.17)) gives
$
    E\big(\tau_1^*(h)\big)
=
    E(\tau)
=
    E[(\tau-\sigma)+\sigma]
\leq
    E(\ell_-+\ell_+)
\leq
    2\k^{-2}e^{\k h}
$.

Finally,
$\inf_{[0,\tau_1^*(h)]}\wk$ is denoted by $\beta$ in (\cite{Faggionato}, eq. (2.2));
its law is given by (\cite{Faggionato}, eq. (2.4) in its Lem. 1), which leads to our \eqref{InegBeta*}.
\hfill $\Box$

\subsection{Definition of valleys, and standard $h_t$-minima \label{sec2.2} $\tm$, $i\in\N^*$ }~\\
We would like to consider some valleys of the potential around the $h_t$-minima $m_i$, $i\in\N^*$,
which would be intervals containing at least $[\tau_i^-((1+\k)h_t), M_i]$. However, these valleys could intersect.
In order to define valleys which are well separated and i.i.d., we introduce the following notation.  This notation is later used to define valleys of the potential around some $\tilde m_i$, which are shown in Lemma \ref{CVs} to be equal to the $m_i$ for $1\leq i \leq n_t$ with large probability.

\begin{figure}[h]
\begin{center}
\scalebox{1.5}{\input{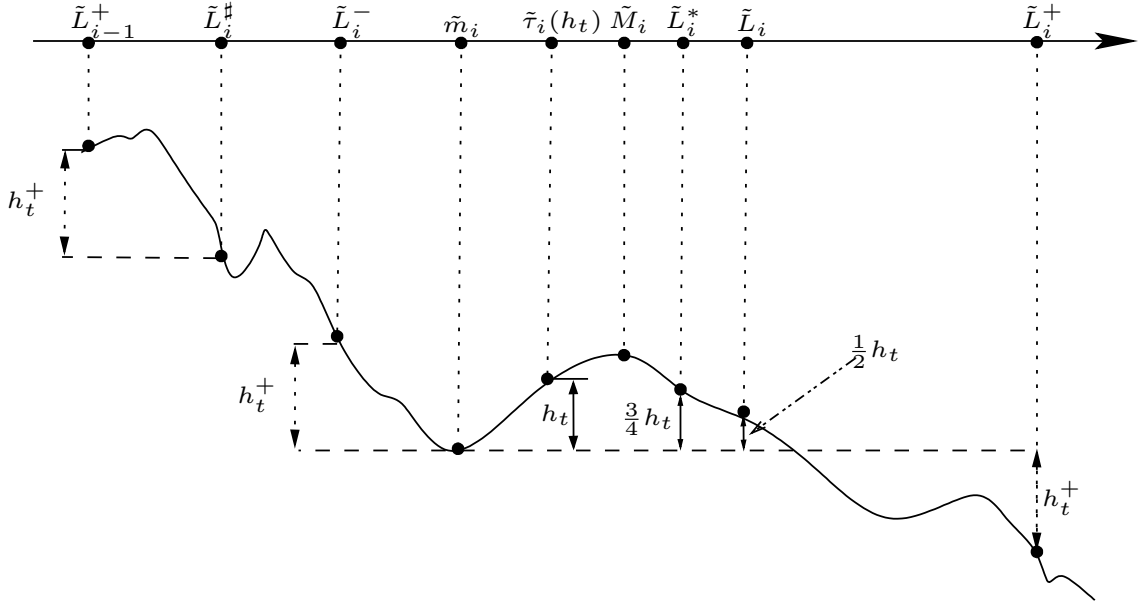}}
\caption{Schema of the potential between $\tilde L_{i-1}^+$ and $\tilde L_i^+$, in the case $\tilde L_i^\sharp <\tilde L_i^-$}\label{fig1}
\end{center}
\end{figure}

\noindent Let
$$
    h_t^+
:=
    (1+ \kappa+ 2\delta) h_t.
$$
We define  $\tilde L_0^+:=0$, $\tilde m_0:=0$,
and recursively for $i \geq 1$ (see Figure \ref{fig1}),
\begin{eqnarray}
     \tilde L_i^{\sharp}
&:= &
    \inf\{x>\tilde L_{i-1}^+,\ \wk(x)\leq\wk(\tilde L_{i-1}^+)-h_t^+   \},
\nonumber\\
    \tilde \tau_i(h_t)
& := &
    \inf \big\{x \geq   \tilde L_i^{\sharp},\ \wk(x)-\inf\nolimits_{[\tilde L_i^{\sharp},x]}\wk \geq h_t\big\},
\label{eqDefTaui1}
\\
    \tm
&:= &
    \inf\big\{x \geq \tilde L_i^{\sharp},\ \wk(x)=\inf\nolimits_{[\tilde L_i^{\sharp},\tilde \tau_i(h_t)]}\wk\big\},
\nonumber\\
    \tilde L_i^{+}
&:= &
    \inf\{x>\tilde \tau_i(h_t),\ \wk(x)\leq\wk(\tilde \tau_i(h_t))-h_t-h_t^+   \}.
\nonumber
\end{eqnarray}
\smallskip
We also introduce the following random variables for $i\in\N^*$:
\begin{eqnarray}
    \tilde M_i
&:=&
    \inf\{ s>\tm,\ \wk(s)=\max\nolimits_{\tm \leq u\leq \tilde L_i^+} \wk(u)  \},
\nonumber\\
    \tilde L_i
& := &
    \inf\{x>\tilde \tau_i(h_t),\ \wk(x)-\wk(\tm)= h_t/2  \},
\label{eqDefLiTilde}
\\
    \tilde \tau_i(h)
& := &
    \inf \{s > \tilde m_i,\  \wk(x)-\wk(\tm)=h\}, \qquad h>0,
\label{eqDefTaui2}
\\
    \tilde \tau_i^-(h)
& := &
    \sup \{s < \tilde m_i,\  \wk(x)-\wk(\tm)=h\},
    \qquad h>0,
\nonumber\\
    \tilde L_i^{-}
&:= &
    \tilde \tau_i^-(h_t^+). 
\nonumber
\end{eqnarray}

We stress that these r.v. depend on $t$, which we do not write as a subscript to simplify the notation.
Notice also that $\tilde \tau_i(h_t)$ is the same in definitions \eqref{eqDefTaui1} and \eqref{eqDefTaui2} with $h=h_t$.
Moreover by continuity of $\wk$, $\wk(\tilde \tau_i(h_t))=\wk(\tilde m_i)+h_t$.
So, the $\tilde m_i$, $i\in\N^*$, are  $h_t$-minima,
since $\wk(\tilde m_i)=\inf_{[\tilde L_{i-1}^+,\tilde \tau_i(h_t)]}\wk$,
$\wk(\tilde \tau_i(h_t))=\wk(\tilde m_i)+h_t$ and
$\wk(\tilde L_{i-1}^+)\geq\wk(\tilde m_i)+h_t$.
Moreover,
\begin{equation}\label{InegRvTilde1}
\tilde L_{i-1}^+<\tilde L_i^{\sharp}\leq \tm < \tilde \tau_i(h_t)<\tilde L_i <\tilde L_i^{+},
\qquad
    i\in\N^*,
\end{equation}
\begin{equation}\label{InegRvTilde2}
\tilde L_{i-1}^+\leq \tilde L_i^{-}< \tm < \tilde \tau_i(h_t)<\tilde M_i<\tilde L_i^{+},
\qquad
    i\in\N^*.
\end{equation}
Furthermore by induction the r.v. $\tilde L_i^{\sharp}$, $\tilde \tau_i(h_t)$ and $\tilde L_i^{+}$, $i\in\N^*$ are stopping times
for the natural filtration of $(\wk(x), \ x\geq 0)$,
and so $\tilde L_i$, $i\in\N^*$, are also stopping times.
Also by induction,
\begin{equation}\label{InegLiPremieresDescentes}
    \wk(\tilde L_i^{\sharp})
=
    \inf_{[0, \tilde L_i^{\sharp}]} \wk,
\qquad
    \wk(\tilde m_i)
=
    \inf_{[0,\tilde \tau_i(h_t)]}\wk,
\qquad
    \wk(\tilde L_i^+)
=
    \inf_{[0, \tilde L_i^+]} \wk
=
    \wk(\tilde m_i)-h_t^+,
\end{equation}
for $i\in\N^*$.
We also introduce the analogue of $V^{(i)}$ for $\tm$ as follows:
$$
    \tV(x)
:=
    W_{\kappa}(x)-W_{\kappa}(\tilde m_i),
\qquad
    x\in\R,
    \ i\in\N^*.
$$

\noindent
We call $i\ th$ \textit{valley} the translated truncated potential
$\big(\tV(x),\ \tilde L_i^- \leq x \leq \tilde L_i \big)$, for $i\geq 1$.


We show in the following lemma that, with an overwhelming high probability, the first $n_t+1$ positive
$h_t$-minima $\mf_i$, $1 \leq i\leq n_t+1$,
coincide with the r.v. $\tm$, $1 \leq i\leq n_t+1$.
We introduce the corresponding event
$\mV_t:=\cap_{i=1}^{n_t+1}\{\mf_i=\tm\}$.

\begin{lemma} \label{CVs}
Assume $0<\delta<1$.
For $t$ large enough,
$$
    P\left(\overline{\mathcal V}_t\right)
\leq
    C_+ n_t e^{-\k h_t/2}
=
    e^{[-\k /2+o(1)]h_t}.
$$
Moreover, the sequence $\left(\big(\tV(x+\tilde {L}^+_{i-1}),\ 0\leq x \leq \tilde {L}^+_i-\tilde {L}^+_{i-1} \big),\ i \geq 1 \right)$, is i.i.d.
\end{lemma}

We notice that consequently, the valleys $\big(\tV(x),\ \tilde L_i^- \leq x \leq \tilde L_i \big)$,  $i\geq 1$,
are also i.i.d.

\noindent{\bf Proof:}
Since  $\tilde m_i$ is an $h_t$-minimum for $\wk$ for every $i\geq 1$,
we have
$\{\tilde m_i,\ i\in\N^*\}\subset \{\mf_i,i\in\N^*\}$.
We now assume that we are on the complement $\overline{ \mV}_t$ of ${\mV}_t$.
So the smallest $j\geq 1$ such that $m_j\neq \tilde m_j$ satisfies $1\leq j\leq n_t+1$.
Due to the previous inclusion and since $\tilde m_0=0$, we have
$\tilde m_{j-1}<\mf_j<\tilde m_{j}$.
If for this $j$, $\tilde L_{j}^{\sharp}<\mf_j<\tilde m_{j}$,
there would be a $v>\mf_j>\tilde L_{j}^{\sharp}$ such that
\begin{eqnarray}
   \wk(\mf_j)
& = &
    \inf_{[\mf_j,v]}\wk
\geq
    \inf_{[\tilde L_{j}^{\sharp},v]}\wk,
\nonumber\\
\wk(v)
& \geq &
    \wk(\mf_j)+h_t
\geq
    \inf_{[\tilde L_{j}^{\sharp},v]}\wk+h_t,
\label{eqLemmaVt11}
\end{eqnarray}
since $m_j$ is an $h_t$-minimum.
The definition of $\tilde \tau_{j}(h_t)$ and \eqref{eqLemmaVt11} would give $\tilde \tau_{j}(h_t)\leq v$,
and then
$\tilde L_{j}^{\sharp}<\mf_j<\tilde m_j\leq \tilde \tau_{j}(h_t)\leq v$, by definition of $\tilde m_j$.
Then
$
    \wk(\mf_j)
=
    \inf_{[\mf_j,v]}\wk
\leq
    \wk(\tilde m_{j})
$,
which contradicts the definition of $\tilde m_{j}$.

Hence, $\tilde m_{j-1}<\mf_j\leq \tilde L_{j}^\sharp$.
Thus by \eqref{InegLiPremieresDescentes},
$
    \wk(\mf_j)
\geq
    \wk\big(\tilde L_j^{\sharp}\big)
=
    \wk\big(\tilde L_{j-1}^+\big)-h_t^+
=
    \wk(\tilde m_{j-1})-2 h_t^+
$
if $j\geq 2$,
whereas
$
    \wk(\mf_j)
\geq
    \wk\big(\tilde L_1^{\sharp}\big)
=
    - h_t^+
$ if $j=1$. Consequently,
\begin{equation}\label{eqDecompositionVbarre}
     \overline{ \mV}_t
\subset
    \{\wk(\mf_{1})\geq -h_t^+\}
    \cup
    \cup_{j=2}^{n_t+1}\{\wk(\mf_j)\geq \wk(\mf_{j-1})-2h_t^+\}.
\end{equation}
For $j\geq 2$, we have by the start of Fact \ref{Fact_Faggio},
\begin{eqnarray}
&&
    P[\wk(\mf_j)\geq \wk(\mf_{j-1})-2h_t^+]
\label{eqProbaVbarreDecomposition1}
\\
& \leq &
    P\big[\wk(M_{j-1})-\wk(m_{j})\leq e^{\k h_t/2}\big]
    +
    P\big[\wk(M_{j-1})-\wk(m_{j-1})>e^{\k h_t/2}-2h_t^+\big]
\leq
    \frac{C_+}{e^{\k h_t/2}}.
\nonumber
\end{eqnarray}
For $j=1$, notice that either
$0< \Mf_0 <\mf_{1}$,
which has probability
\begin{equation}\label{eqProbaM0avantm1}
    P(0< \Mf_0 <\mf_{1})
\leq
    2h_t e^{-\k h_t}
\end{equation}
for large $t$ by \eqref{eqProbaM0avantm1Fact},
either
there is no $h_t$-maximum in $(0,\mf_{1}]$, and so $M_0\leq 0$.

In this last case we show that
$\mf_{1}\leq \tau_1^*(h_t)$, which we defined in  \eqref{tauetoile}.
Indeed, assume $\tau_1^*(h_t)<m_1$.
We consider $u=\inf\big\{x\geq 0,\ \wk(x)=\inf_{[0, \tau_1^*(h_t)]} \wk\big\}$
and $y=\inf\{x\geq 0,\ \wk(x)=\sup_{[\tau_1^*(h_t), m_1]} \wk\}$.
It follows from the definition of $h_t$-extrema
that $\wk(m_1)=\inf_{[M_0, m_1]}\wk$.
Since $M_0\leq 0\leq u <\tau_1^*(h_t)<m_1$,
this would give $\wk(m_1)\leq \wk(u)=\wk[\tau_1^*(h_t)]-h_t\leq \wk(y)-h_t$,
and $\wk(y)=\sup_{[u,m_1]}\wk$,
so $y$ is an $h_t$-maximum and belongs to $(0,m_1]$,
which contradicts this case.
So $\mf_{1}\leq \tau_1^*(h_t)$.


Hence in this case, there exists $v>m_1$ such that $\wk(v)\geq \wk(m_1)+h_t$ and $\wk(m_1)=\inf_{[m_1, v]}\wk$,
so $\tau_1^*(h_t)\leq v$. Thus, $\wk(m_1)=\inf_{[m_1, \tau_1^*(h_t)]}\wk$.
This, $M_0\leq 0$ and $\wk(m_1)=\inf_{[M_0, m_1]}\wk$ yield
$\wk(\mf_{1})=\inf_{[0,\tau_1^*(h_t)]}\wk$.
So, \eqref{InegBeta*} and \eqref{eqProbaM0avantm1} give for large $t$,
$$
    P[\wk(\mf_{1})\geq -h_t^+]
\leq
    2h_t e^{-\k h_t}
    +P\big[\inf\nolimits_{[0,\tau_1^*(h_t)]}\wk\geq - h_t^+\big]
\leq
    C_+ e^{-\k h_t/2}.
$$
Hence, using \eqref{eqDecompositionVbarre} and \eqref{eqProbaVbarreDecomposition1}, we get $P(\overline{{ \mV}_t})\leq  C_+  n_t e^{-\k h_t/2}$
for large $t$.


Finally, the fact that the sequence
$
\big(
    \tV(x+\tilde {L}^+_{i-1}),\ 0\leq x \leq \tilde {L}^+_i-\tilde {L}^+_{i-1}
\big)
$,
 $i \geq 1$, is i.i.d. follows from the strong Markov property applied at times $\tilde {L}^+_{i-1}$,
which are stopping times.
\hfill$\Box$

The following remark will be useful in the rest of the paper.

\begin{rem} \label{RemEgaliteAvecouSansTilde}
On ${ \mV}_t$, we have for every $1\leq i\leq n_t$, $m_i=\tilde m_i$,
and as a consequence,
$\tV(x)=V^{(i)}(x)$, $x\in\R$,
$\tau_i^-(h)=\tilde \tau_i^-(h)$ and
$\tau_i(h)=\tilde \tau_i(h)$ for  $h>0$.
Moreover, $\tilde M_i=M_i$.
Indeed, $\tilde M_i$ is an $h_t$-maximum for $\wk$, which belongs to $[\tilde m_i,\tilde m_{i+1}]=[m_i, m_{i+1}]$ on
${ \mV}_t$, and there is exactly one $h_t$-maximum in this interval since the $h_t$-maxima and minima alternate,
which we defined as $M_i$, so $\tilde M_i=M_i$. So in the following, on ${ \mV}_t$, we can write these r.v. with or without tilde.
\end{rem}


\subsection{Some technical estimates related to the potential}

\smallskip

We first provide  estimates for the distance between some points of a given valley:

\smallskip

\begin{lemma}  \label{lemMi}
Assume $0<\delta<1/2$. For large $t$, for all $1 \leq i \leq n_t$,
\begin{eqnarray}
    P\big[\mt_{i+1}-\Mt_{i} \leq e^{\k(1-\delta) h_t}, \mV_t\big]
& \leq &
    P\big(\mf_{i+1}-\Mf_{i} \leq e^{\k(1-\delta) h_t}\big)
\leq
    C_+ e^{-\k \delta h_t},
\label{2.0}
\\
    P\big[\tt_{i} (h)-\mt_{i} \geq  8h/\k\big]
& \leq &
    C_+ e^{- \kappa h /(2\sqrt{2})}, \qquad 0\leq h\leq h_t,
\label{bessel}
\\
    P\big[\mt_{i}-\tt_{i}^-(h) \geq  8h/\k\big]
& \leq &
    C_+ e^{- \kappa h /(2\sqrt{2})}, \qquad 0\leq h\leq h_t.
\label{eqLongueurTempsAtteinte-}
\end{eqnarray}
\end{lemma}

\smallskip

\noindent{\bf Proof of Lemma \ref{lemMi}:}
First, it follows from Remark \ref{RemEgaliteAvecouSansTilde}  that
$
    \{\mt_{i+1}-\Mt_{i} \leq e^{\k (1-\delta)h_t}\}\cap\mV_t
\subset
    \{\mf_{i+1}-\Mf_{i} \leq e^{\k (1-\delta)h_t}\}
$
for $1\leq i \leq n_t$.
This gives the first inequality of \eqref{2.0}, whereas  for the second one
we can use the results of Faggionato \cite{Faggionato} that we have gathered in
Facts \ref{Fact_Williams} and \ref{Fact_Faggio},
and then $\textnormal{BES}(3,\k/2)$ processes.
This technic is used several times in this proof and later in the paper.

\noindent
Let $i\geq 1$. We apply \eqref{EqLaplaceEll-}
with $\alpha(t)=e^{-\k (1-\delta)h_t}$ so that
$E(e^{-\a(t)(\mf_{i+1}-\Mf_i)})\sim_{t\to+\infty}\kappa^2 e^{-\k\delta h_t}/2$.
This and Markov inequality yield for large $t$,
$$
    P\big(\mf_{i+1}-\Mf_i\leq e^{\k (1-\delta)h_t}\big)
=
    P\big(e^{-\a(t) (m_{i+1}-M_i)}\geq e^{-1}\big)
\leq
    C_+ e^{-\k \delta h_t}.
$$
This ends the proof of \eqref{2.0}.

Applying \eqref{eqFaggioLaplaceHitting} with $\alpha=-\kappa^2/16$ gives
\begin{equation}\label{eqLaplaceTauih}
    E \left(e^{\kappa^2 (\tau_i(h)-\mf_i)/16} \right)
=
    \frac{\sinh(\k h /2)}{\sqrt{2}\sinh(\k  h/ \sqrt{8})}
\sim_{h\to+\infty}
    \frac{1}{\sqrt{2}}e^{\kappa (1- 1/ \sqrt{2})h /2}.
\end{equation}
So, choosing $C_+$ large enough, RHS of \eqref{eqLaplaceTauih} $\leq C_+ e^{\kappa (1- 1/ \sqrt{2})h /2}$
for all $0\leq h \leq h_t$.
Thus,
\begin{equation}\label{eqProbaLaplacetauih}
    P\big(\tau_i(h)-m_i\geq      8h/\k\big)
=
    P\big(e^{(\k^2/16)[\tau_i(h)-m_i]}\geq     e^{\k h/2}\big)
\leq
    C_+ e^{-\k h/(2\sqrt{2})},
\end{equation}
for all $0\leq h\leq h_t$ by Markov inequality.
Since
$\tilde \tau_i(h)-\tilde m_i=\tau_i(h)-m_i$ on $\mV_t$
by Remark \ref{RemEgaliteAvecouSansTilde},
and
$
    P\big(\overline{\mathcal V}_t\big)
\leq
    e^{-\k h_t/(2\sqrt{2})}
\leq
    e^{-\k h/(2\sqrt{2})}
$
for large $t$ for all $0\leq h \leq h_t$ by Lemma \ref{CVs},
$$
    P\big( \tilde \tau_i(h)-\tilde m_i\geq 8h/\k\big)
\leq
    P\big(\tau_i(h)-m_i\geq 8h/\k, \mV_t\big)
+
    P\big(\overline{\mathcal V}_t\big)
\leq
    C_+ e^{-\k h/(2\sqrt{2})},
\qquad
    0\leq h\leq h_t,
$$
which gives \eqref{bessel}.
According to Fact \ref{Fact_Williams} {\bf (ii)},
$m_i-\tau_i^-(h)$ and $\tau_i(h)-m_i$ are equal in law for $i\geq 2$ and $0\leq h\leq h_t$.
So we can replace $\tau_i(h)-m_i$ by $m_i-\tau_i^-(h)$ in \eqref{eqProbaLaplacetauih},
which similarly as before gives \eqref{eqLongueurTempsAtteinte-} for $i\geq 2$.
Finally, $\tilde m_1-\tilde \tau_1^-(h)$ has the same law as
$\tilde m_2-\tilde \tau_2^-(h)$ for $0\leq h\leq h_t$ by Lemma \ref{CVs},
so the case $i=1$ of \eqref{eqLongueurTempsAtteinte-} follows from the case $i=2$.
\hfill$\Box$

\smallskip

We also recall that a scale function of $W_{\zeta}$ $(x\mapsto W(x)-\zeta x/2)$
is for $\zeta\neq 0$ (see e.g. \cite{Faggionato}, (5.1)),
\begin{eqnarray}
    s_{\zeta}(u)
& := &
    e^{\zeta u}-1
=
    2e^{\zeta u /2}\sinh(\zeta u /2),
\qquad
    u\in\R,
\nonumber
\\
    P\big[\tau^{W_{\zeta}^y}(z)<\tau^{W_{\zeta}^y}(x)\big]
& = &
    s_{\zeta}(y-x)/s_{\zeta}(z-x),
\qquad
    x<y<z,
\label{eqScaleFunctionWk}
\end{eqnarray}
for which we remind that $W_{\zeta}^y$ denotes a $(-\zeta/2)$-drifted Brownian motion starting from $y$.

We also need some estimates on hitting times by $\wk$ and a process $R$ with law $\textnormal{BES}(3,\k/2)$
(recall that in \eqref{3.10}, $P^{\alpha h}$ denotes the law of $R$ starting from $\alpha h$,
that is, with law $\textnormal{BES}^{\alpha h}(3,\k/2)$):

\smallskip

\begin{lemma}
Let $0< \gamma< \alpha  <\omega $, and $R$ having law $\textnormal{BES}(3,\k/2)$.
For all $h$ large enough, 
\begin{eqnarray}
    P^{\alpha h}\left(\tau^{R}(\gamma  h)<\tau^{R} (\omega h)  \right)
& \leq  &
    2 \exp[-\kappa (\alpha-\gamma) h ]
    \label{3.10},
\\
    P\left(\tau^{R}( \omega h)-\tau^{R} (\alpha h) \leq 1  \right)
& \leq &
    4\exp\big[- [(\omega-\alpha) h]^2/3\big],
    \label{3.10b}
\\
    P\bigg(\int_0^{\tau^R(h)}e^{R(u)}\dd u \geq e^{(1-\alpha) h}\bigg)
& \geq  &
    1-3\exp[-\k \alpha h /2], \qquad 0<\alpha<1,
\label{MinorationAVallee}
\\
    P(\tau^R(h)>8h/\k )
& \leq &
    C_+ \exp\big[- \kappa h /(2\sqrt{2})\big].
    \label{bessel4}
\end{eqnarray}
\end{lemma}

\noindent{\bf Proof:}
We recall that $R$ has the same law as the $(\k/2)$-drifted Brownian motion
$W_{-\kappa}$ Doob-conditioned to stay positive.
In particular, this proof relies on \eqref{eqLienEntreRetWzPositif} and \eqref{eqLienEntreRetWzNul}.
First, using \eqref{eqLienEntreRetWzPositif} and \eqref{eqScaleFunctionWk},
we have for $0< \gamma< \alpha  <\omega $,
\begin{eqnarray}
  \text{LHS of } \eqref{3.10}
& = &
    P\big[\tau^{W_{-\k}^{\alpha h}}({\gamma h})<\tau^{W_{-\k}^{\alpha h}}({\omega h})
    \mid \tau^{W_{-\k}^{\alpha h}}(\infty)<\tau^{W_{-\k}^{\alpha h}}(0)\big] \nonumber
    \\
& = &
    P\big[\tau^{W_{-\k}^{\a h}}({\gamma h})<\tau^{W_{-\k}^{\a h}}({\omega h})]
    P
    \big[\tau^{W_{-\k}^{\gamma h}}(\infty)<\tau^{W_{-\k}^{\gamma h}}(0)\big]
    /P\big[\tau^{W_{-\k}^{\a h}}(\infty)<\tau^{W_{-\k}^{\a h}}(0)\big]  \nonumber
    \\
& = &
    \left(1-\frac{s_{-\k}[(\a-\gamma)h]}{s_{-\k}[(\o-\gamma)h]}\right)
    \frac{s_{-\k}(\gamma h)}{s_{-\k}(\alpha h)}  \nonumber
    \\
& = &
    \frac{\sinh[\k (\o-\a) h/2]\sinh(\k \gamma h/2)}
    {\sinh[\k (\o-\gamma) h/2]\sinh(\k \a h /2)}, \label{preuve3.2}
\end{eqnarray}
where LHS means left hand side. This gives \eqref{3.10} for large $h$.


Now, we notice that the left hand side of \eqref{3.10b} is thanks to \eqref{eqLienEntreRetWzNul}
for $0<\alpha<\omega$,
\begin{eqnarray*}
&&
    \lim_{z\downarrow 0}
    P\big[
    \tau^{W_{-\k}^z}(\omega h)-\tau^{W_{-\k}^z}(\alpha h)\leq 1
    \mid \tau^{W_{-\k}^z}(0)=\infty
    \big]
\\
& = &
    \lim_{z\downarrow 0}
    P\big[\tau^{W_{-\k}^z}({\a h})<\tau^{W_{-\k}^z}(0)\big]
    P\big[\tau^{W_{-\k}^{\alpha h}}({\o h})\leq 1, \tau^{W_{-\k}^{\alpha h}}(0)=\infty\big]
    /P\big[\tau^{W_{-\k}^z}(0)=\infty\big]
    ,
\end{eqnarray*}
where we applied the strong Markov property at time $\tau^{W_{-\k}^z}(\alpha h)$.
Moreover for large $h$,
\begin{align*}
   P\big[\tau^{ W_{-\kappa}^{\alpha h}}( \omega  h ) \leq 1\big]
=
    P\bigg(\sup_{x\in[0,1]}\big(W(x)+\alpha h+\k x/2\big)\geq \omega h\bigg)
& \leq
    P\bigg(\sup_{[0,1]}W\geq (\omega-\alpha) h-\frac{\k}{2}\bigg)
\\
& \leq
    2\exp\big[-[(\o-\a)h]^2/3\big],
\end{align*}
because
$\sup_{[0,1]} W \egloi |W(1)|$, where $\egloi$ denotes equality in law,
and $P(|W(1)|\geq x)\leq 2e^{-x^2/2}$ for $x\geq 1$.
Since
$\lim_{z\downarrow 0} P\big[\tau^{W_{-\k}^z}({\a h})<\tau^{W_{-\k}^z}(0)\big]/P\big[\tau^{W_{-\k}^z}(0)=\infty\big]=(1-e^{-\k\a h})^{-1}\leq 2$
for large $h$ by \eqref{eqScaleFunctionWk},
we get \eqref{3.10b}.

To prove \eqref{MinorationAVallee}, let $0<\alpha<1$. Notice that the probability of
$
\big\{\inf_{\tau^R[(1-\alpha/2)h]\leq u \leq \tau^R(h)}R(u)\geq (1-\alpha)h\big\}
\cap
\big\{ \tau^R(h)-\tau^R[(1-\alpha/2)h]\geq 1\big\}
$
is at least $1-3e^{-\k \alpha h /2}$ for large $h$ by \eqref{3.10} and \eqref{3.10b}.
Moreover, we have on this event,
$$
    \int_0^{\tau^R(h)}e^{R(u)}\dd u
\geq
    \int_{\tau^R[(1-\alpha/2)h]} ^{\tau^R(h)} e^{R(u)}\dd u
\geq
    \big[\tau^R(h) -\tau^R[(1-\alpha/2)h]\big] e^{(1-\alpha)h}
\geq
    e^{(1-\alpha)h},
$$
which proves \eqref{MinorationAVallee}.
Finally for $0\leq h\leq h_t$,
$\tau^R(h)$ is equal in law to $\tau_1(h)-m_1$ by Fact \ref{Fact_Williams} {\bf (ii)},
so \eqref{bessel4}  follows from \eqref{eqProbaLaplacetauih}
\hfill$\Box$

We also provide in the following lemma some probabilities concerning
$\big(\tilde V^{(i)}(x),\ \tilde \tau_i^-(h_t^+)\leq x \leq \tilde \tau_i^-(h_t)\big)$.
Unfortunately, they cannot be evaluated directly with the help of Fact \ref{Fact_Williams} and BES$(3,\k/2)$ processes,
so we evaluate them with more elementary technics. The proof of this lemma is deferred to Section \ref{SectionAnnexe}.

\begin{lemma}\label{LemmaProbasSansFaggionato}
With probability $P$ at least $1-e^{-\k h_t/8}$ for large $t$, we have for every $1\leq i\leq n_t$,
\begin{eqnarray}
\label{eqLemmaSansFaggio1}
    \tilde m_i-\tilde L_i^-
\leq
    \tilde L_i^+-\tilde L_i^-
=
    \tilde L_i^+-\tilde \tau_i^-(h_t^+)
& \leq &
    40 h_t^+/\k,
\\
    \int_{\tilde L_i^-}^{\tilde m_i} e^{\tilde V^{(i)}(u)}\dd u
& \geq &
    e^{h_t^+-\k h_t/2},
\label{eqLemmaSansFaggio2}
\\
    \inf_{[\tilde \tau_i^-(h_t^+), \tilde \tau_i^-(h_t)]}\tilde V^{(i)}
& \geq &
    h_t/2.
\label{eqLemmaSansFaggio3}
\end{eqnarray}
\end{lemma}

\noindent We end this section with the following basic result and its proof:
\smallskip

\begin{lemma}\label{lemmaProbaAtteintewk}
Let $0< \alpha  <\omega $.
We have for large $h$,
\begin{equation}
    P\big[\tau^{ W_{\kappa}}(- \alpha  h ) \geq 2 \omega h/ \kappa  \big]
\leq
    \exp\big[-\kappa  (\omega-\alpha)^2 h/(4 \omega)\big].
\label{3.14b}
\end{equation}
\end{lemma}

\noindent{\bf Proof:} We have,
$$
    P\big[\tau^{\wk}(-\alpha h) \geq 2 \omega h/ \kappa\big]
\leq
    P\big[W(2 \omega h/ \kappa) \geq (\omega- \alpha) h\big]
=
    P\big[W(1) \geq \sqrt{\kappa h}(\omega- \alpha)/ \sqrt {2 \omega} \big].
$$
Since $P(W(1)\geq x)\leq e^{-x^2/2}$ for $x\geq 1$, this proves the lemma.
\hfill$\Box$

\mysection{Quasi-Independence in the trajectories of $X$}\label{SectionQuasiIndep}

We now assume that $0<\k<1$.
In this section we provide some information on the typical trajectories of $X$. We also
show that the times spent in  the different valleys are, asymptotically in $t$,
nearly independent under the annealed probability $\P$ (see Proposition \ref{propind}).
Then we prove that the time spent by $X$ between the valleys is negligible.

We start with some classical formulae about the diffusion $X$, its hitting times and local times, which will be
important in the rest of the paper.

\subsection{Some formulas related to the diffusion $X$}\label{SubSectNotationX}
We first introduce
\begin{equation}
    \label{eqDefA}
    A(r)
:=
    \int_0^r e^{\wk(x)}\dd x,\qquad r\in\R,
\qquad
    A_\infty
:=
    \int_0^\infty e^{\wk(x)}\dd x<\infty \ a.s.
\end{equation}
We recall that $A$ is a scale function of $X$ under $\P^{\wk}$
(see e.g. \cite{Shi1} formula (2.2)), that is,
\begin{eqnarray}
    \P^{\wk}_y\big[H(z)<H(x)\big]
& = &
    [A(y)-A(x)]/[A(z)-A(x)],
\qquad
    x<y<z.
\label{eqScaleFunctionX}
\end{eqnarray}
Here $\P^{W_{\kappa}}_y$ denotes the law of the diffusion $X$ in the potential $\wk$, starting from $y$ instead of $0$, conditionally on $\wk$,
and $\E^{W_{\kappa}}_y$ is the corresponding expectancy.
As in (Brox \cite{Brox}, eq. (1.1) or Shi \cite{Shi}, eq. (2.2)), there exists a Brownian motion
$(B(s), \ s\geq 0)$, independent of $\wk$,
such that $X(t)=A^{-1}[B(T^{-1}(t))]$ for every $t\geq 0$, where
\begin{align}
    T(r)
:=
    \int_0^r \exp\{-2\wk[A^{-1}(B(s))]\}\text{d}s,
\qquad
    0\leq r< \tau^B(A_{\infty}) \label{T}.
\end{align}
\noindent As a consequence, with this notation, we have
(see e.g. Shi \cite{Shi1} eq. (4.5) and (4.6)),
\begin{equation}
    \label{eqFormuleH}
    H(r)
=
    T[\tau^B(A(r))]
=
    \int_{-\infty}^r \exp[-\wk(u)]\mathcal{L}_B[\tau^B(A(r)),A(u)]\dd u,
\qquad
     r\geq 0.
\end{equation}
Moreover the local time of $X$ is $\mathcal{L}_X(t,x)=e^{-\wk(x)}\mathcal{L}_B[T^{-1}(t),A(x)]$,
$t\geq 0$, $x\in\R$, as proved in Shi (\cite{Shi}, eq. (2.5)). So,
\begin{equation}\label{eqLocalTime}
    \mathcal{L}_X(H(r),x)
=
    e^{-\wk(x)}\mathcal{L}_B[\tau^B(A(r)),A(x)]
\qquad
    r\geq 0,\ x\in\R.
\end{equation}
It will sometimes be useful to notice that $H(r)= H_-(r)+H_+(r)$, where
\begin{equation}\label{eqDefH-H+}
    H_-(r)
:=
    \int_{-\infty}^0 \mathcal{L}_X(H(r),x) \dd x,
\qquad
    H_+(r)
:=
    \int_0^r \mathcal{L}_X(H(r),x) \dd x,
\qquad
    r\geq 0,
\end{equation}
are the time spent by $X$ respectively in $\R_-$ and in $\R_+$ before it first hits $r$.
We will sometimes need the following result (see Dufresne \cite{Dufresne}, or Borodin et al. \cite{BorodinSalminem} IV.48 p. 78):

\begin{fact}\label{FactDufresnes} (Dufresne)
The r.v. $2/A_\infty$ is a gamma variable of parameter $(\k,1)$,
and so has a density equal to  $e^{-x}x^{\k-1}\un_{\R_+}(x)/\Gamma(\k)$.

\end{fact}


\subsection{Probability that the diffusion $X$ leaves the valleys on the right}
In this subsection, we prove that for most environments,
with a large quenched probability,
the diffusion $X$, after first hitting $\tilde m_i$, leaves the valley $[\tilde L_i^-, \tilde L_i]$ on its right,
for every $1\leq i \leq n_t$.
More precisely, we introduce
$
    H_{x\rightarrow y }
:=
    \inf\{s>H(x), \ X(s)=y\}-H(x)
$
for $(x,y)\in\R_+^2$,
which is equal to
$
    H(y)-H(x)
$
if $x<y$. We also introduce
\begin{equation}\label{eqDefUi}
    U_i
 :=
    H(\tilde L_i)-H(\tilde m_i)
=
    H_{\tilde m_i\rightarrow \tilde L_i}
    ,
\qquad
    \B_i^{}
:=
    \big\{ U_i<H_{\tilde m_i \rightarrow \tilde L_i^-}\big\},
\qquad
    i\geq 1.
\end{equation}

We have (notice that $n_t e^{-\k\delta h_t}=o(1)$ since $\phi(t)=o(\log t)$),

\begin{lemma}\label{LemmaProbaQuitterVallesaDroite}
Assume $0<\delta< 1/8$.
We have, for large $t$,
\begin{equation*}
P\left[
    \bigcap_{i=1}^{n_t}
        \left\{
                \P^{\wk}_{\tilde m_i}\big[H\big(\tilde L_i\big)>H\big(\tilde L_i^-\big)\big]
            =
                \P^{\wk}\big({\overline{\B_i}}\big)
            \leq
                e^{-\k h_t/2}
        \right\}
\right]
\geq
    1-C_+ n_t e^{-\k\delta h_t}.
\end{equation*}
\end{lemma}

\smallskip

\noindent{\bf Proof:}
By the strong Markov property and then by \eqref{eqScaleFunctionX}, we have
for all $1\leq i \leq n_t$,
\begin{equation}
    \P^{\wk}(\overline{\B_i})
=
    \P^{\wk}_{\tilde m_{i}} \left(H(\tilde L_{i}) > H(\tilde L_{i}^-) \right)
=
    \Big(\int_{\tilde m_i}^{\tilde L_i} e^{\tV(x)}\textnormal{d}x\Big)
    \Big( \int_{\tilde L_i^-}^{\tilde L_i} e^{\tV(x)}  \textnormal{d}x\Big)^{-1}
\leq
    Q_i/D_i,
\label{eqProbaSortirDuBonCote}
\end{equation}
where,
since
$
    \sup_{[\tilde m_i,\tilde L_i]}\tilde V^{(i)}
\leq
    \sup_{[\tilde m_i,\tilde L_i^+]}\tilde V^{(i)}
=
    \tilde V^{(i)}(\tilde M_i)
$
and
$
    \tilde L_i^- <  \tilde m_i
$,
\begin{align*}
    Q_i
:=
    (\tilde L_i- \tilde L_i^-)  e^{\tilde V^{(i)}( \tilde M_i)},
\qquad
    D_i
:=
    \int_{\tilde L_i^-}^{  \tilde m_i } e^{\tilde V^{(i)}(x)}\dd x.
\end{align*}
We start with the denominator $D_i$.
By \eqref{eqLemmaSansFaggio2}, we have for all $1\leq i \leq n_t$, since $\delta<1/8$,
\begin{equation}
\label{eqProbaInegDi}
    P\big[ D_i \geq e^{h_t^+-\k h_t/2}\big]
\geq
    1- e^{-\k h_t/8}
\geq
    1-e^{-\k\d h_t}.
\end{equation}
We now consider the numerator $Q_i$ for some $1\leq i \leq n_t$.
First by \eqref{eqLemmaSansFaggio1} and since $\tilde L_i< \tilde L_i^+$, we have
$
    P(\tilde L_{i}-\tilde L_{i}^- > 40 h_t^+ / \kappa)
\leq
    e^{- \k h_t/8}
$.
Moreover, since $\tilde \tau_i(h_t)$ is a stopping time,
$(\Vi[x+\tilde\tau_i(h_t)]-h_t,\ x\geq 0)$ is equal in law to $(\wk(x),\ x\geq 0)$, so
\begin{align*}
    P\big[\tV( \tilde M_i) >h_t(1+ \delta)\big]
=
    P\bigg(\sup_{[\tilde\tau_i(h_t),\tilde L_i^+]}\big(\tV-h_t\big)\geq \d h_t\bigg)
\leq
    P\left(\sup_{s \geq 0} \wk(s) >\delta h_t \right)
\leq
    e^{- \delta \kappa h_t},
\end{align*}
since $P[\sup_{[0,\infty)}\wk\geq x]=P[\inf_{[0,\infty)}W_{-\k}\leq -x]=e^{-\k x}$,
$x\geq 0$, e.g. by \eqref{eqScaleFunctionWk}.
Finally
$$
    P\big[ Q_i \leq 40 h_t^+\k^{-1}e^{(1+\delta)h_t}\big]
\geq
    1-2e^{- \delta \kappa h_t}
$$
for $\delta<1/8$ and $t$ large enough. This combined with \eqref{eqProbaSortirDuBonCote}
and \eqref{eqProbaInegDi} gives for large $t$,
\begin{align}\label{eq3.13}
    P\big[
         \P^{\wk}({\overline{\B_i}})
      \leq
         C_+ h_t e^{(1+\delta)h_t-(h_t^+-\k h_t/2)}
      \leq  e^{-(\k/2) h_t}
    \big]
\geq
    1- 3 e^{-\k\delta h_t}.
\end{align}
This proves the lemma.
\hfill $\Box$


\subsection{Probability that the diffusion $X$ goes back to $\tilde m_i$ after leaving the $i$-th valley}\label{SubSectRetourLietoile}
We introduce for $i\in\N^*$  (see Figure \ref{fig1}),
\begin{equation}\label{eqDefLi*}
    \tilde L_i^*
:=
    \inf\{x \geq \tt_i(h_t),\  \wk(x)- W_{\kappa}(\mt_i) \leq 3 h_t/4 \}
<
    \tilde L_i.
\end{equation}

In the next lemma, we show that with a large probability, after hitting $\tilde L_i$,
$X$ hits $\tilde m_{i+1}$ before (maybe) going back to $\tilde L_i^*$ for $1\leq i<n_t$.
This is helpful in
the proof of Lemma \ref{lemtps}
and in Subsection \ref{SubSectProofLocalization}.

\begin{lemma}\label{LemmaProbaRetourEnmi}
We have for large $t$,
\begin{equation}\label{eqProbaSortirDesValleesaGauche}
    P\bigg[
                \P^{\wk}
                \big(
                    \cap_{i=1}^{n_t-1}
                    \big\{H_{\tilde L_i\rightarrow \mt_{i+1}}
                        <H_{\tilde L_i \rightarrow \tilde L_i^*}
                    \big\}
                \big)
            \geq
                1-2n_t e^{-h_t/8}
     \bigg]
\geq
    1-C_+ n_t e^{-\k h_t/16}.
\end{equation}
\end{lemma}
This, combined with Lemma \ref{LemmaProbaQuitterVallesaDroite}, gives a picture of the typical
trajectories of $X$. In particular, with large probability, up to time $t$,
$X$ visits successively the different valleys $[\tilde L_i^-, \tilde L_i]$, $1\leq i< n_t$; it
exits each one on its right $\tilde L_i$, then does not go back to $\tilde L_i^*$ and then to $\tilde m_i$, but goes to $\tilde m_{i+1}$.
That is, with a large probability, the diffusion $X$ hits successively
each $h_t$-minimum $m_i$, $1\leq i\leq n_t$ and does not come back to the previously visited $m_j$.

\noindent {\bf Proof of Lemma \ref{LemmaProbaRetourEnmi}:}
We define $A_\infty^x:=\int_x^{\infty}e^{\wk(u)-\wk(x)}\dd u$, $x\in\R$.
Let
\begin{eqnarray}
    \label{eqDefE1SubScetRetourLietoile}
    \B_1^{\ref{SubSectRetourLietoile}}
& :=  &
    \cap_{i=1}^{n_t-1}\{H_{\tilde L_i \rightarrow \tilde \tau_{i+1}(h_t)}<
                        H_{\tilde L_i \rightarrow \tilde L_i^*}\}
\subset
    \cap_{i=1}^{n_t-1}\big\{H_{\tilde L_i\rightarrow \mt_{i+1}}<H_{\tilde L_i \rightarrow \tilde L_i^*}\big\},
\\
        \B_2^{\ref{SubSectRetourLietoile}}
& := &
    \cap_{i=1}^{n_t-1}\big\{A_{\infty}^{\tilde L_i}\leq e^{h_t/16},\
    A_{\infty}^{\tt_{i+1}(h_t)}\leq e^{h_t/16}, \ A_{\infty}^{\tilde L_i^*}\geq e^{-h_t/16}\big\}.
\nonumber
\end{eqnarray}
Applying Fact \ref{FactDufresnes} (Dufresne) we have
$P(A_{\infty}\geq y)\leq C_+ y^{-\k}$
for $y>0$, and
$P(A_{\infty}\leq y)\leq e^{-1/y}$
for small $y>0$.
Moreover, since $\tilde L_i$, $\tt_{i+1}(h_t)$ and $\tilde L_i^*$ are stopping times for the natural filtration of $\wk$,
the r.v. $A_{\infty}^{\tilde L_i}$, $A_{\infty}^{\tilde \tau_{i+1}(h_t)}$
and $A_{\infty}^{\tilde L_i^*}$ have the same law as $A_\infty$ under $P$ by the strong Markov property.
Consequently for large $t$,
\begin{equation}\label{eqProbaE2Lem35}
    P(\overline{\B_2^{\ref{SubSectRetourLietoile}}})
\leq
    n_t\big(2C_+e^{-\k h_t/16}+e^{-e^{h_t/16}}\big)
\leq
    C_+ n_t e^{-\k h_t/16}.
\end{equation}
Moreover, applying the strong Markov property and \eqref{eqScaleFunctionX},
we have on $\B_2^{\ref{SubSectRetourLietoile}}$,
$$
    \P^{\wk}\big[H_{\tilde L_i \rightarrow \tilde \tau_{i+1}(h_t)}>
            H_{ \tilde L_i \rightarrow  \tilde L_i^*}\big]
=
    \P_{\tilde L_i}^{\wk}\big[H(\tilde\tau_{i+1}(h_t))>H(\tilde L_i^*)\big]
=
    {Q_i^*}/D_i^*,
\qquad
    1\leq i \leq n_t-1,
$$
where, recalling that $\wk(\tilde L_i)=\wk(\tilde m_i)+h_t/2$ and $\wk(\tilde L_i^*)=\wk(\tilde m_i)+3h_t/4$,
\begin{eqnarray}
    Q_i^*
& := &
    \int_{\tilde L_i}^{\tilde\tau_{i+1}(h_t)} e^{\wk(x)}\dd x
\leq
    e^{\wk(\tilde L_i)}A_{\infty}^{\tilde L_i}
\leq
    \exp(\wk(\mt_i)+h_t/2+h_t/16),
\label{Ni} \\
    D_i^*
& := &
    \int_{\tilde L_i^*}^{\tilde\tau_{i+1}(h_t)} e^{\wk(x)}\dd x
=
    e^{\wk(\tilde L_i^*)}A_{\infty}^{\tilde L_i^*}
    -e^{\wk(\tilde \tau_{i+1}(h_t))}A_\infty^{\tilde \tau_{i+1}(h_t)}
\nonumber
\end{eqnarray}
on $\B_2^{\ref{SubSectRetourLietoile}}$.
Moreover,
$\wk(\mt_{i+1})\leq \wk(\tilde L_{i+1}^{\sharp})\leq \wk(\tilde L_i^+) \leq \wk(\mt_i)-h_t$,
so on $\B_2^{\ref{SubSectRetourLietoile}}$ for large $t$,
$$
    D_i^*
\geq
    e^{\wk(\tilde m_i)+11h_t/16}
    -
    e^{\wk(\tilde m_{i+1})+h_t}e^{h_t/16}
\geq
    e^{\wk(\tilde m_i)}\big[e^{11h_t/16}-e^{h_t/16}\big]
\geq
    e^{\wk(\mt_i)+11h_t/16}/2
$$
for all $1\leq i<n_t$, so
$
    Q_i^*/D_i^*
\leq
    2e^{-h_t/8}
$.
Thus
$
    P^{\wk}\big(\overline{\B_1^{\ref{SubSectRetourLietoile}}}\big)
    \un_{\B_2^{\ref{SubSectRetourLietoile}}}
\leq
    2n_t e^{-h_t/8}
$.
This and \eqref{eqProbaE2Lem35} give
\begin{equation}\label{eqProbaE1}
    \P\big(\overline{\B_1^{\ref{SubSectRetourLietoile}}}\big)
\leq
    \P\big(\overline{\B_1^{\ref{SubSectRetourLietoile}}}\cap  \B_2^{\ref{SubSectRetourLietoile}}\big)
    +P\big(\overline{ \B_2^{\ref{SubSectRetourLietoile}}}\big)
\leq
    C_+ n_t e^{-\k h_t/16},
\end{equation}
which we need in the proof of
Lemma \ref{SubSectRetourLietoile}. Moreover, we get
$
    \text{LHS of }\eqref{eqProbaSortirDesValleesaGauche}
\geq
    P\big(\P^{\wk}(\B_1^{\ref{SubSectRetourLietoile}})\geq 1-2n_t e^{-h_t/8}\big)
\geq
    P\big(\B_2^{\ref{SubSectRetourLietoile}}\big)
\geq 1-C_+ n_t e^{-\k h_t/16}
$.
This proves \eqref{eqProbaSortirDesValleesaGauche}.
\hfill$\Box$


\subsection{Independence in a trajectory of $X$}
We are interested in the law of $\bU$, defined as follows:
\begin{eqnarray}
    \tilde{A}_i(z)
& := &
    \int_{\tilde m_i}^{z} e^{\tilde V^{(i)}(u)}\textnormal{d}u,
\qquad
    z \in\R,\ i\in\N^*,
\nonumber\\
\label{eqDefU}
    \bU
& := &
    \int_{\tilde L_2^-}^{\tilde L_2}
    e^{-\tilde V^{(2)}(u)}
    \lo_B\big[\tau^B(\tilde{A}_2(\tilde L_2)),\tilde{A}_2(u)\big] \textnormal{d}u,
\end{eqnarray}
where $(B(s),\ s\geq 0)$ is a Brownian motion independent of $\wk$ and then of $\tilde V^{(2)}$.
As explained below in \eqref{4.38}, $\bU$ is equal in law to the exit time of $X$ from
the valley $\big[\tilde L_2^-, \tilde L_2\big]$ under $\P^{{W_{\kappa}}}_{\tilde m_{2}}$
if $X$ leaves this valley on its right.
Notice that we have chosen the second valley in the definition of $\bU$
because Fact \ref{Fact_Williams} provides the law of $V^{(i)}$ near $m_i$ for $i\geq 2$ but not for $i=1$.
Moreover we stress that $\bU$, as well as the r.v. $U_i$, $i\geq 1$, depend on $t$,
since the $\tilde m_i$, $\tilde L_i^-$ and $\tilde L_i$ depend on $h_t$.


We now prove that the law of the sum $U_1+\dots +U_n$
(defined in \eqref{eqDefUi})
can be approximated by the law of the sum of $n$ independent copies of $\bU$,
in the following sense:


\begin{prop} \label{propind}
Assume $0<\delta<1/8$.
There exists a constant $C_2>0$ such that for large $t$,
\begin{equation}
\label{EqPropPseudoInd}
    \forall \lambda>0,\ \forall 1 \leq n \leq n_t,
\qquad
    \left| \E\left(e^{- \lambda \sum_{i=1}^{ n} U_i } \right)- \left[\E\left(e^{- \lambda \bU} \right) \right]^{ n } \right|
\leq
    C_2 n_t e^{- \delta \kappa h_t},
\end{equation}
where $n_t e^{- \delta \kappa h_t}=o(1)$ since $\phi(t)=o(\log t)$.
Moreover for all $n\geq 2$, $[a,b]\subset[0,1]$ and $\alpha>0$,
\begin{align}
    \left|
        \P\Big(\sum_{i=1}^{n-1}\frac{U_i}{t}\in[a,b], \sum_{i=1}^{n}\frac{U_i}{t}\geq \alpha\Big)
        -
        \int_a^b\P\Big(\sum_{i=1}^{n-1} \frac{U_i}{t}\in\dd x\Big)\P(\bU/t\geq \alpha-x)
    \right|
\leq
    C_+ n_t e^{-  \kappa \delta h_t}.
\label{DifferenceProbasDansLemmeLaplace}
\end{align}
Similarly for $n=1$,
$
    \left|
        \P(U_1/t\geq \alpha)
        -
        \P(\bU/t\geq \alpha)
    \right|
\leq
    C_+ n_t e^{- \kappa \delta  h_t}.
$
\end{prop}


\medskip

\noindent{\bf Proof:}
We fix $0<\delta<1/8$ and $\lambda>0$.
We also introduce $\mathcal{G}_s:=\sigma(X(u),\ 0\leq u\leq s,\ \wk(x),\ x\in\R)$ for $s\geq 0$.
For  $n\in\N^*$ and $1\leq i< n$, we have $\tilde L_i^-<\tilde m_i<\tilde L_i<\tilde m_n$  by \eqref{InegRvTilde1} and \eqref{InegRvTilde2},
so $U_i$ and $\un_{\B_i}$ (defined in \eqref{eqDefUi})
are  $\mathcal{G}_{H(\tilde m_{n})}$-measurable.
Hence
for $t$ and $n$ such that $2\leq n\leq n_t$,
\begin{align*}
    \E^{W_{\kappa}}\bigg(e^{- \lambda \sum_{i=1}^{n}U_i} \prod_{i=1}^{n} \un_{\B_i^{}} \bigg)
& =
    \E^{W_{\kappa}}
        \bigg(
            \E^{W_{\kappa}}
                \bigg(
                    e^{- \lambda U_{n}}
                    \un_{\B_{n}^{}}
                \Big| \mathcal{G}_{H(\tilde m_{n})}
                \bigg)
            e^{- \lambda \sum_{i=1}^{n-1}U_i} \prod_{i=1}^{n-1} \un_{\B_i}
        \bigg)
    \\
&=
    \E^{W_{\kappa}}
        \bigg(
            \E^{W_{\kappa}}_{\tilde m_{n}}
                \bigg(
                    e^{- \lambda H(\tilde L_{n})} \un_{\{H(\tilde L_{n}) <  H(\tilde L_{n}^-)\}}
                \bigg)
            e^{- \lambda \sum_{i=1}^{n-1}U_i} \prod_{i=1}^{n-1} \un_{\B_i}
        \bigg),
\end{align*}
by the strong Markov property, applied at time $H(\tilde m_{n})$ to $X$
which is a Markov process under the quenched probability $\P^{\wk}$.
Hence we obtain by induction on $n$,
\begin{equation}\label{eqProduitLaplacesHLi}
    \E^{W_{\kappa}}
    \bigg(e^{- \lambda \sum_{i=1}^{n}U_i} \prod_{i=1}^{n} \un_{\B_i}
    \bigg)
=
    \prod_{i=1}^{n}  \E^{W_{\kappa}}_{\tilde m_{i}} \left( e^{- \lambda H(\tilde L_{i})}
    \un_{\{H(\tilde L_{i}) < H( \tilde L_{i}^-)\}}   \right).
\end{equation}
Under $\P^{{W_{\kappa}}}_{\tilde m_{i}}$, $(X(u)-\tilde m_i,\ u\geq 0)$
is a diffusion in the potential $\wk(x+\tilde m_i)-\wk(\tilde m_i)=\tV(x+\tilde m_i)$,
$x\in\R$, and starting from $0$. So, applying \eqref{eqDefA} and \eqref{eqFormuleH}
with $\tilde{A}_i (.+\tilde m_i)$ instead of $A(.)$,
there exists a Brownian motion $(\tilde B_i(s),\ s\geq 0)$, independent of $\tV$,  such that
the hitting time of $r\geq 0$ by $X-\tilde m_i$ is under $\P^{{W_{\kappa}}}_{\tilde m_{i}}$,
\begin{equation}\label{eqHittingTimePartantdemi}
\int_{-\infty}^{r} e^{-\tV(x+\tilde m_i)}
    \lo_{\tilde B_i}\big[\tau^{\tilde B_i}
    (\tilde{A}_i (r+\tilde m_i)),\tilde{A}_i(x+\tilde m_i)\big] \textnormal{d}x,
\end{equation}
and is in fact $H(r+\tilde m_i)$, hitting time of $r+\tilde m_i$ by $X$. So
under $\P^{{W_{\kappa}}}_{\tilde m_{i}}$ on $\{H(\tilde L_{i}) < H(\tilde L_{i}^-) \}$,
\begin{equation}
    H(\tilde L_i)
=
    \bU_i
:=
    \int_{\tilde L_i^-}^{\tilde L_i} e^{-\tV(u)}
    \lo_{\tilde B_i}\big[\tau^{\tilde B_i}(\tilde{A}_i (\tilde{L}_i)),\tilde{A}_i(u)\big]
    \textnormal{d}u,
\label{4.38}
\end{equation}
where we replaced the born $-\infty$ in the integral by $\tilde L_i^-$ because
$
e^{-\tV(u)}\lo_{\tilde B_i}[\tau^{\tilde B_i}(\tilde{A}_i (\tilde{L}_i)),\tilde{A}_i(u)]
$
is by \eqref{eqLocalTime} equal
under $\P^{{W_{\kappa}}}_{\tilde m_{i}}$
to
$
\mathcal{L}_X(H(\tilde L_i), u)
$,
which is $0$ for $u<\tilde L_i^-$ on $\{H(\tilde L_{i}) < H(\tilde L_{i}^-) \}$.
For the same reason, we also have,
under $\P^{{W_{\kappa}}}_{\tilde m_{i}}$,
\begin{equation}\label{eqUiLocalTime}
    \bU_i
=
    \mathcal{L}_X(H(\tilde L_i), [\tilde L_i^-, \tilde L_i])
,
\end{equation}
where
$\mathcal{L}_X(s,\Delta):=\int_\Delta \mathcal{L}_X(s,x)\dd x$, $\Delta\subset \R$
is the total time spent by $X$ in $\Delta$ until time $s$.
Also, let $\mathcal{L}_X([u,v],\Delta):=\mathcal{L}_X(v,\Delta)-\mathcal{L}_X(u,\Delta)$ for $0\leq u<v$.
We have by \eqref{4.38},
\begin{equation}\label{eqetInegProduitsLaplace}
    \eqref{eqProduitLaplacesHLi}
=
    \prod_{i=1}^{n}  \E^{W_{\kappa}}_{\tilde m_{i}}
                \left(
                        \exp(- \lambda \bU_i)
                        \un_{\{H(\tilde L_{i}) < H( \tilde L_{i}^-)\}}
                \right)
\leq
    \prod_{i=1}^{n}  \E^{W_{\kappa}}_{\tilde m_{i}}
                \left(
                        \exp(- \lambda \bU_i)
                \right).
\end{equation}
%
We notice that on $\{H(\tilde L_{i}) > H( \tilde L_{i}^-)\}$ under $\P_{\tilde m_i}^{\wk}$, we have
thanks to  \eqref{eqUiLocalTime},
$
   \bU_i
\geq
    \mathcal{L}_X\big([H(\tilde L_i^-)+H_{\tilde L_i^-\to \tilde m_i},H(\tilde L_i)], [\tilde L_i^-, \tilde L_i]\big)
$,
which is the time spent in $ [\tilde L_i^-, \tilde L_i]$ by $X$
between times $H(\tilde L_i^-)+H_{\tilde L_i^-\to \tilde m_i}$
and $H(\tilde L_i)$. So, we get
\begin{eqnarray*}
&&
        \E^{W_{\kappa}}_{\tilde m_{i}}
                \Big(
                        e^{- \lambda \bU_i}
                        \un_{\{H(\tilde L_{i}) > H( \tilde L_{i}^-)\}}
                \Big)
\\
& \leq &
        \E^{W_{\kappa}}_{\tilde m_{i}}
            \Big(
                \un_{\{H(\tilde L_{i}) > H( \tilde L_{i}^-)\}}
                \E^{W_{\kappa}}_{\tilde m_{i}}
                \Big(
                        e^{- \lambda
                            \mathcal{L}_X([H(\tilde L_i^-)+H_{\tilde L_i^-\to \tilde m_i},H(\tilde L_i)],
                            [\tilde L_i^-, \tilde L_i])
                        }
                \Big|
                \mathcal{G}_{[H(\tilde L_i^-)+H_{\tilde L_i^-\to \tilde m_i}]\wedge H(\tilde L_i)}
                \Big)
            \Big)
\\
& = &
    P^{W_{\kappa}}_{\tilde m_{i}}\big[H(\tilde L_{i}) > H( \tilde L_{i}^-)\big]
    \E^{W_{\kappa}}_{\tilde m_{i}}
                \big(
                        e^{- \lambda \bU_i}
                \big)
=
    P^{W_{\kappa}}\big(\overline{\B_i}\big)
    \E^{W_{\kappa}}_{\tilde m_{i}}
                \big(
                        e^{- \lambda \bU_i}
                \big)
\end{eqnarray*}
by the strong Markov property.
This and \eqref{eqetInegProduitsLaplace} yield
\begin{equation*}
    \eqref{eqProduitLaplacesHLi}
 =
    \prod_{i=1}^{n}
        \bigg[
        \E^{W_{\kappa}}_{\tilde m_{i}}
                \left(
                        e^{- \lambda \bU_i}
                \right)
        -
        \E^{W_{\kappa}}_{\tilde m_{i}}
                \left(
                        e^{- \lambda \bU_i}
                        \un_{H(\tilde L_{i}) > H( \tilde L_{i}^-)}
                \right)
        \bigg]
 \geq
        \prod_{i=1}^{n}
        \big[
            1
            -
            P^{W_{\kappa}}\big(\overline{\B_i}\big)
        \big]
        \E^{W_{\kappa}}_{\tilde m_{i}}
                \left(
                        e^{- \lambda \bU_i}
                \right).
\end{equation*}
Consequently by Lemma \ref{LemmaProbaQuitterVallesaDroite},
$$
    \E\left(e^{- \lambda \sum_{i=1}^{ n} U_i } \right)
\geq
    \E\big[\eqref{eqProduitLaplacesHLi}\big]
\geq
        \Big[
            1
            -
            e^{-\k h_t/2}
            \Big]^n
        E\bigg[
            \prod_{i=1}^{n}
                \E^{W_{\kappa}}_{\tilde m_{i}}
                \left(
                        e^{- \lambda \bU_i}
                \right)
        \bigg]
    -
    \frac{C_+n_t}{e^{\k\delta h_t}}.
$$
Since $(\tV(x),\ \tilde {L}^-_i\leq x \leq \tilde {L}_i )$, $i\geq 1$ are i.i.d.
by Lemma \ref{CVs}, the r.v. $\E^{W_{\kappa}}_{\tilde m_{i}}\left(e^{- \lambda \bU_i}\right)$
are also i.i.d, and $\bU_n\egloi \bU$ under $\P$.
Moreover,  $(1-x)^n\geq 1- x n $ for $0\leq x \leq 1$, so this gives for large $t$ since $\delta<1/2$,
\begin{equation}\label{InegUpperBoundLaplaceSommeUi}
    \E\left(e^{- \lambda \sum_{i=1}^{ n} U_i } \right)
\geq
    \Big[\E\big(e^{-\lambda \bU}\Big)\Big]^n
    -
    C_+ n_t e^{-\k \delta h_t},
\end{equation}
for all $1\leq n\leq n_t$.
Similarly, using \eqref{eqetInegProduitsLaplace} and  Lemma \ref{LemmaProbaQuitterVallesaDroite},
$$
    \E\left(e^{- \lambda \sum_{i=1}^{ n} U_i } \right)
\leq
    \E\big[\eqref{eqProduitLaplacesHLi}\big]
    +
    \P\big[\cup_{i=1}^n \overline{\B_i}\big]
\leq
    \Big[\E\big(e^{-\lambda \bU}\Big)\Big]^n
    +
    C_+ n_t e^{-\k \delta h_t},
$$
for every $1\leq n\leq n_t$.
This together with \eqref{InegUpperBoundLaplaceSommeUi} proves \eqref{EqPropPseudoInd}.

For \eqref{DifferenceProbasDansLemmeLaplace}, we obtain for $n\geq 2$, $[a,b]\subset[0,1]$ and $\alpha>0$,
\begin{equation*}
    \P\bigg(\sum_{i=1}^{n-1}\frac{U_i}{t}\in[a,b], \sum_{i=1}^{n}\frac{U_i}{t}\geq \alpha,\ \B_n\bigg)
 =
    \E\Big[\un_{\sum_{i=1}^{n-1}\frac{U_i}{t}\in[a,b]}
    \E^{\wk}\Big(\un_{\frac{U_n}{t}\geq \alpha-\sum_{i=1}^{n-1}\frac{U_i}{t}}\un_{\B_n}
    \Big|\mathcal{G}_{H(\tilde m_{n})}\Big)
    \Big].
\end{equation*}
Since $U_i$ is, for $1\leq i \leq n-1$, $\mathcal{G}_{H(\tilde m_{n})}$-measurable,
whereas $U_n$ and $\un_{\B_n}$ are under $\P^{\wk}$ independent of $\mathcal{G}_{H(\tilde m_{n})}$
by strong Markov property, this is equal to
\begin{equation*}
    \E\bigg[\un_{\sum_{i=1}^{n-1}\frac{U_i}{t}\in[a,b]}
        \bigg(\E^{\wk}\Big(\un_{\frac{U_n}{t}\geq \alpha-x}\un_{\B_n}\Big)_{\big|x=\sum_{i=1}^{n-1}\frac{U_i}{t}}\bigg)
    \bigg]
 =
    \int_a^b
    \P\bigg(\sum_{i=1}^{n-1} \frac{U_i}{t}\in\dd x\bigg)
    \P(U_n/t\geq \alpha-x, \B_n),
\end{equation*}
since $\sum_{i=1}^n U_i/t$ is  measurable with respect to
$\sigma(\wk(v),\ v\leq \tilde L_{n-1},\ X(u),\ u\leq H(\tilde L_{n-1}))$,
and so is independent of
$\E^{\wk}\big(\un_{U_n/t\geq \alpha-x}\un_{\B_n}\big)$
which is for every $x\in\R$
measurable with respect to
$\sigma\big(\wk\big(v+\tilde L_{n-1}^+\big)-\wk\big(\tilde L_{n-1}^+\big),\ v\geq 0\big)$
with $\tilde L_{n-1}\leq \tilde L_{n-1}^+\leq \tilde L_n^-$
by \eqref{InegRvTilde1}.

By \eqref{4.38},  $U_n=H(\tilde L_n)=\bU_n$ on $\B_n$ under $\P_{\tilde m_n}^{\wk}$, and $\bU_n\egloi \bU$ under $\P$ by Lemma \ref{CVs}.
So,
\begin{equation}\label{InegProbasUn}
    \P(\bU/t\geq \alpha-x)-\P(\overline{ \B_n})
\leq
    \P(U_n/t\geq \alpha-x, \B_n)
\leq
    \P(\bU/t\geq \alpha-x),
\qquad
    x\in\R.
\end{equation}
Since $\P(\overline{\B_n})\leq C_+ n_t e^{-\k\delta h_t}$ by Lemma \ref{LemmaProbaQuitterVallesaDroite},
we get \eqref{DifferenceProbasDansLemmeLaplace}  for $n\geq 2$.
Finally the case $n=1$ follows from \eqref{InegProbasUn}.
\hfill$\Box$

\subsection{Negligible parts in the trajectory of $X$}

\noindent We now prove that the total time spent by the diffusion $X$ between the first $n_t$ valleys is negligible compared to $t$.

We first give some estimates about hitting times.
To this aim, we recall the notation of Subsection \ref{SubSectNotationX} and in particular $H_-(r)$
and $H_+(r)$, which are defined in \eqref{eqDefH-H+} and \eqref{eqLocalTime}.
We start with an estimate concerning
the total time spent by $X$ in $\R_-$, that is,
$
    H_-(+\infty)
:=
    \lim_{r\to+\infty}H_-(r)
$.


\begin{lemma}\label{remarquenegatif}
For $z$ large enough (this lemma is in fact true for every $\k\in(0,\infty)$),
 \begin{eqnarray}
    \P(H_-(+\infty)>z)
& \leq &
    C_+[(\log z)/z]^{\k/(\k+2)}.
    \label{eqHmoins}
 \end{eqnarray}
\end{lemma}

\noindent{\bf Proof:}
For $a>0$, $\a>0$ and $b>0$, let
\begin{equation*}
    \B_{1}^{\ref{remarquenegatif}}
:=
    \big\{\sup\nolimits_{x<0} e^{-\wk(x)}\leq a
    \big\},
\quad
    \B_{2}^{\ref{remarquenegatif}}
:=
    \left\{A_{\infty}\leq\a\right\},
\quad
    \B_{3}^{\ref{remarquenegatif}}
:=
    \left\{\sup\nolimits_{y<0} \lo_B[\tau^B(\a),y]\leq b\right\},
\end{equation*}
 \begin{equation*}
     \lneg
   :=   \sup\limits_{r\geq 0} \ \sup\limits_{x<0} \lx(H(r),x).
\end{equation*}
We first prove an inequality with regards to $\lneg$.
We notice that by \eqref{eqLocalTime},
\begin{equation*}
    \lneg
=
    \sup\limits_{r\geq 0} \ \sup\limits_{x<0}
    \Big\{e^{-\wk(x)}\lo_B\big[\tau^B(A(r)),A(x)\big]\Big\}
\leq
    \Big(\sup\limits_{x<0} e^{-\wk(x)}\Big)
    \sup\limits_{y<0} \lo_B[\tau^B(A_{\infty}),y].
\end{equation*}
By the first Ray--Knight theorem (see e.g.\ Revuz and Yor \cite{RevuzYor}, chap.~XI), there exist two Bessel processes $R_2$ and
$R_0$, of dimensions $2$ and $0$ respectively, starting from
$0$ and $R_2(\a)$, such that
$\lo_B(\tau^B(\a),x)$ is equal to $R_2^2(\a-x)$ for $x\in[0,\a]$ and
to $R_0^2(-x)$ for $x<0$. Hence, for $\a \le b$,
\begin{equation*}
    \P\big(\overline{ \B_{3}^{\ref{remarquenegatif}}}\big)
=
    \P\big(R_2^2(\a)>b\big)
    +\int_0^b \P\left( \sup_{y>0}R_0^2(y)>b \big| R_0^2(0)=x\right)
    \P\big( R_2^2(\a) \in [x,x+\text{d}x] \big).
\end{equation*}
Since $R_2$ is equal in law to the euclidian norm of a $2$-dimensional Brownian motion
and so $R_2^2(\alpha)$ is exponentially distributed with mean $2\alpha$,
since $\P(\sup_{\R_+}R_0^2>b  |R_0^2(0)=x)=x/b$, $0\leq x\leq b$ (a scale function of $R_0^2$ being $x\mapsto x$, see e.g.
Revuz and Yor, \cite{RevuzYor} p. 442 with $\nu=-1$),
and since $e^{-x}\leq 1/x$, $x> 0$,
this gives by scaling
\begin{equation}\label{InegProbaMaxLocBsurR-}
    \P(\overline{ \B_{3}^{\ref{remarquenegatif}}})
=
    \exp\bigg(-\frac{b}{2\a}\bigg)
    +
    \E\bigg(\frac{R_2^2(\a)}{b} {\bf 1}_{ \{ R_2^2(\a)\leq b\}}\bigg)
\leq
    \left[2+\E\left(R_2^2(1)\right)\right]\frac{\a}{b}
=
    \frac{4\a}{b}.
\end{equation}
Now, let $z>0$,
$a:=z^{\frac{1}{\k+2}}$, $\a:=z^{\frac{1}{\k+2}}$ and
$
    b
:=
    z^{\frac{\k+1}{\k+2}}
$.
Notice that on
$\B_{1}^{\ref{remarquenegatif}}\cap \B_{2}^{\ref{remarquenegatif}}\cap \B_{3}^{\ref{remarquenegatif}}$,
$\lneg\leq z$.
Moreover,
$P[\sup_{[0,\infty)}\wk\geq x]=e^{-\k x}$,
$x\geq 0$ by \eqref{eqScaleFunctionWk} as in the proof of Lemma \ref{LemmaProbaQuitterVallesaDroite}.
So, using Fact \ref{FactDufresnes} (Dufresne) for $\P\big(\overline{ \B_{2}^{\ref{remarquenegatif}}}\big)$,
we get for $z$ large enough,
\begin{equation}
    \P[\lneg>z]
\leq
    \P(\overline{ \B_{1}^{\ref{remarquenegatif}}})
    +\P(\overline{ \B_{2}^{\ref{remarquenegatif}}})
    +\P(\overline{ \B_{3}^{\ref{remarquenegatif}}})
 \leq
    a^{-\k}
    +[2/\a]^\k/(\k\Gamma(\k))
     +4\a/b
\leq
    C_+ z^{-\frac{\k}{\k+2}}.
\label{eqprobasEMNO}
\end{equation}
We now turn back to $H_-(+\infty)
$.
Define for $c>0$,
\begin{equation*}
    \B_{4}^{\ref{remarquenegatif}}
:=
    \Big\{\min_{0\leq s\leq \tau^B(\a)}B(s)>-\a z^{\frac{\k+1}{\k+2}}\Big\},
\qquad
    \B_{5}^{\ref{remarquenegatif}}
:=
    \Big\{\big|A^{-1}(-z)\big|\leq c\log z\Big\}.
\end{equation*}
On
$\cap_{i=1}^5 \B_{i}^{\ref{remarquenegatif}}$,
notice that by \eqref{eqLocalTime}, for $r\geq 0$,
$\mathcal{L}_X(H(r),x)=0$  if $A(x)\leq \min_{0\leq s \leq \tau^B(A(r))} B(s)$,
and in particular if $A(x)<-\a z^{\frac{\k+1}{\k+2}}=-z$ since $\a\geq A_\infty\geq A(r)$, and so
$$
    H_-(r)
=
     \int_{A^{-1}\left(-z\right)}^0
     \mathcal{L}_X(H(r),x)
     \text{d}x
\leq
     \int_{A^{-1}\left(-z\right)}^0
     \lneg
     \text{d}x.
$$
This gives on $\cap_{i=1}^5 \B_{i}^{\ref{remarquenegatif}}$,
since $\lneg\leq z$ on this event,
\begin{equation}\label{exp47}
    H_-(+\infty)
\leq
    \big|A^{-1}(-z)\big|\lneg
\leq
    c z\log z.
\end{equation}

\noindent Moreover, for $c>2/\k$, small $\e>0$, and all large $z$,
using $(W(-u), \ u\in\R)\egloi (W(u), \ u\in\R)$,
\begin{align}
&
    \P\big(\overline{ \B_{5}^{\ref{remarquenegatif}}}\big)
=
    \P\big[-z<A(-c\log z)\big]
=
    \P\left(z>\int_0^{c\log z}e^{W(-u)+\k u/2}\text{d}u\right)
\nonumber\\
& \leq
    \P\bigg[z>\exp\bigg(\inf_{[0, c\log z]}W\bigg)\frac{2}{\k}(z^{\frac{\k c}{2}}-1)\bigg]
\leq
    \P\left[\exp\left(\inf_{[0, c\log z]}W\right)<z^{1-\frac{\k c}{2}+\e}\right]
\leq
    2z^{-\frac{1}{2c}\left(\frac{\k c}{2}-1-\e\right)^2},
\label{Maine11}
\end{align}
since $\inf_{[0,c\log z]}W\egloi -\sqrt{c\log z} |W(1)|$ and $\P(W(1)\geq x)\leq e^{-x^2/2}$ for large $x$.
Moreover,
$
\P(\overline{ \B_{4}^{\ref{remarquenegatif}}})=\a/[\a+\a z^{\frac{\k+1}{\k+2}}]
\leq z^{-\frac{\k+1}{\k+2}}
$.
Choosing $c$ large enough, this, together with
\eqref{eqprobasEMNO}, \eqref{exp47} and \eqref{Maine11}
gives
$
    \P[H_-(+\infty)>c z\log z]
\leq
    C_+ z^{-\k/(\k+2)}
$,
which proves
\eqref{eqHmoins}.
\hfill$\Box$


\medskip

\begin{lemma}\label{LemTpsVal}
There exists a constant $C_3>0$ such that for every $h>0$,
\begin{equation}\label{eqLemmaInegH1}
    \E\big[H_+(\tau_1^*(h))\big]
\leq
    C_3 e^h,
\\
\end{equation}
where
$\tau^*_1(h)=\inf\big\{u\geq 0,\ \wk(u)-\inf_{[0,u]}\wk\geq h\big\}$
as in \eqref{tauetoile}. Moreover,
\begin{equation}\label{eqLemmaInegH2}
    \P[H_-(\tilde m_1) \geq t/\log h_t ]
\leq
    C_+ [(\log t)^2/t]^{\k/(\k+2)}.
\end{equation}
\end{lemma}

\noindent{\bf Proof:}
First, \eqref{eqLemmaInegH2}  comes directly from Lemma \ref{remarquenegatif}
since $\log h_t\sim_{t\to +\infty} \log\log  t$.

For \eqref{eqLemmaInegH1},
we notice that
by the scale property of $B$,
recalling that $A(u) \geq 0$ for all $u \geq 0$ and $A$ is independent of $B$, we have for every $r\geq 0$, which can depend on the environment $\wk$,
$$
    \eo[H_+(r)]
=
    \eo\left(\int_0^{r} e^{-\wk(u)}A(r)\mathcal{L}_B\left(\tau^B(1), A(u)/A(r)\right)\dd u\right).
$$
We remind that
$
    \E[\mathcal{L}_B(\tau^B(1), y)]
=
    \E[R_2^2(1-y)]
=
    2(1-y)
$
for $0\leq y\leq 1$,
by
the first Ray--Knight theorem, $R_2$ being a $2$-dimensional Bessel process starting from $0$
as in the proof of Lemma \ref{remarquenegatif}.
So by Fubini and due to the independence of $B$ and $\wk$,
\begin{equation}
    \eo[H_+(r)]
=
    \int_0^{r} e^{-\wk(u)}2\left(A(r)-A(u)\right)\dd u
=
    2 \int_0^{r} \int_u^{r} e^{\wk(v)-\wk(u)}\dd v\dd u.\label{eqEsperanceH+}
\end{equation}
Hence, applying this to $r=\tau_1^*(h)$
for $h>0$,
we get
\begin{eqnarray*}
    \E[H_+(\tau_1^*(h))]
& = &
    2 E\bigg(\int_0^{\tau_1^*(h)}
    \int_u^{\tau_1^*(h)} e^{\wk(v)-\wk(u)}\dd v\dd u\bigg)
\\
& \leq &
    2 E\bigg(\int_0^{\infty}\un_{u\leq\tau_1^*(h)}
    \int_u^{\tau_1^*(u,h)} e^{\wk(v)-\wk(u)}\dd v\dd u\bigg),
\end{eqnarray*}
where
$\tau_1^*(u,h):=\inf\{x\geq u, \ \wk(x)-\inf_{[u,x]}\wk\geq h\}\geq \tau_1^*(h)$.
Applying Fubini followed by the Markov property at time $u$, we get
\begin{eqnarray}
    \E[H_+(\tau_1^*(h))]
& \leq &
    2 \int_0^{\infty}E\bigg(\un_{u\leq\tau_1^*(h)}
    \int_0^{\tau_1^*(u,h)-u} e^{\wk(\a+u)-\wk(u)}\dd \a\bigg)\dd u
\nonumber\\
& = &
    2 \int_0^{\infty}E(\un_{u\leq\tau_1^*(h)})
    E\bigg(\int_0^{\tau_1^*(h)} e^{\wk(\a)}\dd \a\bigg)\dd u
=
    2E\big(\tau_1^*(h)\big)\beta_0(h),
\label{eqBeta1_et_Beta2}
\end{eqnarray}
where, similarly as in Enriquez et al. {(\cite{ESZ2}, Lem. 4.9)},
\begin{equation*}
    \beta_0(h)
:=
    E\Big(\int_0^{\tau_1^*(h)} e^{\wk(u)}\dd u\Big).
\end{equation*}

We now prove that $\beta_0(h)\leq C_+ e^{(1-\k)h}$.
We notice that $\wk(u)\leq h$ for all $0\leq u\leq \tau_1^*(h)$ and so
$\mathcal{L}_{\wk}(\tau_1^*(h), x)=0$ for all $x>h$. Consequently,
by the occupation time formula and Fubini,
$$
    \beta_0(h)
=
    E\bigg(\int_{-\infty}^h e^x \mathcal{L}_{\wk}(\tau_1^*(h), x) \dd x\bigg)
\leq
    E\bigg(\int_{-\infty}^h e^x \mathcal{L}_{\wk}(\infty, x) \dd x\bigg)
=
    \int_{-\infty}^h e^x E\big[\mathcal{L}_{\wk}(\infty, x) \big]\dd x,
$$
where
$
    \mathcal{L}_{\wk}(\infty, x)
=
    \lim_{u\to +\infty}\mathcal{L}_{\wk}(u, x)
$.
Moreover,
$E\big[\mathcal{L}_{\wk}(\infty, 0)\big]=2/\k<\infty$,
since $\mathcal{L}_{\wk}(\infty, 0)$ is an exponential variable of mean $2/\k$
(see e.g. Borodin et al. \cite{BorodinSalminem}, p. 90 at the end of paragraph V.11).
Furthermore  by the strong Markov property,
$$
    E\big[\mathcal{L}_{\wk}(\infty, x)\big]
=
    P\big[\tau^{\wk}(x)<\infty\big]E\big[\mathcal{L}_{\wk}(\infty, 0)\big]
=
    \big[\un_{(-\infty, 0)}(x)+e^{-\k x}\un_{(0,\infty)}(x)\big]
    2/\k,
\qquad
    x\in\R,
$$
since $P\big[\tau^{\wk}(x)<\infty\big]=e^{-\k x}$ for $x>0$ by \eqref{eqScaleFunctionWk}, and is $0$ for $x\leq 0$
since $\lim_{+\infty}\wk =-\infty$ a.s.
So for large $h$,
$$
    \beta_0(h)
\leq
    \int_{-\infty}^0 \frac{2}{\k}e^x \dd x
    +\int_0^h \frac{2}{\k}e^{(1-\k)x}\dd x
\leq
    C_+ e^{(1-\k)h}.
$$
This, together with \eqref{eqBeta1_et_Beta2}
and
$
    E\big(\tau_1^*(h)\big)
\leq
    C_+ e^{\kappa h}
$
provided in Fact \ref{Fact_Faggio}
gives \eqref{eqLemmaInegH1}.
\hfill$\Box$

{

We now have all the tools needed to bound the time spent by $X$ between the valleys. We recall that
$
U_i= H(\tilde L_i)-H(\tilde m_i)
$
for $i\geq 1$.
More precisely, we prove the following lemma:

\begin{lemma} \label{lemtps}
Assume $0<\delta<2^{-3/2}$ and $(1+2\delta)\k<1$.
For $t$ large enough,
$$
    \P\left(H(\tilde m_{1} )\leq \frac{2t}{\log h_t} \right)
\geq
    \P\left(\bigcap_{k=1}^{ n_t}\left\{0\leq H(\tilde m_{k})-\sum_{i=1}^{k-1}U_i \leq \frac{2t}{\log h_t} \right\}\right)
\geq
    1- C_+ n_t (\log h_t) e^{- \phi(t)},
$$
where
$\sum_{i=1}^0\dots =0$ by convention.
Notice
that
$
    n_t (\log h_t) e^{- \phi(t)}
\leq
    (\log \log t)e^{[(1+\delta)\k-1]\phi(t)}
=
    o(1)
$ as $t\to+\infty$ since $\log \log t=o(\phi(t))$.
\end{lemma}

\noindent {\bf Proof :}	
We use the notation $\tilde L_i^*$ defined in \eqref{eqDefLi*}.
We introduce
\begin{equation}\label{eqDefXi}
    X_i(t)
:=
    X\big(t+H\big(\tilde L_i\big)\big),
\qquad
    X_i^*(t)
:=
    X\big(t+H\big(\tilde L_i^*\big)\big),
\qquad
    t\geq 0,
\end{equation}
which are diffusions in the environment $\wk$, starting respectively from $\tilde L_i$ and $\tilde L_i^*$, by the strong Markov property.
We also denote by $H_{X_i}(r)$ the hitting time of $r$ by $X_i$ for $r\in\R$.

We first notice that since $U_i=H(\tilde L_i)-H(\tilde m_i)$, $i\in\N^*$,
\begin{align}\label{eqSommeHmk}
    H(\tilde m_k)
=
    H(\tilde m_1)+\sum_{i=1}^{ k-1}U_i+
    \sum_{i=1}^{ k-1}\left(H(\tilde m_{i+1})-H(\tilde{L}_i)\right),
\qquad
    1\leq k\leq  n_t,
\end{align}
and $H(\tilde m_{i+1})-H(\tilde L_i)\geq 0$
since $\tilde m_{i+1}>\tilde L_i$ by \eqref{InegRvTilde1}.
So, we just have to prove that
$H(\tilde m_1)+\sum_{i=1}^{ n_t-1}\big(H(\tilde m_{i+1})-H(\tilde{L}_i)\big)\leq 2t/\log h_t$ with large probability.

The idea of the proof is to use Lemma \ref{LemmaProbaRetourEnmi},
which says that on some large event $\B_1^{\ref{SubSectRetourLietoile}}$,
the diffusion $X_i$ starting from $\tilde L_i$ hits $\tilde m_{i+1}$ before $\tilde L_i^*$.
This  allows us to write  (see step 2)
$
    E^{\wk}[H_{X_i}(\tilde m_{i+1})\un_{\B_1^{\ref{SubSectRetourLietoile}}}]
\leq
    E^{\wk}_{\tilde L_i^*}[H_+(\tilde \tau_{i+1}(h_t))]
$.
Thanks to some large event studied in Step 1, we can compare the
expectancy of this last quantity with
$
    \E[H_+(\tilde \tau_{1}^*(h_t))]
$,
which we
can
bound by Lemma \ref{LemTpsVal}.

\noindent {\bf Step 1: } In this step, we prove that $P(\overline{\B_2^{\ref{lemtps}}})\leq C_+ n_t e^{-\k h_t/2}$,
where,
\begin{eqnarray*}
    \tilde \tau_{i+1}^*(h_t)
& := &
    \inf\{u\geq \tilde L_i^*,\ \wk(u)-\inf_{[\tilde L_i^*,u]}\wk\geq h_t\}
\leq
    \tilde\tau_{i+1}(h_t),
\qquad
    i\geq 1,
\\
    \B_2^{\ref{lemtps}}
& := &
    \cap_{i=1}^{n_t-1}\big\{\tilde \tau_{i+1}^*(h_t)=\tilde \tau_{i+1}(h_t)\big\},
\end{eqnarray*}
where we used, for the inequality,  $\tilde L_i^*<\tilde L_i< \tilde L_{i+1}^\sharp$,
coming from  \eqref{eqDefLi*} and  \eqref{InegRvTilde1}.
By definition of $\tilde \tau_{i+1}(h_t)$ and \eqref{InegLiPremieresDescentes}, we observe that
$
    \big\{\tilde \tau_{i+1}^*(h_t)\neq \tilde \tau_{i+1}(h_t)\big\}
=
    \big\{\tilde \tau_{i+1}^*(h_t)\leq \tilde L_{i+1}^{\sharp}\big\}
=
    \big\{\inf_{[\tilde L_i^*, \tilde \tau_{i+1}^*(h_t)]}\wk -\wk(\tilde L_i^*)\geq -2h_t^+-3h_t/4\big\}
$.
So, applying the strong Markov property at stopping time $\tilde L_i^*$ yields
$$
    P[\tilde \tau_{i+1}^*(h_t)\neq\tilde \tau_{i+1}(h_t) ]
=
    \P\big(\inf\nolimits_{[0,\tau_1^*(h_t)]}\wk \geq -2h_t^+-3h_t/4\big)
\leq
    C_+ h_t e^{-\k h_t}
$$
by \eqref{InegBeta*}.
Then $P\big(\overline{\B_2^{\ref{lemtps}}}\big)\leq C_+ n_t h_t e^{-\k h_t}\leq C_+ n_t e^{-\k h_t/2}$.

\noindent {\bf Step 2:}
On $\B_2^{\ref{lemtps}}$,
$
    H(\mt_{i+1})-H(\tilde L_i)
=
    H_{X_i}(\mt_{i+1})
\leq
    H_{X_i}(\tt_{i+1}(h_t))
=
    H_{X_i}(\tt_{i+1}^*(h_t))
$,
which is, on $\B_1^{\ref{SubSectRetourLietoile}}$ (see \eqref{eqDefE1SubScetRetourLietoile}),
the total time spent by $X_i$ in $[\tilde L_i^*,+\infty)$ before hitting $\tt_{i+1}(h_t)=\tilde \tau_{i+1}^*(h_t)$.
This last quantity is less than or equal to the total time spent in $[\tilde L_i^*,+\infty)$
by $X_i^*$ before hitting $\tt_{i+1}^*(h_t)$.
This is the total time spent in $[0,+\infty)$ by $X_i^*-\tilde L_i^*$
before it first hits $\tilde \tau_{i+1}^*(h_t)-\tilde L_i^*$,
which has the same law as $H_+(\tau_1^*(h_t))$ under the annealed probability $\P$,
since $X_i^*-\tilde L_i^*$ is a diffusion in the environment
$(\wk(\tilde L_i^*+x)-\wk(\tilde L_i^*),\ x\in\R)$, which has on $[0,+\infty)$
the same law as $(\wk(x),\ x\geq 0)$ because
$\tilde L_i^*$ is a stopping time for $\wk$.
Consequently,
\begin{equation}\label{H+tauEtoile}
    \E[(H(\tilde m_{i+1})-H(\tilde L_i))
    1_{\B_1^{\ref{SubSectRetourLietoile}}\cap \B_2^{\ref{lemtps}}  }]
\leq
    \E[H_+(\tau_1^*(h_t))],
\qquad
    1\leq i \leq n_t-1.
\end{equation}
A Markov inequality, this last inequality \eqref{H+tauEtoile}, and then Lemma \ref{LemTpsVal}  lead to
\begin{align}
\label{InegProbaSumUi}
&
    \P\bigg( H_+(\tilde m_1)+\sum_{i=1}^{ n_t-1}\left(H(\tilde m_{i+1})-H(\tilde{L}_i)\right)
    \geq \frac{t}{\log h_t} ,
    \ \B_1^{\ref{SubSectRetourLietoile}}, \ \B_2^{\ref{lemtps}}, \ \B_3^{\ref{lemtps}},\ \mV_t \bigg)
\\
&\leq
    \frac{\log h_t}{t}
    \Big[\E\big[H_+(m_1)\un_{\B_3^{\ref{lemtps}}}\big]
          +(n_t-1)\E\big[H_+(\tau_1^*(h_t))\big]
    \Big]
\leq
    \frac{\log h_t}{t} n_t C_3  e^{h_t} ,
\nonumber
\end{align}
where $\B_3^{\ref{lemtps}}:=\{\mf_{1}\leq \tau_1^*(h_t)\}$ and since $\tilde m_1=m_1$ on $\mV_t$.
Recall that $\phi(t)=o(\log t)$ and $\log\log t=o(\phi(t))$ as $t\to+\infty$, and then
$\log h_t\sim_{t\to+\infty} \log \log t$.
This and \eqref{eqSommeHmk} lead to
\begin{eqnarray}
&&
    \P\bigg(H(\tilde m_{n_t})-\sum_{i=1}^{n_t-1}U_i\geq \frac{2t}{\log h_t}\bigg)
 =
    \P\bigg(
            H(\tilde m_1)
            +\sum_{i=1}^{ n_t-1}\left(H(\tilde m_{i+1})-H(\tilde{L}_i)\right)
            \geq \frac{2t}{\log h_t}
    \bigg)
\nonumber\\
& \leq &
    \P\big(H_-(\tilde m_1)\geq t/\log h_t\big)
    +\eqref{InegProbaSumUi}
    +\P\big(\overline{\B_1^{\ref{SubSectRetourLietoile}}}\big)
    +\P\big(\overline{\B_2^{\ref{lemtps}}}\big)
+\P\big(\overline{\B_3^{\ref{lemtps}}}\big)
    +P\big(\overline{\mV_t}\big)
\nonumber\\
& \leq &
    C_+ n_t (\log h_t)e^{-\phi(t)}.
\label{InegProbaSommeTempsEntreValleesPourLemma35}
\end{eqnarray}
Indeed, we used in the last inequality \eqref{eqLemmaInegH2}, \eqref{eqProbaE1}, Step 1,
Lemma \ref{CVs}, and the fact that
$
    \P(\overline{\B_3^{\ref{lemtps}}})
\leq
    \P(0\leq \Mf_0 <\mf_{1})
\leq
    2h_t e^{-\k h_t}
$
by \eqref{eqProbaM0avantm1} and the definition of $h_t$ maximum and $M_0$.
As explained after \eqref{eqSommeHmk},
this concludes the proof.
\hfill$\Box$


\mysection{Time spent in a standard valley}\label{SectionTimeStandardValley}

The aim of this section is to prove the following proposition,
which gives the second order of the Laplace transform of $\bU$, which is defined in
\eqref{eqDefU} and is useful because of Proposition \ref{propind}:

\begin{prop} \label{proplap}
Assume $\k\in(0,1)$ and $0<\delta<\inf(2/27,\k^2/2)$.
Let $\lambda>0$. As $t\to+\infty$,
\begin{align*}
    e^{ \kappa \phi(t)} \left(1-\E\left(e^{-\lambda \bU/t}\right) \right)
=
    C_{\kappa}\lambda^{\kappa}+o(1),
\end{align*}
where $C_\k:=8^\k(C_0+ |\Upsilon_0|)>0$, with $C_0:={\Gamma(1-\kappa)\Gamma(\k+2)  }/{(1+ \kappa)^{\kappa}}$ and
\begin{equation}\label{eqDefUpsilon0}
    \Upsilon_0
:=
    \k
    \int_0^\infty
    \bigg[
        \frac{ y^{\k}}{[\Gamma(\k+1)I_\k(2\sqrt{ y})]^2}
        -\left(1+\frac{ y}{\kappa+1}\right)^{-2}
    \bigg]
    y^{-\k-1}
    \dd y
<
    0,
\end{equation}
$I_\k$ being the modified Bessel function of the first kind.
\end{prop}

Before proving this in Subsection \ref{SectionRenewal},  we need  additional estimates given below.

\subsection{Some technical estimates}

\noindent Recall that $(R(s),\ s\geq 0)$ is a process with law $\textnormal{BES}(3,\k/2)$, and that
for $a<b$,
$
\big(W_{\kappa}^{b} (s),\ 0 \leq s \leq \tau^{W_{\kappa}^b}(a)\big)
$
is a $(-\kappa/2)$-drifted Brownian motion starting from $b$
and killed when it first hits $a$.
We now introduce
\begin{equation}\label{eqDefF+G+}
    F^{\pm}(x)
:=
    \int_0^{\tau^{R}(x)} \exp(\pm R_{}(s))ds,
\ \
    x>0,
\qquad\quad
    G^{\pm}(a,b)
:=
    \int_0^{\tau^{W_{\kappa}^{b}}(a)} \exp\big(\pm W_{\kappa}^{b} (s)\big)ds,
\quad a<b.
\end{equation}

The following technical lemma is useful to estimate the Laplace transform appearing
in Proposition \ref{proplap}:

\begin{lemma} \label{lem4.5}
There exists $C_4>0$, $M>0$ and $\eta_1\in(0,1)$ such that $\forall y>M, \forall \gamma\in(0, \eta_1]$,
\begin{eqnarray}
\left|E \left( e^{- \gamma F^-(y)}  \right)-[1+2\gamma/(\k+1)]^{-1}\right|
    &\leq&
    C_4\max(e^{-\k y},\gamma^{3/2})
    ,\label{Fmoins}\\
\left|E \left(e^{-\gamma F^+(y)/e^y}\right)-[1-2\gamma/(\k+1)]\right|
    & \leq &
    C_4\max(e^{-\k y},\gamma^{3/2})
    ,\label{Ff2}\\
\left|E \left(e^{-\gamma G^+(y/2,y)/e^y}\right)-[1-\Gamma(1-\kappa)(2\gamma)^{\k}/\Gamma(1+\k)]\right|
    & \leq &
    C_4\max(\gamma^\k e^{-\k y/2},\gamma)
    .\label{Gpf2}
\end{eqnarray}
Moreover, there exists $c_1>0$, such that for all $y >0$,
$E \left[ {F^+(y)}/{e^y}  \right]\leq c_1$.
Finally,
\begin{equation}\label{eqLimiteLaplaceF-}
    \lim_{x\to+\infty} E\left(e^{-\gamma F^-(x)}\right)
=
    \frac{ (2\gamma)^{{\k/2}}}{\k\Gamma(\k)I_\k(2\sqrt{2\gamma})}
\qquad
    \gamma>0.
\end{equation}
\end{lemma}

The proof of this lemma is deferred to Section \ref{SectionAnnexe}.


Before proving Proposition \ref{proplap},
we also need to introduce the following technical lemma,
which is useful  to approximate  $\bU$, and in particular the local time appearing in its expression
\eqref{eqDefU}:

\begin{lemma}\label{LemmaContinuiteEn0}
$(B(t),\ t\in\R)$ being a standard two-sided Brownian motion,
there exists a constant $c_2$ such that for every $0<\e<1$, $0<\eta<1$ and $x>0$,
\begin{align}
&
    \P\bigg(
        \sup_{u\in[-\eta, \eta]}
        \Big|\mathcal{L}_B\big(\tau^B(1), u\big)-\mathcal{L}_B\big(\tau^B(1), 0\big)\Big|
        >\e \mathcal{L}_B\big(\tau^B(1), 0\big)
    \bigg)
\leq\
    c_2\frac{\eta^{1/6}}{\e^{2 /5}}, \label{Dev} \\
&
    \P\Big(\sup_{u\in[0,1]} \mathcal{L}_B\big(\tau^B(1), u\big) \geq x \Big)
\leq\
    4e^{-x/2}, \label{Diel} \\
&
    \P\Big(\sup_{u\leq 0} \mathcal{L}_B\big(\tau^B(1), u\big) \geq x \Big)
\leq\
    \frac{4}{x}. \label{Diel2}
\end{align}
\end{lemma}

\noindent {\bf Proof:}
\noindent
First, \eqref{Diel2} is the particular case $\alpha=1$ of \eqref{InegProbaMaxLocBsurR-},
which we proved in Lemma \ref{remarquenegatif}.

Second, by the first Ray-Knight theorem (see e.g.\ Revuz and Yor \cite{RevuzYor}, chap.~XI),
$\mathcal{L}_B(\tau^B(1),u)=R_2^2(1-u)$ for $u\in[0,1]$, where $R_2^2$ is a
$2$-dimensional squared Bessel process starting from $0$,
so \eqref{Diel}  follows directly from Diel (\cite{Diel} Lem. 2.3 (iii)).

We now turn to the proof of \eqref{Dev}.
Let $0<\e<1$, $0<\eta<1$ and
\begin{equation*}
    \B^{\ref{LemmaContinuiteEn0}}
:=
    \bigg\{
        \sup_{u\in[-\eta, \eta]}
        \Big|\mathcal{L}_B\big(\tau^B(1), u\big)-\mathcal{L}_B\big(\tau^B(1), 0\big)\Big|
        >
        \e \mathcal{L}_B\big(\tau^B(1), 0\big)
    \bigg\}.
\end{equation*}
We have, for $\alpha>0$ and $x>0$,
\begin{eqnarray}
    \P(\B^{\ref{LemmaContinuiteEn0}})
& = &
    \P\big(\B^{\ref{LemmaContinuiteEn0}}\cap\big\{\mathcal{L}_B\big(\tau^B(1), 0\big)\geq \alpha\big\}\big)
    + \P\big(\B^{\ref{LemmaContinuiteEn0}} \cap\big\{\mathcal{L}_B\big(\tau^B(1), 0\big)< \alpha\big\}\big)
    \nonumber\\
& \leq &
    \P\bigg(
        \sup_{u\in[-\eta, \eta]}
        \Big|\mathcal{L}_B\big(\tau^B(1), u\big)-\mathcal{L}_B\big(\tau^B(1), 0\big)\Big|
        >\e \alpha
      \bigg)
    +\P\big[\mathcal{L}_B\big(\tau^B(1), 0\big)< \alpha\big]
    \nonumber\\
& \leq &
    \P\big[\tau^B(1)\geq x\big]
    +\P\bigg(
        \sup_{u\in[-\eta, \eta], \ 0\leq s\leq x}
        \Big|\mathcal{L}_B(s, u)-\mathcal{L}_B(s, 0)\Big|>\e \alpha\bigg)
    +\frac{\alpha}{2},
    \label{eqCalculProbaContinuite0}
\end{eqnarray}
since $\mathcal{L}_B(\tau^B(1), 0)=R_2^2(1)$ is an exponential variable with mean $2$.
Now, notice that
$$
    \P\big[\tau^B(1)\geq x\big]
=
    \P\Big(\sup_{0\leq u\leq x}B(u)<1\Big)
=
    \P\big(|B(x)|<1\big)
\leq
    2/{\sqrt{2\pi x}}.
$$
Let $0<\e_0<1/2$.
The second term of \eqref{eqCalculProbaContinuite0} is less than or equal to
\begin{eqnarray}
&&
    \P\bigg(\sup_{u\in[-\eta, \eta]-\{0\}, \ 0\leq s\leq x}\frac{|\mathcal{L}_B(s, u)-\mathcal{L}_B(s, 0)|}{|u|^{1/2-\e_0}}>\frac{\e \alpha}{\eta^{1/2-\e_0}}\bigg)\nonumber\\
& \leq &
    \frac{\eta^{1/2-\e_0}}{\e \alpha}
    \E\bigg(\sup_{a\neq b, \ 0\leq s\leq x}\frac{|\mathcal{L}_B(s, b)-\mathcal{L}_
    B(s, a)|}{|a-b|^{1/2-\e_0}}\bigg),
\label{EspBarlowYor}
\end{eqnarray}
the last inequality being a consequence of Markov inequality. Now, applying Barlow and Yor (\cite{BarlowYor}, (ii) p. 199 with $\gamma=1$) to
the continuous local martingale $(B(t\wedge x),\ t\geq 0)$
and its jointly continuous local time $(\mathcal{L}_B(t\wedge x,a),\ t\geq 0,\ a\in\R)$
proves that the expectancy in \eqref{EspBarlowYor} is less than or equal to
$
    C_+
    \E\big[\big(\sup_{s\geq 0} |B(s\wedge x)|\big)^{1/2+\e_0}\big]
\leq
    C_+(\sqrt{x})^{1/2+\e_0}
$.
Consequently, \eqref{eqCalculProbaContinuite0} leads to
$$
    \P\big(\B^{\ref{LemmaContinuiteEn0}}\big)
\leq
    2/{\sqrt{2\pi x}}
    +C_+(\sqrt{x})^{1/2+\e_0}\eta^{1/2-\e_0}(\e \alpha)^{-1}
    +\alpha/2.
$$
Now, we choose $\alpha=\e^{-2/5}\eta^{1/5}$, $x=\e^{4/5}\eta^{-2/5}$ 
and $\e_0<1/36$; we get $\P(\B^{\ref{LemmaContinuiteEn0}})\leq C_+\eta^{1/6}\e^{-2/5}$, which  concludes the proof.
\hfill$\Box$

\subsection{Approximation of the exit time from a typical valley}
\noindent \\
We now prove that the standard exit time $\bU$, defined in \eqref{eqDefU},
can be approximated by a product of (sums of) independent r.v,
$(\Ip_1+\Ip_2)(\Im_1+\Im_2){\bf e}_1$.
We need this later to approximate the Laplace transform of $\bU$ and then prove Proposition \ref{proplap},
in particular because we have estimates of the Laplace transforms of these r.v.
$\mathcal{I}_1^{\pm}$ and $\mathcal{I}_2^{\pm}$ in Lemma \ref{lem4.5}.

\begin{prop} \label{lemU1}
Assume
$0<\delta<\inf(2/27,\k^2/2) $
and let $\e_t:= 3e^{-(1-3\delta)h_t/6}$.
Possibly on an enlarged probability space,
there exist random variables $\Ip_1, \Ip_2, \Im_1$ and $\Im_2$, depending on $t$ and
${\bf e}_1$, such that

\noindent {\bf (i)} $\Ip_1$, $\Ip_2$, $\Im_1$, $\Im_2$ and ${\bf e}_1$ are independent;

\noindent {\bf (ii)}
${\bf e}_1$ is exponentially distributed with mean $2$, and
$$
    \Ip_1 \egloi F^+(h_t),
\qquad
    \Ip_2 \egloi G^+(h_t/2,h_t),
\qquad
    \Im_1\egloi \Im_2 \egloi F^-(h_t/2),
$$
where $\egloi$ denotes equality in law, \\
\noindent {\bf (iii)}
for $t$ large enough,
$\P(\mathcal{A}_t)\geq 1-C_+e^{-(c_-)\delta  h_t}$, where
\begin{equation}\label{eqDefAt}
    \mathcal{A}_t
:=
    \left\{\left|\bU-(\Ip_1+\Ip_2)(\Im_1+\Im_2){\bf e}_1\right|
    \leq (\Ip_1+\Ip_2)(\Im_1+\Im_2){\bf e}_1 \e_t \right\}.
\end{equation}
\end{prop}


The proof of this proposition involves 3 lemmas. The first two are straightforward consequence of
what we have already proved and the last one is more technical.


The expressions of $\Im_1$, $\Im_2$, $\Ip_2$ and some intermediate r.v. $\Ip_0$ are given in the following lemma, which also provides their laws.
The random variables
${\bf e_1}$ and $\Ip_1$ are defined later, respectively in Lemma \ref{propIV} and in \eqref{eqDefIp1}.

\begin{lemma}  \label{lemFB}
We have with the notations $F^\pm$ and $G^\pm$ introduced in \eqref{eqDefF+G+},
\begin{align}
    &
    \Ip_0
:=
    \int^{ \tau_2(h_t)}_{ \mf_2} e^{ V^{(2)}(x)} \textnormal{d}x
    \egloi F^{+}(h_t),
\quad
&
    \Im_1
:=
    \int^{ \tau_2(h_t/2)}_{\mf_2} e^{-V^{(2)}(x)} \textnormal{d}x
    \egloi F^{-}(h_t/2),
\label{FB1}
\\
&
    \Ip_2
:=
    \int_{ \tau_2(h_t)}^{{ L}_2} e^{V^{(2)}(x)}
    \textnormal{d}x  \egloi  G^+(h_t/2,h_t) ,
\quad
&
    \Im_2
:=
    \int^{\mf_2}_{ {\tau}_2^-(h_t/2)} e^{ -V^{(2)}(x)} \dd x
\egloi
    F^{-}(h_t/2).
\nonumber
\end{align}
where
$L_2
:=
    \inf\{x>\tau_2(h_t),\ V^{(2)}(x)= h_t/2  \}
$
is defined similarly as $\tilde L_2$ in \eqref{eqDefLiTilde} without tilde,
so that $L_2=\tilde L_2$ on $\mV_t$.

\end{lemma}

\noindent{\bf Proof:}
This is a direct consequence of Fact \ref{Fact_Williams} {\bf (ii)} for $\Ip_0$, $\Im_1$ and $\Im_2$.
and of Fact \ref{Fact_Williams} {\bf (iii)} for $\Ip_2$.
\hfill$\Box$


Recall the notation $\tilde{A}_2(z)=\int_{\tilde m_2}^{z} e^{\tilde V^{(2)}(x)}\textnormal{d}x$ introduced
just after \eqref{eqProduitLaplacesHLi}. We have,

\begin{lemma} \label{mMA}
For all $0<\zeta \leq 1$ and $0<\e<1/2$,
for $t$ large enough,
\begin{align}
    P\big[e^{\zeta h_t (1- \e)}\leq \tilde A_2( \tilde \tau_2 (\zeta h_t)) \leq e^{\zeta  h_t (1+ \e)}\big]
\geq
    1-4 e^{-   \kappa \e\zeta h_t/2}. \label{4.6bb}
\end{align}
\end{lemma}

\noindent{\bf Proof:}
First, notice that on $\mV_t$, by Remark \ref{RemEgaliteAvecouSansTilde},
$
    \tilde A_2( \tilde \tau_2 (\zeta h_t))
=
    \int_{m_2}^{\tau_2 (\zeta h_t)} e^{V^{(2)}(x)}\textnormal{d}x
$,
which is equal in law to $F^+(\zeta h_t)$
thanks to Fact \ref{Fact_Williams} {\bf (ii)}.
Hence by Lemma \ref{CVs},
\begin{eqnarray}
    \nonumber
    \textnormal{LHS of }\eqref{4.6bb}
& \geq &
    P\big[e^{\zeta h_t (1- \e)}\leq  F^{+}(\zeta h_t) \leq e^{\zeta h_t (1+ \e)}\big]-P(\overline{\mV}_t).
\\
& \geq &
    P\big[F^{+}(\zeta h_t) \geq e^{\zeta h_t (1- \e)}\big]
    -P\big[F^{+}(\zeta h_t) > e^{\zeta h_t (1+ \e)}\big]
    -e^{[-\k /2+o(1)]h_t}.~~~~
\label{eqLemmaA1}
\end{eqnarray}
Since $F^+( \zeta h_t)\leq \tau^R(\zeta h_t) e^{\zeta h_t}$,
we have by \eqref{bessel4} for large $t$,
\begin{equation}
    P\big[F^+(\zeta h_t)>e^{\zeta h_t (1+\e)}\big]
\leq
    P\big[\tau^R(\zeta h_t)>e^{\e \zeta h_t}\big]
\leq
    C_+ e^{-\k \zeta h_t/(2\sqrt {2})}.
\label{eqLemmaA2}
\end{equation}
For the lower bound, notice that by \eqref{MinorationAVallee},
$$
    P\big[F^+(\zeta h_t)\geq e^{(1-\e)\zeta h_t}\big]
\geq
    1-3\exp(-\k  \e \zeta h_t/2).
$$
This together with \eqref{eqLemmaA1} and \eqref{eqLemmaA2} proves the lemma.
\hfill$\Box$


\medskip
In the following lemma, we provide an approximation of $\bU$ by the product
$\tA2\mathcal{I}^-{  \bf} {\bf e}_1 $, where $\tA2$ and $\mathcal{I}^-$
depend only on the potential $\wk$, whereas ${\bf e}_1 $ is independent of the potential.

\medskip

\begin{lemma} \label{propIV}
For all $0<\e<\inf(2/27,\k^2/2) $, and $t$ large enough,
\begin{equation}\label{eqPropIV}
    \P\left(
        \left|\bU-\tA2\mathcal{I}^-{\bf e}_1 \right|
            \leq
        2e^{-(1-3 \e)h_t/6}\tA2\mathcal{I}^-{  \bf} {\bf e}_1
    \right)
\geq
    1- C_+e^{-   (c_-)\e h_t},
\end{equation}
where
$$
    \mathcal{I}^-
:=
    \int_{\tilde{\tau}_2^-(h_t/2)}^{ \tilde{\tau}_2(h_t/2)} e^{-\tilde V^{(2)}(u)}  \textnormal{d}u,
\qquad
    {\bf  e}_1
=
    \lo_B\big[\tau^B\big(\tilde{A}_2\big(\tilde L_2\big)\big),0\big]/\tilde{A}_2\big(\tilde L_2\big).
$$
Moreover,
$
    {\bf  e}_1
$ is independent of $\wk$, and exponentially distributed with mean $2$.
\end{lemma}


\noindent {\bf Proof } 
Let $0<\e<\inf(2/27,\k^2/2) $. We first notice that
\begin{align*}
    \bU
=
    \int_{\tilde L_2^-}^{\tilde L_2}
        e^{-\tilde V^{(2)}(u)} \tA2 \lo_{B'}\big[\tau_{B'}(1),\tilde{A}_2(u)/\tA2\big] \textnormal{d}u,
\end{align*}
where  $B'(u):=B\big(\big[\tA2\big]^2 u\big)/\tA2$ for $u\geq 0$, and therefore $(B'(u),\ u\geq 0)$
is by scaling, as $B$,  a standard Brownian motion independent of $\wk$, that is, of $\tilde V^{(2)}$.
The idea of the proof is that, loosely speaking, for $u$ close to $0$,
and more precisely for $u$ between $\tilde \tau_2^-(h_t/2)$ and $\tilde \tau_2(h_t/2)$,
$\lo_{B'}\big[\tau_{B'}(1),\tilde{A}_2(u)/\tA2\big]$
is nearly
$\lo_{B'}\big[\tau_{B'}(1),0\big]={\bf e_1}$,
whereas for $u$ far from $0$, that is $u\notin [\tilde \tau_2^-(h_t/2), \tilde \tau_2(h_t/2)]$,
$e^{-\tilde V^{(2)}(x)}$ is "nearly" $0$, with large probability.

We first notice that ${\bf e_1}=\mathcal{L}_{B'}(\tau^{B'}(1), 0)$
is an exponential r.v. with mean $2$ by the first Ray-Knight theorem, and is independent of $\wk$.
We cut $\bU/\tA2$ into three integrals:
\begin{eqnarray}
    \frac{\bU}{\tA2}
& = &
    \int_{\tilde L_2^-}^{\tilde \tau_2^-(h_t/2)}
    +\int_{\tilde \tau_2^-(h_t/2)}^{\tilde \tau_2(h_t/2)}
    +\int_{\tilde \tau_2(h_t/2)}^{\tilde L_2}
        e^{-\tilde V^{(2)}(u)}  \lo_{B'}\big[\tau_{B'}(1),\tilde{A}_2(u)/\tA2\big] \textnormal{d}u
    \nonumber\\
& =: &
    \mathcal{J}_0+\mathcal{J}_1+\mathcal{J}_2.
    \label{eqDecoupageenJi}
\end{eqnarray}
In what follows, we show that the main contribution comes from $\mathcal{J}_1$.

\noindent{\bf Step 1: study of
$\tilde A_2(u)/\tilde A_2(\tilde L_2)$.
}
We introduce
\begin{equation}\label{eqDefDeltat}
    \delta_t
:=
    e^{ -h_t (1-3\e)/2},
\end{equation}
$$
    \B_1^{\ref{propIV}}
:=
    \big\{\delta_t \tA2\geq  \tilde A_2( \tts2)\big\},
\qquad
    \B_2^{\ref{propIV}}
:=
    \big\{\delta_t \tA2\geq  -\tilde A_2( \tilde\tau_2^-(h_t/2))\big\},
$$
so that on $\B_1^{\ref{propIV}}\cap \B_2^{\ref{propIV}}$,
$\tilde A_2(u)/\tilde A_2(\tilde L_2)\in[-\delta_t,\delta_t]$ for all
$u\in [\tilde \tau_2^-(h_t/2), \tilde \tau_2(h_t/2)]$.
We first prove that
$P\big(\overline{\B_1^{\ref{propIV}}}\big)\leq  C_+ e^{-\k\e h_t/4}$.
By Lemma \ref{mMA}, 
\begin{equation}
    P\big[ \tilde A_2(\tts2) \leq e^{h_t (1+\e)/2} \big]
\geq
    1-C_+ e^{- \k\e  h_t/4}.
\label{eqE4part1}
\end{equation}
Notice that $ \tA2 = \Ip_0 +\Ip_2 $ on $\mV_t$,
and that $\Ip_0\egloi F^+(h_t)$ by Lemma \ref{lemFB}. So
\begin{equation}\label{InegProbaAL2}
    P\big[ \tA2 \geq e^{h_t(1-\e)}\big]
\geq
    P\big[F^+(h_t) \geq   e^{h_t(1-\e)}\big]-P(\overline \mV_t)
\geq
    1-4e^{-\k \e h_t /2},
\end{equation}
where we used \eqref{MinorationAVallee} and Lemma \ref{CVs} in the last inequality. This, together with \eqref{eqE4part1} gives
\begin{align}\label{xo+}
    P(\B_1^{\ref{propIV}})
\geq
   P\big[\delta_t \tA2\geq \delta_t e^{h_t(1-\e)}=e^{h_t(1+\e)/2}\geq \tilde A_2( \tts2)\big]
\geq
    1- C_+ e^{-\k \e h_t/4}.
\end{align}
Similarly on $\mV_t$,
$
    -\tilde A_2(\tilde \tau_2^-(h_t/2))
=
    \int_{\tau_2^-(h_t/2)}^{m_2}\exp(V^{(2)}(s))\dd s
\egloi
F^+(h_t/2)
$
by Fact \ref{Fact_Williams},
and so
$$
    P\big[
        -\tilde A_2(\tilde \tau_2^-(h_t/2))
        \geq
        e^{h_t(1+\e)/2}
    \big]
\leq
    P\big[
        F^+(h_t/2)
        \geq
        e^{h_t(1+\e)/2}
    \big]
+
    P\big(\overline{\mV_t}\big)
\leq
    C_+ e^{-\k h_t/(4\sqrt{2})},
$$
by \eqref{eqLemmaA2} and Lemma \ref{CVs}.
This and \eqref{InegProbaAL2} give for large $t$,
$$
    P\big(\B_2^{\ref{propIV}}\big)
\geq
    1-5e^{-\k\e h_t/2}.
$$

\noindent \textbf{Step 2: study of $\mathcal{J}_2$.}
We prove in this step that for $t$ large enough,
\begin{equation}
    \P\big[\mathcal{J}_2 \geq c_+ h_t^2 e^{-(1- \e) h_t/2}  \big]
\leq
    C_+e^{- \k\e h_t/2  },
\label{eqApproxJ2}
\end{equation}
for some constants $c_+$ and $C_+$.
Let
$
    \B_3^{\ref{propIV}}
:=
    \{\sup_{u\in[0,1]} \mathcal{L}_{B'}(\tau^{B'}(1), u) \leq h_t\}
$, and define
$$
    \B_4^{\ref{propIV}}
:=
    \big\{\inf\nolimits_{[\tts2,\tilde{\tau}_2(h_t)]}\tilde V^{(2)}>(1-\e) h_t/2\big\},
\qquad
    \B_5^{\ref{propIV}}
:=
    \big\{\tilde L_2^{+}-\tilde L_2^- \leq 40 h_t^+ / \kappa\big\}.
$$
We have on
$
\B_3^{\ref{propIV}}\cap \B_4^{\ref{propIV}}\cap \B_5^{\ref{propIV}}
$,
\begin{equation}\label{eqApproxJ2Bis}
    \mathcal{J}_2
\leq
    h_t
    \int_{\tts2}^{\tilde L_2} e^{-\tilde V^{(2)}(u)}\dd u
\leq
    h_t
    e^{- (1- \e) h_t/2} \big[\tilde L_2-\tts2 \big]
\leq
    \frac{40 h_t^+ h_t}{\kappa} e^{- (1- \e) h_t/2}.
\end{equation}
Now, Fact \ref{Fact_Williams}, equation  (\ref{3.10}) with $\alpha=1/2$, $\gamma=(1-\e)/2$ and $\omega=1$,
and Lemma \ref{CVs} give
$$
    P\big(\overline{\B_4^{\ref{propIV}}}\big)
\leq
    P\big[\inf\nolimits_{[\tau_2(h_t/2),\tau_2(h_t)]} V^{(2)}\leq(1-\e) h_t/2,\mV_t\big]
    + P(\overline{\mV}_t)
\leq
    3e^{-\k \e h_t/2}.
$$
Moreover,
$
    P\big(\overline{\B_5^{\ref{propIV}}}\big)
\leq
    e^{-\k h_t/8}
\leq
    e^{-\k \e h_t/2}
$
by \eqref{eqLemmaSansFaggio1} since $\e<1/4$,
and $\P\big(\overline{\B_3^{\ref{propIV}}}\big)\leq 4e^{-h_t/2}$ by (\ref{Diel}).
This, together with \eqref{eqApproxJ2Bis}
proves \eqref{eqApproxJ2}.

\noindent \textbf{Step 3: study of $\mathcal{J}_0$.} We prove that for $t$ large enough,
\begin{equation}\label{eqApproxJ0}
    \P\big[\mathcal{J}_0 \geq 40\k^{-1}  h_t^+ e^{-(1- 4\e) h_t/2}  \big]
\leq
    C_+e^{- \k \e h_t/2}.
\end{equation}

\noindent
Similarly as in Step 2, we introduce
\begin{eqnarray*}
    \B_6^{\ref{propIV}}
& := &
    \big\{\inf\nolimits_{\tilde L_2^- \leq u \leq \tilde \tau^-_2(h_t/2)}\tilde V^{(2)}(u)\geq (1/2-\e)h_t\big\},
\\
    \B_7^{\ref{propIV}}
& := &
    \big\{\sup\nolimits_{  s \leq 0 } \lo_{B'}(\tau^{B'}(1),s) \leq e^{\e h_t}\big\}.
\end{eqnarray*}
Lemma \ref{CVs}, equation \eqref{3.10} with $\gamma=1/2-\e$, $\alpha=1/2$ and $\omega=1$,
and \eqref{eqLemmaSansFaggio3}
give
\begin{eqnarray*}
    P\big(\overline{\B_6^{\ref{propIV}}}\big)
& \leq &
    P\big(
        \inf\nolimits_{\tilde L_2^- \leq u \leq \tilde \tau^-_2(h_t)}\tilde V^{(2)}(u)< (1/2-\e)h_t
      \big)
\\
&&
    +
    P\big(
        \inf\nolimits_{\tilde \tau_2^-(h_t) \leq u \leq \tilde \tau^-_2(h_t/2)}\tilde V^{(2)}(u)< (1/2-\e)h_t
        ,\mV_t
      \big)
    +P\big(\overline{ \mV_t}\big)
    \ \leq \ 3e^{-\kappa \e h_t}.
\end{eqnarray*}
Moreover, by (\ref{Diel2}),
$\P(\overline{\B_7^{\ref{propIV}}} ) \leq 4e^{- \e h_t}$.

\noindent
Therefore, on
$
\B_5^{\ref{propIV}}\cap\B_6^{\ref{propIV}}\cap\B_7^{\ref{propIV}}
$,
i.e with a probability larger than $1- C_+e^{- \k \e h_t/2}$, we obtain
\begin{equation}\label{InegJ0Fin}
    \mathcal{J}_0
\leq
     \sup_{  s \leq 0 } \lo_{B'}[\tau^{B'}(1),s]
    \int^{ \tilde \tau^-_2(   h_t/2 )}_{\tilde L_2^-}e^{ -\tilde V^{(2)}(u)} \dd u
\leq
    40\k^{-1} h_t^+   e^{- (1/2- 2\e)h_t } ,
\end{equation}
which yields \eqref{eqApproxJ0}.

\smallskip
\noindent \textbf{Step 4: study of $\mathcal{J}_1$.} We prove that for $t$ large enough,
\begin{equation}\label{eqProbaApproxJ1}
    \P\big[\mathcal{J}_1 \leq  {\bf e}_1  e^{-\e h_t}/2 \big]
\leq
    C_+ e^{-(c_-)\e h_t}.
\end{equation}
First, recall that ${\bf e_1}=\mathcal{L}_{B'}(\tau^{B'}(1), 0)$ and let
\begin{equation*}
    \B_8^{\ref{propIV}}
 :=
    \left\{
        \sup\nolimits_{s\in[-\delta_t, \delta_t]}
        \Big|\mathcal{L}_{B'}(\tau^{B'}(1), s)-\mathcal{L}_{B'}(\tau^{B'}(1), 0)\Big|
    \leq
        \delta_t^{1/3} \mathcal{L}_{B'}(\tau^{B'}(1), 0)
    \right\}.
\end{equation*}
We know that $\P(\overline{\B_8^{\ref{propIV}}})\leq C_+ \delta_t^{1/30}$ by (\ref{Dev}).
Since on $\B_1^{\ref{propIV}}\cap \B_2^{\ref{propIV}}$,
$\tilde A_2(u)/\tilde A_2(\tilde L_2)\in[-\delta_t,\delta_t]$
for all
$u\in [\tilde \tau_2^-(h_t/2), \tilde \tau_2(h_t/2)]$,
we get on $\B_1^{\ref{propIV}}\cap \B_2^{\ref{propIV}}\cap\B_8^{\ref{propIV}}$,
\begin{equation}\label{eqx0}
    \big(1-\delta_t^{1/3}\big) \mathcal{I}^- {\bf e_1}
\leq
    \mathcal{J}_1
\leq
    \big(1+\delta_t^{1/3}\big) \mathcal{I}^- {\bf e_1}.
\end{equation}
We finally prove that $\mathcal{I}^-$ is not too small,  with a similar argument as before.
First, we have
\begin{equation}\label{InegI-}
    \mathcal{I}^-
\geq
    \int_{\tilde m_2}^{\tilde \tau_2(\e h_t)} e^{-\tilde V^{(2)}(u)}  \textnormal{d}u
\geq
    [\tt_2( \e h_t)-\mt_2]e^{-\e h_t}
\geq
    e^{-\e h_t}
\end{equation}
on $\B_9^{\ref{propIV}}$,
where $\B_9^{\ref{propIV}}:=\{\tt_2(\e h_t)-\mt_2 \geq 1\}$.
Moreover for large $t$,
$$
    P\big(\B_9^{\ref{propIV}}\big)
\geq
    P\big(\tau_2(\e h_t)-\tau_2(\e h_t/2)\geq 1,\mV_t\big)-P\big(\overline{\mV_t}\big)
\geq
    1
    -e^{-\kappa h_t/4}
$$
by Fact \ref{Fact_Williams}, \eqref{3.10b} and Lemma \ref{CVs}.
Let $\B_{10}^{\ref{propIV}}:=\{{\bf e}_1 \geq e^{- \e h_t/2}\}$, and observe that $\P(\B_{10}^{\ref{propIV}}) \geq 1- e^{-\e h_t/2}$ since ${\bf e_1}$ is exponentially distributed with mean $2$.
Since
on
$\B_1^{\ref{propIV}}\cap \B_2^{\ref{propIV}}\cap \B_8^{\ref{propIV}}\cap \B_9^{\ref{propIV}}
\cap \B_{10}^{\ref{propIV}}
$,
\begin{equation}\label{InegJ1}
    \mathcal{J}_1
\geq
    \big(1-\delta_t^{1/3}\big) \mathcal{I}^- {\bf e_1}
\geq
    (1/2)   e^{-\e h_t} {\bf e_1}
\geq
    e^{-2\e h_t}
\end{equation}
for large $t$ by \eqref{eqx0} and \eqref{InegI-}, this gives \eqref{eqProbaApproxJ1}.

\smallskip
\noindent \textbf{Step 5: end of the proof.}
We have on $\cap_{i=1}^{10} \B_i^{\ref{propIV}}$, for $t$ large enough,
by \eqref{eqDecoupageenJi}, \eqref{eqx0} and \eqref{eqDefDeltat},
\begin{equation}\label{InegMinorationU}
    \bU/\tA2
\geq
    {\mathcal{J}}_1
\geq
    (1-e^{ -h_t (1-3\e)/6})\mathcal{I}^-{\bf e_1}.
\end{equation}
Moreover
on $\cap_{i=0}^{10} \B_i^{\ref{propIV}}$, for $t$ large enough,
\begin{equation}\label{InegJ0PlusJ2}
    \mathcal{J}_0+\mathcal{J}_2
\leq
    e^{(-1/2+3\e)h_t}
\leq
    e^{(-1/2+5\e)h_t}    {\mathcal{J}}_1
\end{equation}
by \eqref{InegJ0Fin} and \eqref{eqApproxJ2Bis} for the first inequality,
and \eqref{InegJ1} for the second one.
As a consequence,
\begin{equation}\label{EncadrementU1parI1}
    \bU/\tA2
\leq
    \big[1+e^{(-1/2+5\e)h_t}\big]{\mathcal{J}}_1
\leq
    \big[1+2e^{ -h_t (1-3\e)/6}\big]\mathcal{I}^-{\bf e_1}
\end{equation}
for large $t$, where we used \eqref{eqDecoupageenJi} and \eqref{InegJ0PlusJ2} in the first inequality,
and \eqref{eqx0}, \eqref{eqDefDeltat} and $\e\leq 2/27$ in the second one.
Finally, by \eqref{InegMinorationU} and \eqref{EncadrementU1parI1},
RHS of \eqref{eqPropIV} $\geq 1-\sum_{i=1}^{10}P\big(\overline{\B_i^{\ref{propIV}}}\big)$,
where RHS means right hand side,
which proves the lemma.
\hfill$\Box$

\smallskip

We are now ready prove Proposition \ref{lemU1}, for which we use the notation of Lemma \ref{lemFB}.

\smallskip

\noindent{\bf Proof of Proposition \ref{lemU1}:}
The idea of the proof is that thanks to Lemma \ref{propIV},
we can already approximate  $\bU$
by $\tilde A_2(\tilde L_2) \mathcal{I}^- {\bf e_1}$, which is equal to
$(\Ip_0+\Ip_2)(\Im_1+\Im_2){\bf e}_1$ on $\mV_t$.
However, $\Ip_0$ is not independent of $\Im_1$, so we would like to replace it by
a r.v $\Ip_1$ with the same law and independent of $\Im_1$, $\Im_2$, $\Ip_2$ and ${\bf e_1}$.
We do this by replacing in $\Ip_0$ the small quantity
$\int_{m_2}^{\tau_2(h_t/2)} e^{V^{(2)}(s)} \dd s$ by an independent copy of it.


More precisely,
we define
\begin{equation}\label{eqDefI3+}
        \Ip_3
    :=
        \int_{\tau_2(h_t/2)}^{\tau_2(h_t)} e^{V^{(2)}(s)} \dd s
\leq
        \Ip_0.
\end{equation}
By Fact \ref{Fact_Williams}  {\bf (ii)}, $(V^{(2)}(m_2+s),\ 0\leq s\leq \tau_2(h_t)-m_2)$
is a Markov process, so $\Ip_3$ and $\Im_1$ are independent by the strong Markov property.
Moreover by Fact \ref{Fact_Williams} {\bf (ii)}, $\Im_2$ is independent of this Markov process and then
is independent of $\Ip_3$ and $\Im_1$.
Also, by Fact \ref{Fact_Williams}  {\bf (iii)},
$\Ip_2$ is independent of $(\wk(s),\ s\leq \tau_2(h_t))$ and then of
$\Im_2$, $\Ip_3$ and $\Im_1$. Finally by Lemma \ref{propIV}, ${\bf e_1}$ is independent of $\wk$ and then
$\Ip_2$,  $\Ip_3$, $\Im_1$, $\Im_2$, and ${\bf e_1}$ are independent.

Furthermore, by Lemma \ref{lemFB},
$\Im_1\egloi F^-(h_t/2)$, $\Im_2\egloi F^-(h_t/2)$
and $\Ip_2\egloi G^+(h_t/2,h_t)$,
${\bf e}_1$ is by Lemma \ref{propIV} independent of $\wk$ and exponentially distributed with mean $2$.
Moreover as before, by  Fact \ref{Fact_Williams}, Lemma \ref{lemFB} and the strong Markov property,
these r.v. $\Im_1$, $\Im_2$, $\Ip_2$ and ${\bf e}_1$ are independent of
$(V^{(2)}(s+\tau_2(h_t/2)),\ 0\leq s\leq \tau_2(h_t)-\tau_2(h_t/2))$
which has the same law as a
$\textnormal{BES}(3,\k/2)$
starting from $h_t/2$ and stopped when it first hits $h_t$.

We now consider, possibly on an enlarged probability space, a process
$(R^{(1)}(s),\ 0\leq s\leq \tau^{R^{(1)}}(h_t/2))$, independent of $\wk$ and ${\bf e_1}$ and then independent of
$\Im_1$, $\Im_2$, $\Ip_2$, ${\bf e_1}$ and $(V^{(2)}(s+\tau_2(h_t/2)),\ 0\leq s\leq \tau_2(h_t)-\tau_2(h_t/2))$,
and distributed as $(R(s),\ 0\leq s \leq \tau^{R}(h_t/2))$, $R$ being a  $\textnormal{BES}(3,\k/2)$ process.
We now extend this process by setting
$R^{(1)}(u):=V^{(2)}[u-\tau^{R^{(1)}}(h_t/2)+\tau_2(h_t/2)]$
for $\tau^{R^{(1)}}(h_t/2)\leq u \leq \tau^{R^{(1)}}(h_t/2)+\tau_2(h_t)-\tau_2(h_t/2)$.
By the Strong Markov property, $(R^{(1)}(s),\ 0\leq s\leq \tau^{R^{(1)}}(h_t))$
has the same law as $(R(s),\ 0\leq s\leq \tau^{R}(h_t))$, and then
\begin{equation}\label{eqDefIp1}
    \Ip_1
:=
    \int_0^{\tau^{R^{(1)}}(h_t)} e^{R^{(1)}(s)}\dd s
\egloi
    F^+(h_t),
\qquad
    \Ip_3
=
    \int_{\tau^{R^{(1)}}(h_t/2)}^{\tau^{R^{(1)}}(h_t)} e^{R^{(1)}(s)}\dd s.
\end{equation}
Furthermore, since $(R^{(1)}(s),\ 0\leq s\leq \tau^{R^{(1)}}(h_t))$ is obtained by gluing
two processes independent of  $\Im_1$, $\Im_2$, $\Ip_2$ and ${\bf e_1}$,
it is also independent of these r.v., and so $\Ip_1$ is also
independent of these r.v.  $\Im_1$, $\Im_2$, $\Ip_2$ and ${\bf e_1}$.
This already proves affirmations {\bf (i)} and {\bf (ii)} of Proposition \ref{lemU1}.

Moreover, with the same notation as in Lemma \ref{propIV},
we have on $\mV_t$,
$\Io^-=\Im_1+\Im_2$ and
$
    \tilde A_2(\tilde L_2)=\int_{\tilde m_2}^{\tilde L_2} e^{\tilde V^{(2)}(x)}\textnormal{d}x
=
    \Ip_0+\Ip_2
$,
where
$
    \Ip_0
=
    \int^{ \tau_2(h_t)}_{ \mf_2} e^{ V^{(2)}(x)} \textnormal{d}x
$ as defined in \eqref{FB1}.

We now    approximate $\Ip_0$ by $\Ip_1$.
Since $\Ip_0-\Ip_3=\int_{\mf_2}^{\tau_2(h_t/2)}e^{V^{(2)}(s)}\dd s\egloi F^+(h_t/2)$
by \eqref{eqDefI3+} and Fact \ref{Fact_Williams},
and since
$\Ip_1-\Ip_3=\int_{0}^{\tau^{R^{(1)}}(h_t/2)}e^{R^{(1)}(s)}\dd s\egloi F^+(h_t/2)$,
we get
\begin{eqnarray*}
    P\left(\Ip_1-\Ip_3>e^{(1+\delta) h_t/2} \right)
=
    P\left(\Ip_0-\Ip_3>e^{(1+\delta) h_t/2} \right)
=
    P\left( \frac{F^+(h_t/2)}{e^{h_t/2}} > e^{\delta h_t/2} \right)
\leq
    \frac{c_1 }{e^{\delta h_t/2}},
\end{eqnarray*}
by Markov inequality since
$E \left[ {F^+(y)}/{e^y}  \right]\leq c_1$ for all $y >0$
by  Lemma \ref{lem4.5}.
Moreover, $\Ip_3\leq \Ip_1$,
and by \eqref{MinorationAVallee}, with a probability larger than $1-3e^{-\k \delta h_t/2}$ for large $t$,  $\Ip_1 \geq e^{(1-\delta)h_t} $.
Therefore,  with a probability greater than $1-4 e^{-\k\delta h_t/2}$ for large $t$,
\begin{align*}
&
    \Ip_0
=
    \Ip_3+(\Ip_0-\Ip_3)
\leq
    \Ip_1+e^{(1+ \delta)h_t/2}
\leq
    (1+e^{-(1-3 \delta)h_t/2})\Ip_1,
\\
&
    \Ip_0
\geq
    \Ip_3
=
    \Ip_1-(\Ip_1-\Ip_3)
\geq
    \Ip_1-e^{(1+\delta)h_t/2}
\geq
    (1-e^{-(1-3 \delta)h_t/2}) \Ip_1,
\end{align*}
and then
\begin{align*}
    (1-e^{-(1-3 \delta)h_t/2}) \Ip_1
\leq
    \Ip_0
\leq
    (1+e^{-(1-3 \delta)h_t/2}) \Ip_1.
\end{align*}
This and Lemma \ref{CVs} give with probability at least $1-5 e^{-\k\delta h_t/2}$,
$$
    (1-e^{-(1-3 \delta)h_t/2})
\leq
    \frac{\tilde A(\tilde L_2)\Im}{(\Ip_1+ \Ip_2)(\Im_1+ \Im_2)}
=
    \frac{\Ip_0 +\Ip_2}{\Ip_1+ \Ip_2}
\leq
    (1+e^{-(1-3 \delta)h_t/2}).
$$
since on $\mV_t$, $\Io^-=\Im_1+\Im_2$ and
$
    \tilde A(\tilde L_2)
=
    \Ip_0+\Ip_2
$.
Finally, this and Lemma \ref{propIV} applied with $\e =\delta$
prove affirmation {\bf (iii)} of the proposition.
\hfill$\Box$



\subsection{Second order of the Laplace transform of a standard exit time}\label{SectionRenewal}

\noindent \\
We are now ready
to prove Proposition \ref{proplap}.
This proof is quite technical and can be skipped at first reading.


\noindent \textbf{Idea of the proof of Proposition \ref{proplap}:}
In this proof we use Proposition \ref{lemU1} to approximate the Laplace transform of $\bU$ by that of
$(\Ip_1+\Ip_2)(\Im_1+\Im_2){\bf e_1}$, and take advantage of the fact that the r.v. that appear in this product are independent. This allows us to condition first by ${\bf e_1}$, and then by $\sigma(\Ip_1,\Ip_2)$. We then cut this into 2 parts :
one (studied in Lemma \ref{LemmaPartA}) for which $\Ip_1+\Ip_2$ is "small", which allows us to use the approximations of
the Laplace transforms provided by Lemma \ref{lem4.5},
and one (studied in Lemma \ref{LemmaPartB}) for which $\Ip_1+\Ip_2$ is "big", which allows us to approximate
$\Ip_1+\Ip_2$ by a r.v. having the same law as $A_\infty$, for which we know the density.

\noindent \textbf{Proof of Proposition \ref{proplap}:}
We fix $\k\in(0,1)$, $0<\delta<\inf(2/27,\k^2/2)$  and  $\lambda>0$. We have for every $t>0$, $\mathcal{A}_t$ being defined in \eqref{eqDefAt},
and with $\e_t=3e^{-(1-3\delta)h_t/6}$ as defined in Proposition \ref{lemU1},
\begin{align*}
    \E\left(e^{-\lambda \bU/t }\right)
=
    \E\left(e^{-\lambda \bU/t}\un_{\mathcal{A}_t}\right)
    +\E\left(e^{-\lambda \bU/t}\un_{\overline{{\mathcal{A}}_t}}\right)
\leq
    \E\left(e^{-\lambda \bU/t } \un_{\mathcal{A}_t}\right)
    +\P\left(\overline{{\mathcal{A}}_t}\right).
\end{align*}
Hence by the definition \eqref{eqDefAt} of $\mathcal{A}_t$ and Proposition \ref{lemU1}, we get with $\lambda_t^{\pm}:=2\lambda \left(1\pm \e_t\right)/t$,
\begin{equation}\label{InegLaplaceUS0}
    S_{0,t}^+
    -C_+e^{- (c_-)\delta h_t}
\leq
\E\left(e^{-\lambda \bU/t}\right)
\leq
    S_{0,t}^-+C_+e^{- (c_-)\delta h_t},
\end{equation}
where
$
    S_{0,t}^{\pm}
:=
    \E\left(e^{-(\lambda^\pm_t/2) (\Ip_1+\Ip_2)(\Im_1+\Im_2){\bf e_1}}\right).
$
Let $\theta\in(3\k/4, \k)$,
\begin{equation}\label{eqDefBi}
    \mathcal B_1
:=
    \{\Ip_1+ \Ip_2> t e^{- \theta \phi(t)}  \},
\qquad
    \mathcal B_2
:=
    \{\Ip_1+ \Ip_2 \leq t e^{ - \theta \phi(t)}  \}=\overline{\mathcal B_1}.
\end{equation}
Since ${\bf e_1}/2$ is an exponential r.v. with mean $1$ and is independent of the r.v. $\Io^\pm_i$,
$i\in\{1,2\}$ by Proposition \ref{lemU1},
we have $S_{0,t}^{\pm}=S_{1,t}^{\pm}+S_{2,t}^{\pm}$, where for $i\in\{1,2\}$,
\begin{eqnarray}
    S_{i,t}^{\pm}
& := &
    \E\left(e^{-(\lambda^\pm_t/2) (\Ip_1+\Ip_2)(\Im_1+\Im_2){\bf e_1}}\un_{\mathcal B_i}\right)
\nonumber\\
& = &
    \int_0^\infty E\left(\un_{\mathcal B_i}e^{-z\lambda^\pm_t (\Ip_1+\Ip_2)(\Im_1+\Im_2)}\right) e^{-z}\dd z
    \label{eqKiIntermediaire}
    \\
& = &
    \int_0^{\infty}
    E\Big(  \un_{\mathcal B_i}\left[\left.E\left(
    e^{-z \lambda^{\pm}_t(\Ip_1+\Ip_2)\Im_1}\right|\Ip_1,\Ip_2\right)\right]^2 \Big)
    e^{-z}\textnormal{d}z  ,\label{eqKi}
\end{eqnarray}
since $\Im_1$ and $\Im_2$ are independent and independent of $\Ip_1$ and $\Ip_2$ and have the same law,
once more by Proposition \ref{lemU1}.
We also define,
\begin{equation}\label{eqDefS3}
    Z_t(x)
:=
    \bigg[1+{\frac{2\lambda^\pm_t x}{\kappa+1} (\Ip_1+\Ip_2)}\bigg]^{-1},
    \quad x>0,
\qquad
    S_{3,t}^{\pm}
:=
    \int_0^{ \infty}
    E\left( \un_{\mathcal{B}_1}Z_t^2(z)\right)e^{-z}\textnormal{d}z.
\end{equation}

We start with the study in the case $\Ip_1+\Ip_2$ is "small", that is, $\Ip_1+\Ip_2\leq t e^{-\theta \phi(t)}$
as is the case on $\mathcal{B}_2$. More precisely, we prove

\begin{lemma} \label{LemmaPartA}
As $t\to+\infty$,
with $C_0={\Gamma(1-\kappa)\Gamma(\k+2)  }/{(1+ \kappa)^{\kappa}}$ as defined in
Proposition \ref{proplap},
\begin{equation}\label{eqConclusionPartA}
    S_{2,t}^{\pm}+S_{3,t}^{\pm}
=
    1-C_0 8^{\k}\lambda^\k e^{-\k \phi(t)}+o(e^{-\k\phi(t)}).
\end{equation}
\end{lemma}

\noindent{\bf Proof:}
Let $a(t):=e^{-(3/4)\k \phi(t)}$.
Now, consider
$
    0
\leq
    z
\leq
    \eta_1 a(t) e^{\theta \phi(t)} /(t \lambda_t^\pm)
=
    \eta_1 a(t) e^{\theta \phi(t)}/[2\lambda(1\pm\e_t)]
$,
where $\eta_1\in(0,1)$ is defined in Lemma \ref{lem4.5}.
We have on $\mathcal B_2$ for such $z$,
\begin{equation}\label{zLambdaI}
0\leq z \lambda^{\pm}_t (\Ip_1+ \Ip_2) \leq \eta_1 a(t).
\end{equation}
This gives by \eqref{Fmoins} applied to $\Im_1\egloi F^-(h_t/2)$
for $t$ so large that $h_t/2\geq M$ and $a(t)\leq 1$,
\begin{eqnarray*}
    \left. E\left(e^{-z \lambda^{\pm}_t (\Ip_1+\Ip_2)\Im_1}\right|\Ip_1,\Ip_2\right)
& = &
    E\left(e^{-\gamma F^-(h_t/2)}\right)_{\big|\gamma=z \lambda^{\pm}_t (\Ip_1+\Ip_2)}
\\
& \leq &
    Z_t(z)
    +C_4 \max\big(e^{-\k h_t/2}, \big[z\lambda_t^\pm(\Ip_1+\Ip_2)\big]^{\frac{3}{2}}\big),
\end{eqnarray*}
on $\mathcal B_2$ for such $z$, thanks to the independence of $\Im_1$ and $(\Ip_1, \Ip_2)$ for the equality.
Therefore for large $t$, by \eqref{eqKi},
and since $0\leq Z_t(z)\leq 1$,
\begin{eqnarray*}
    S_{2,t}^{\pm}
& \leq &
    \int_0^{\frac{\eta_1 a(t) e^{\theta \phi(t)}}{t\lambda_t^\pm}}
    E\left( \un_{\mathcal B_2}Z_t^{2}(z)
    +\un_{\mathcal B_2} 3 C_4 \Big([z\lambda_t^\pm(\Ip_1+\Ip_2)]^{\frac{3}{2}}+e^{-\k h_t/2}\Big)
    \right)
    e^{-z}\textnormal{d}z\\
&&
    +
    \int_{\eta_1 a(t) e^{\theta \phi(t)}/(t\lambda_t^\pm)}^{ \infty}e^{-z} \textnormal{d}z.
\end{eqnarray*}
We notice that
$\eta_1 a(t) e^{\theta \phi(t)}/(t\lambda_t^\pm)\geq \phi(t)$
for large $t$ since $\theta>3\k/4$, and so
$
    \int_{\eta_1 a(t) e^{\theta \phi(t)}/(t\lambda_t^\pm)}^{\infty} e^{-z}\textnormal{d}z
\leq
e^{-\phi(t)}
$.
Moreover by \eqref{zLambdaI},
$
    \int_0^{\eta_1 a(t) e^{\theta \phi(t)}/(t\lambda^\pm_t)}
    E\left(\un_{\mathcal B_2} [z \lambda^{\pm}_t (\Ip_1+\Ip_2)]^{3/2}\right)e^{-z}\textnormal{d}z
\leq
    [\eta_1 a(t)]^{3/2}=o(e^{-\k\phi(t)})
$, and $e^{-\k h_t/2}=o(e^{-\k \phi(t)})$ since $\phi(t)=o(\log t)$.
So for large $t$  since $0<\k<1$,
\begin{equation}\label{InegS2}
    S_{2,t}^{\pm}
\leq
    \int_0^\infty
    E\left[ \un_{\mathcal B_2} Z_t^2(z)
    \right]
    e^{-z}\textnormal{d}z
    +o(e^{-\k\phi(t)}).
\end{equation}
Recall that for
 any random variable $Y\geq 0$, we have by Fubini,
\begin{align*}
    E\left(\left(1+Y\right)^{-2} \right)
& =
    \int_0^{ \infty}\textnormal{d}u
    \int_u^{ \infty}
    E\left(e^{-x(1+ Y)} \right)\textnormal{d}x.
\end{align*}
So, by independence of $\Ip_1$ and $\Ip_2$, we have for every $z> 0$,
\begin{equation}
    E\left(Z_t^2(z)\right)
=
    \int_0^{ \infty} \textnormal{d}u \int_u^{ \infty}  e^{-x}
    E\left[\exp\left( -\frac{ \rho_{\kappa}^\pm z x \lambda \Ip_1}{t}\right)\right]
    E\left[\exp\left( -\frac{ \rho_{\kappa}^\pm z x \lambda \Ip_2}{t}\right)\right]
    \textnormal{d}x,
\label{eqDoubleIntetEsp}
\end{equation}
where
$$
    \rho_{\kappa}^\pm
:=
    \frac{4 \left(1\pm \e_t\right)}{\kappa+1}.
$$
Recall that $t=e^{h_t} e^{\phi(t)}$,
$\Ip_1 \egloi F^+(h_t)$ and
$\Ip_2 \egloi G^+(h_t/2,h_t)$ by Proposition \ref{lemU1}.
So applying \eqref{Ff2} and \eqref{Gpf2},
we have whenever
$ \rho_{\kappa}^\pm z x\lambda /e^{\phi(t)}  \leq \eta_1$ and $h_t\geq M$,
\begin{align*}
&
    E\left(e^{ -\frac{ \rho_{\kappa}^\pm z x \lambda \Ip_1}{t}}\right)
\leq
    1-\frac{2\rho_{\kappa}^\pm z x\lambda}{(\k+1)e^{\phi(t)}}
    +C_4\max\left(e^{-\k h_t}, \left(\rho_{\kappa}^\pm z x\lambda/e^{\phi(t)}\right)^{3/2}\right),
    \\
&
    E\left(e^{ -\frac{ \rho_{\kappa}^\pm z x \lambda \Ip_2}{t}}\right)
\leq
    1-\frac{\Gamma(1-\kappa)}{\Gamma(1+\kappa)}\frac{ (2\rho_{\kappa}^\pm  z x \lambda)^{\kappa}    }
    {e^{{\kappa}\phi(t)} }
    +C_4\max\left(\left(\rho_{\kappa}^\pm z x\lambda/e^{\phi(t)}\right)^\k e^{-\k h_t/2},
\frac{\rho_{\kappa}^\pm z x\lambda}{e^{\phi(t)}}\right).
\end{align*}
So, \eqref{eqDoubleIntetEsp} gives  for large $t$,
\begin{align}
    E\left(Z_t^2(z)\right)
 \leq &
    \int_0^{+ \infty} \textnormal{d}u
        \int_u^{+ \infty} e^{-x}
        \left(1-\frac{\Gamma(1-\kappa)}{\Gamma(1+\kappa) }\frac{(2\rho_{\kappa}^\pm  z x\lambda )^{\kappa}}
        { e^{{\kappa}\phi(t)} }\right)
        \un_{\rho_{\kappa}^\pm z x \lambda/e^{\phi(t)}  \leq \eta_1}
        \textnormal{d}x
\nonumber\\
 & +
    C_+\int_0^{+ \infty} \textnormal{d}u
        \int_u^{+ \infty}  e^{-x}
        \left(\left(\frac{\rho_{\kappa}^\pm z x\lambda}{e^{\phi(t)}}\right)^\k e^{-\k h_t/2}+
        \frac{\rho_{\kappa}^\pm z x\lambda}{e^{\phi(t)}}+\frac{1}{e^{\k h_t/2}}\right)
        \un_{\rho_{\kappa}^\pm z x \lambda/e^{\phi(t)}  \leq \eta_1}
        \textnormal{d}x
\nonumber\\
&
    + \int_0^{+ \infty} \textnormal{d}u
        \int_u^{+ \infty} e^{-x}
        \un_{\rho_{\kappa}^\pm z x\lambda /e^{\phi(t)} > \eta_1}\dd x.
\label{EsperanceZ2Bis}
\end{align}
Now, notice that
$
\int_0^{+ \infty} e^{-z}\textnormal{d}z  \int_0^{+ \infty} \textnormal{d}u \int_u^{+ \infty}
 e^{-x} {\rho_{\kappa}^\pm z x\lambda}\textnormal{d}x/{e^{\phi(t)}}
= 2 \rho_{\kappa}^\pm \lambda/e^{\phi(t)}
$.
Moreover, we have
$
\int_0^{+ \infty} e^{-z}\textnormal{d}z  \int_0^{+ \infty} \textnormal{d}u \int_u^{+ \infty}
 e^{-x} \left[{\rho_{\kappa}^\pm z x\lambda}/{e^{\phi(t)}}\right]^\k \textnormal{d}x
= O(e^{-\k\phi(t)})
$,
and furthermore
$
    \un_{\rho_{\kappa}^\pm z x\lambda /e^{\phi(t)} > \eta_1}
\leq
    \rho_{\kappa}^\pm z x\lambda /\left[\eta_1 e^{\phi(t)}\right]
$,
so that
$
\int_0^{+ \infty}  \textnormal{d}z e^{-z} $ $\int_0^{+ \infty} \textnormal{d}u \int_u^{+ \infty} e^{-x}
\un_{\rho_{\kappa}^\pm z x\lambda /e^{\phi(t)} > \eta_1} \dd x
=O(e^{-\phi(t)})
$,
whereas
$
\int_0^{+ \infty}  \textnormal{d}z e^{-z} $ $\int_0^{+ \infty} \textnormal{d}u \int_u^{+ \infty} e^{-x}
\dd x
=1
$.
This together with \eqref{EsperanceZ2Bis} and $\phi(t)=o(h_t)$ gives
\begin{align}
&
    \int_0^{+ \infty} E\left(Z_t^2(z)\right)e^{-z} \textnormal{d}z
\nonumber \\
&
    \leq  \int_0^{+ \infty} \textnormal{d}z  e^{-z}  \int_0^{+ \infty} \textnormal{d}u
    \int_u^{+ \infty} e^{-x}
    \left(1-\frac{\Gamma(1-\kappa)}{\Gamma(1+\kappa) }
    \frac{(2\rho_{\kappa}^\pm  z x\lambda )^{\kappa}    }{ e^{{\kappa}\phi(t)} }\right)
    \un_{\rho_{\kappa}^\pm z x \lambda/e^{\phi(t)}  \leq \eta_1}\textnormal{d}x
    + O(e^{- \phi(t)})
\nonumber \\
&=   1
    -\frac{\Gamma(1-\kappa)}{\Gamma(1+\kappa) }
    \frac{ 8^{\kappa}\lambda^{\kappa}}{(1+\kappa)^{\kappa}e^{\k\phi(t)}}
    (1\pm \e_t)^{\kappa} \int_0^{+ \infty}  e^{-z}z^{\kappa}\textnormal{d}z
    \int_0^{+ \infty} e^{-x} x^{\kappa+1}\textnormal{d}x
    \nonumber \\
&   \ \   +
    \frac{\Gamma(1-\kappa)(2\lambda\rho_{\kappa}^\pm)^\k}{\Gamma(1+\kappa)e^{{\kappa}\phi(t)} }
    \int_0^{+ \infty} \textnormal{d}z  e^{-z}  \int_0^{+ \infty} \textnormal{d}u
    \int_u^{+ \infty}  e^{-x}
    (  z x )^{\kappa}
    \un_{\rho_{\kappa}^\pm z x \lambda/e^{\phi(t)}  >\eta_1}\textnormal{d}x
    + O(e^{- \phi(t)}) \label{IntegraleEta1AEstimer}\\
&= 1-\frac{\Gamma(1-\kappa)\Gamma(\k+2) 8^{\kappa}\lambda^{\kappa} }{(1+ \kappa)^{\kappa}e^{\kappa\phi(t)}   }
    (1\pm \e_t)^{\kappa}+ o(e^{- \k\phi(t)}),
    \label{eq204}
\end{align}
since, by the dominated convergence theorem, the integral in Line \eqref{IntegraleEta1AEstimer} goes to $0$ as $t\to+\infty$
and then Line $\eqref{IntegraleEta1AEstimer}=o(e^{-\k\phi(t)})$.
Combining equations \eqref{InegS2}, \eqref{eqDefS3}, then \eqref{eq204} and $\lim_{t\to+\infty}\e_t=0$, we get
\begin{equation}\label{InegS2PlusS3}
    S_{2,t}^{\pm}+S_{3,t}^{\pm}
\leq
        \int_0^{+ \infty} E\left(Z_t^2(z)\right)e^{-z} \textnormal{d}z
        +o(e^{-\k\phi(t)})
\leq
    1-C_0 8^{\k}\lambda^\k e^{-\k \phi(t)}+o(e^{-\k\phi(t)}),
\end{equation}
where $C_0={\Gamma(1-\kappa)\Gamma(\k+2)  }/{(1+ \kappa)^{\kappa}}$ as defined in
Proposition \ref{proplap}.
We prove similarly that $S_{2,t}^{\pm}+S_{3,t}^{\pm} \geq 1-C_0 8^{\k}\lambda^\k e^{-\k \phi(t)}+o(e^{-\k\phi(t)})$.
This proves \eqref{eqConclusionPartA} and then the lemma.
\hfill$\Box$


We now turn to the case $\Ip_1+\Ip_2$ is "big",
as is the case on $\mathcal{B}_1$. More precisely, we prove

\medskip

\begin{lemma}\label{LemmaPartB} We have, with $\Upsilon_0$ as defined in Proposition \ref{proplap},
\begin{equation}\label{eqLaplacePartB}
    (S_{1,t}^{\pm}-S_{3,t}^{\pm})
\sim_{t\to+\infty}
    (8\lambda)^\k \Upsilon_0e^{-\k \phi(t)}.
\end{equation}
\end{lemma}

\medskip

\noindent {\bf Proof:}
We introduce $0<\e<1/2$.

\noindent {\bf Step 1: } Approximation of $\Ip_1+\Ip_2$.
Since
$\theta<\k$,
$\e e^{  (1-\theta) \phi(t)}\geq \e(\log t)^2 \geq 8 h_t/\k$
for $t$ large enough so that $\phi(t)\geq 2(\log\log t)/(1-\k)$.
This gives as $t\to+\infty$,
\begin{equation}\label{InegProba1+}
    P\big(\Ip_1\geq \e  t e^{ - \theta \phi(t)}
    \big)
\leq
    P\big(\tau^R(h_t)\geq \e  e^{(1-\theta) \phi(t)}\big)
\leq
    P\big(\tau^R(h_t) \geq  8 h_t/\k\big)
\leq
    C_+ e^{- \kappa h_t /2\sqrt{2}}
=
    o\big(e^{-\phi(t)}\big)
\end{equation}
since $\Ip_1 \egloi F^+(h_t)\leq e^{h_t}\tau^R(h_t)$ by Proposition \ref{lemU1}, $\phi(t)=o(\log t)$,
and thanks to \eqref{bessel4}.
Moreover,
$\Ip_2 \egloi G^+(h_t/2,h_t)$ and  $\Im_1$, $\Im_2$, $\Ip_1$ and $\Ip_2$ are independent by Proposition \ref{lemU1}.
So, possibly on an enlarged probability space, there exists a random variable $\tilde A_\infty$ equal in law to $A_\infty$
($A_\infty$ being defined in \eqref{eqDefA} and $G^+$ in \eqref{eqDefF+G+})
and such that $\tilde A_\infty$,  $\Im_1$ ,$\Im_2$, $\Ip_1$ and $\Ip_2$ are independent under $P$.
We now introduce $\widehat A_\infty:=e^{-h_t}\Ip_2+e^{-h_t/2}\tilde A_\infty$. By Markov property,
\begin{equation}\label{eqDefAinftyChapeauEtc}
    \Ip_2= e^{h_t}\widehat A_\infty -e^{h_t/2}\tilde A_\infty,
\qquad
    \widehat A_\infty\egloi \tilde A_\infty\egloi A_\infty,
\end{equation}
and $\widehat A_\infty$, $\Im_1$, $\Im_2$ and $\Ip_1$ are independent.
We have
\begin{equation}\label{InegTildeAinfty}
    P\big(e^{h_t/2} \tilde A_\infty\geq \e  t e^{ - \theta \phi(t)}\big)
\leq
    P\big(A_\infty\geq t^{1/3}\big)
\leq
    C t^{-\k/3}
=
    o\big(e^{-\phi(t)}\big)
\end{equation}
as $t\to+\infty$,
since $\phi(t)=o(\log t)$ and
$P(A_{\infty}\geq y)\leq C y^{-\k}$ as explained before \eqref{eqProbaE2Lem35}.
Now, we have on
$
\mathcal B_1
\cap
\{\Ip_1< \e  t e^{ - \theta \phi(t)}\}
\cap
\{e^{h_t/2}\tilde A_\infty< \e  t e^{ - \theta \phi(t)}\}$,
$$
    t e^{ - \theta \phi(t)}
<
    \Ip_1+\Ip_2
=
    \Ip_1+e^{h_t}\widehat A_\infty -e^{h_t/2}\tilde A_\infty
<
    e^{h_t}\widehat A_\infty+\e t e^{ - \theta \phi(t)}.
$$
This yields
$
    e^{h_t}\widehat A_\infty
\geq
    (1-\e)t e^{ - \theta \phi(t)}
$,
and then
$\Ip_1+\Ip_2
\leq (1+\frac{\e}{1-\e})e^{h_t}\widehat A_\infty
\leq (1+2\e)e^{h_t}\widehat A_\infty$.
Similarly,
$
    \Ip_1+\Ip_2
\geq
    e^{h_t}\widehat A_\infty-e^{h_t/2}\tilde A_\infty
\geq
    e^{h_t}\widehat A_\infty-\e t e^{ - \theta \phi(t)}
\geq
    (1-\frac{\e}{1-\e})e^{h_t}\widehat A_\infty
\geq
    (1-2\e)e^{h_t}\widehat A_\infty
$
on the same event.
Consequently, replacing $\e$ by $\e/2$, we get for $0<\e<1$ for large $t$,
\begin{eqnarray}
&&
    P\big(\mathcal B_1\cap\{(1-\e)e^{h_t}\widehat A_\infty\leq \Ip_1+\Ip_2\leq (1+\e)e^{h_t}\widehat  A_\infty\}^c\big)
\nonumber\\
& \leq &
    P\big(\Ip_1\geq (\e/2)  t e^{ - \theta \phi(t)}\big)
    +
    P\big(e^{h_t/2} \tilde A_\infty \geq (\e/2)  t e^{ - \theta \phi(t)}\big)
\leq
    e^{-\phi(t)},
\label{eqApproxI1PlusI2}
\end{eqnarray}
by \eqref{InegProba1+} and \eqref{InegTildeAinfty}.
Similarly, on
$
    \{(1+\e)e^{h_t}\widehat A_\infty\geq t e^{-\theta\phi(t)}\}
\cap
    \{\Ip_1+\Ip_2> (1+\e)e^{h_t}\widehat  A_\infty\}
$,
we have by \eqref{eqDefAinftyChapeauEtc},
$$
    \e t e^{-\theta\phi(t)}/(1+\e)
\leq
    \e e^{h_t}\widehat A_\infty
\leq
    \Ip_1+\Ip_2-e^{h_t}\widehat A_\infty
\leq
    \Ip_1.
$$
Moreover, we have on
$
    \{(1+\e)e^{h_t}\widehat A_\infty\geq t e^{-\theta\phi(t)}\}
\cap
    \{\Ip_1+\Ip_2< (1-\e)e^{h_t}\widehat  A_\infty\}
$, again by \eqref{eqDefAinftyChapeauEtc},
$$
    \e t e^{-\theta\phi(t)}/(1+\e)
\leq
    \e e^{h_t}\widehat A_\infty
\leq
    e^{h_t}\widehat A_\infty-\Ip_1-\Ip_2
\leq
    e^{h_t/2}\tilde A_\infty.
$$
Consequently for large $t$, the last inequality being obtained as in
\eqref{InegProba1+} and \eqref{InegTildeAinfty},
\begin{eqnarray}
&&
    P\big(\{(1+\e)e^{h_t}\widehat A_\infty\geq t e^{-\theta\phi(t)}\}
        \cap\{(1-\e)e^{h_t}\widehat A_\infty\leq \Ip_1+\Ip_2\leq (1+\e)e^{h_t}\widehat  A_\infty\}^c\big)
\nonumber\\
& \leq &
    P\big(\Ip_1\geq [\e/(1+\e)]  t e^{ - \theta \phi(t)}\big)
    +
    P\big(e^{h_t/2} \tilde A_\infty\geq [\e/(1+\e)]  t e^{ - \theta \phi(t)}\big)
\leq
    e^{-\phi(t)}.
\label{eqApproxI1PlusI2Bis}
\end{eqnarray}


\noindent {\bf Step 2:} Simplification.
Thanks to \eqref{eqKiIntermediaire}, \eqref{eqDefS3} first
and then  the definition \eqref{eqDefBi} of $\mathcal B_1$, \eqref{eqApproxI1PlusI2}
and $|e^{-(\dots)}-Z_t^2(u))|\leq 1$, we get
\begin{eqnarray}
S_{1,t}^{\pm}-S_{3,t}^{\pm}
& = &
    \int_0^\infty
    E\left[
        \left(e^{-\lambda_t^\pm u (\Ip_1+\Ip_2)(\Im_1+\Im_2)} -Z_t^2(u) \right)
        \un_{\mathcal B_1}
    \right]e^{-u}\dd u
\nonumber\\
& \leq &
    \int_0^\infty
    E\left[\left(e^{-\lambda_t^\pm u (\Ip_1+\Ip_2)(\Im_1+\Im_2)}
    -Z_t^2(u) \right)
    \un_{\{\Ip_1+\Ip_2\geq t e^{-\theta\phi(t)}\}
    }\right.\nonumber\\
&&
    \left. \qquad\quad\ \times\un_{\{(1-\e)e^{h_t}\widehat A_\infty\leq \Ip_1+\Ip_2\leq (1+\e)e^{h_t}\widehat A_\infty\}}\right]e^{-u} \dd u
    +e^{-\phi(t)}.\nonumber
\end{eqnarray}
Consequently, using $0\leq Z_{t,\e,\infty}(u)\leq Z_t(u)$ in the expectation for $Z_{t,\e,\infty}$ defined below,
\begin{eqnarray}
    S_{1,t}^{\pm}-S_{3,t}^{\pm}
& \leq &
        \int_0^\infty
    E\left[\left(e^{-\lambda_t^\pm u (1-\e)e^{h_t}\widehat A_\infty(\Im_1+\Im_2)}
            -Z_{t,\e,\infty}^2(u)
            \right)
    \un_{\{\Ip_1+\Ip_2\geq t e^{-\theta\phi(t)}\}
    }\right.
\nonumber\\
&&
    \left. \qquad\quad\ \times\un_{\{(1-\e)e^{h_t}\widehat A_\infty\leq \Ip_1+\Ip_2\leq (1+\e)e^{h_t}\widehat A_\infty\}}\right]e^{-u} \dd u
    +e^{-\phi(t)}
\nonumber\\
& = &
    S_{4,t,\e}^{\pm}+S_{5,t,\e}^{\pm}
    +e^{-\phi(t)},\label{eqK1K4}
\end{eqnarray}
where
\begin{eqnarray*}
    Z_{t,\e,\infty}(u)
&:=&
    \left[1+{2 \lambda^\pm_t u (1+\e)e^{h_t}\widehat A_\infty}/({\kappa+1})\right]^{-1},
\nonumber\\
    S_{4,t,\e}^{\pm}
& := &
    \int_0^\infty
    e^{-u}
    E\left[\left(e^{-\lambda_t^\pm u (1-\e)e^{h_t}\widehat A_\infty(\Im_1+\Im_2)}
    - Z_{t,\e,\infty}^2(u)\right)
    \un_{\{(1+\e)e^{h_t}\widehat A_\infty\geq t e^{-\theta\phi(t)}\}}\right.\nonumber\\
&&
    \left.
    \qquad\qquad\qquad
    \un_{\{(1-\e)e^{h_t}\widehat A_\infty\leq \Ip_1+\Ip_2\leq (1+\e)e^{h_t}\widehat A_\infty\}}\right] \dd u,
\nonumber\\
    S_{5,t,\e}^{\pm}
& := &
    \int_0^\infty
    e^{-u}
    E\left[\left(e^{-\lambda_t^\pm u (1-\e)e^{h_t}\widehat A_\infty(\Im_1+\Im_2)}
    -Z_{t,\e\infty}^2(u)\right)
    \right.
\nonumber\\
&&
    \left.
    [\un_{\{\Ip_1+\Ip_2\geq t e^{-\theta\phi(t)}\}}-\un_{\{(1+\e)e^{h_t}\widehat A_\infty\geq t e^{-\theta\phi(t)}\}}]
     \un_{\{(1-\e)e^{h_t}\widehat A_\infty\leq \Ip_1+\Ip_2\leq (1+\e)e^{h_t}\widehat A_\infty\}}\right] \dd u.
\nonumber
\end{eqnarray*}
We will prove in our Step 5 that $S_{5,t,\e}^{\pm}$ is negligible.
So, we start with $S_{4,t,\e}^{\pm}$.
Using \eqref{eqApproxI1PlusI2Bis}
and $|e^{-(\dots)}-Z_{t,\e,\infty}(u)^2)|\leq 1$
leads to
\begin{equation}
\left|S_{4,t,\e}^{\pm} -
    \int_0^\infty
    E\left[\left(e^{-\lambda_t^\pm u (1-\e)e^{h_t}\widehat A_\infty(\Im_1+\Im_2)}
    -Z_{t,\e,\infty}^{2}(u)\right)
    \un_{\{(1+\e)e^{h_t}\widehat A_\infty\geq t e^{-\theta\phi(t)}\}}\right] \frac{\dd u}{e^u}
\right|
\leq
    \frac{1}{e^{\phi(t)}}.
    \label{eqTranfoIntegraleK1K4}
\end{equation}
Recall that by \eqref{eqDefAinftyChapeauEtc} and Proposition \ref{lemU1},
$\Im_1$, $\Im_2$ and $\widehat A_\infty $ are independent, and $\Im_1\egloi \Im_2$.
So,  the integral in \eqref{eqTranfoIntegraleK1K4} can be written as
\begin{align}
  &
    \int_0^\infty
        E\left[\left(
                    E\left(e^{-\lambda_t^\pm u (1-\e)e^{h_t}\widehat A_\infty \Im_1}| \widehat A_\infty\right)^2
                    -Z_{t,\e\infty}^{2}(u)
               \right)
            \un_{\widehat A_\infty \geq  \frac{e^{(1-\theta)\phi(t)}}{1+\e}}
        \right]
    e^{-u} \dd u
\nonumber\\
 \leq &
    \int_0^\infty
        E\left[\left(
                    E\left(e^{-4\lambda (1-\e)^2 u e^{-\phi(t)}\gamma_\k^{-1}\Im_1}| \gamma_\k\right)^2
                    -
                    \bigg[1+\frac{8\lambda(1+\e)^2 u}{(\k+1)e^{\phi(t)}\gamma_\k}\bigg]^{-2}
               \right)
            \un_{\gamma_\k\leq \frac{2(1+\e)}{e^{(1-\theta)\phi(t)}}}
    \right]
    \frac{\dd u}{e^u},
\label{IntegrabiliteFacile}
\end{align}
where $\gamma_\k:=2/\widehat A_\infty$,
for $t$ large enough so that
$2\lambda (1-\e)/t\leq \lambda_t^\pm\leq 2\lambda(1+\e)/t$.
Since $\gamma_\k$ has density $x^{\k-1} e^{-x}\un_{\R_+}(x)/\Gamma(\k)$
by Fact \ref{FactDufresnes} (Dufresne),
and is, as $\widehat A_\infty$,  independent of $\Im_1$, the RHS of \eqref{IntegrabiliteFacile} is equal to
\begin{equation}
    \int_0^\infty e^{-u} \dd u
    \int_0^{\frac{2(1+\e)}{e^{(1-\theta)\phi(t)}}}
    \bigg[
        E\left(e^{-4\lambda (1-\e)^2 u e^{-\phi(t)} x^{-1}\Im_1}\right)^2
        -\left(1+\frac{8 \lambda(1+\e)^2 u  }{(\kappa+1)e^{\phi(t)}x}\right)^{-2}
    \bigg]
    \frac{x^{\k-1}e^{-x}}{\Gamma(\k)}\dd x . \nonumber\\
\end{equation}
With the change of variables $y=8 u \lambda e^{-\phi(t)}x^{-1}$, this is equal to
$(8\lambda )^\k e^{-\k \phi(t)}\Upsilon_{t,\epsilon}$, with
$
    \Upsilon_{t,\epsilon}
:=
    \int_0^\infty\int_0^\infty f_{t,\e}(u,y)\dd y \dd u
$
and
$$
    f_{t,\e}(u,y)
:=
    \un_{\{ y\geq \frac{ 4\lambda u e^{-\theta\phi(t)}}{1+\e}\}}
    \frac{u^\k e^{-u}}{\Gamma(\k)}
    \bigg[
        E\left(e^{- (1-\e)^2y \Im_1/2}\right)^2
        -\left(1+\frac{ (1+\e)^2 y  }{\kappa+1}\right)^{-2}
    \bigg]
    \frac{e^{-8 \lambda  u y^{-1} e^{-\phi(t)}}}{y^{\k+1}}
    .
$$
Consequently by \eqref{eqTranfoIntegraleK1K4},
\begin{equation}\label{eqStep2}
    S_{4,t,\e}^{\pm}
\leq
    (8\lambda )^\k e^{-\k \phi(t)}\Upsilon_{t,\epsilon}+e^{-\phi(t)}.
\end{equation}

\smallskip
\noindent {\bf Step 3:} pointwise convergence. \\
Notice that thanks to
\eqref{eqLimiteLaplaceF-},
$
    \lim_{x\to+\infty} E\big(e^{-\gamma F^-(x)}\big)
=
    \frac{ (2\gamma)^{{\k/2}}}{\k\Gamma(\k)I_\k(2\sqrt{2\gamma})}
$
for $\gamma>0$,
and recall that $\Im_1\egloi F^-(h_t/2)$ by Proposition \ref{lemU1}.
Hence, for every $0<\e<1$, $u>0$ and $y>0$,
\begin{equation}\label{eqDeffepsilon}
    f_{t,\e}(u,y)
\to_{t\to+\infty}
    f_\e(u,y)
:=
    \frac{u^\k e^{-u}}{\Gamma(\k)}
    \bigg[
        \frac{ (1-\e)^{2\k} y^{\k}}{[\Gamma(\k+1)I_\k(2(1-\e)\sqrt{ y})]^2}
        -\left(1+\frac{(1+\e)^2 y}{\kappa+1}\right)^{-2}
    \bigg]
    y^{-\k-1}.
\end{equation}

\noindent{\bf Step 4:} dominated convergence.
We notice that
$f_{t,\e}(u,y)=a_{t,\e}(u,y)+b_{t,\e}(u,y)$ and $f_\e(u,y)=a_\e(u,y)+b_\e(u,y)$
for every $0<\e<1$, $u>0$ and $y>0$, where
\begin{align*}
    a_{t,\e}(u,y)
& :=
    \un_{\{ y\geq \frac{ 4\lambda u e^{-\theta\phi(t)}}{1+\e}\}}
    \frac{u^\k e^{-u}}{\Gamma(\k)}
    \bigg[
        E\left(e^{- (1-\e)^2y \Im_1/2}\right)^2
        -\left(1+\frac{ (1-\e)^2 y  }{\kappa+1}\right)^{-2}
    \bigg]
    \frac{e^{-8\lambda  u y^{-1} e^{-\phi(t)}}}{y^{\k+1}},
\\
    b_{t,\e}(u,y)
& :=
    \un_{\{ y\geq \frac{ 4\lambda u e^{-\theta\phi(t)}}{1+\e}\}}
    \frac{u^\k e^{-u}}{\Gamma(\k)}
    \bigg[
        \left(1+\frac{ (1-\e)^2 y  }{\kappa+1}\right)^{-2}
        -\left(1+\frac{ (1+\e)^2 y  }{\kappa+1}\right)^{-2}
    \bigg]
    \frac{e^{-8\lambda  u y^{-1} e^{-\phi(t)}}}{y^{\k+1}},
\end{align*}
and their pointwise limits on $(\R_+^*)^2$ as $t\to+\infty$ are respectively
\begin{eqnarray}
    a_\e(u,y)
& := &
    \frac{u^\k e^{-u}}{\Gamma(\k)}
    \bigg[
        \frac{ (1-\e)^{2\k} y^{\k}}{[\Gamma(\k+1)I_\k(2(1-\e)\sqrt{ y})]^2}
        -\left(1+\frac{(1-\e)^2 y}{\kappa+1}\right)^{-2}
    \bigg]
    y^{-\k-1},
\nonumber
\\
    b_\e(u,y)
& := &
    \frac{u^\k e^{-u}}{\Gamma(\k)}
    \bigg[
        \left(1+\frac{(1-\e)^2 y}{\kappa+1}\right)^{-2}
        -\left(1+\frac{(1+\e)^2 y}{\kappa+1}\right)^{-2}
    \bigg]
    y^{-\k-1}.
\label{eqDefBEpsilon}
\end{eqnarray}
Since $\forall x> 0,\ I_\k(x)> \frac{(x/2)^{\k}}{\Gamma(\k+1)}+\frac{(x/2)^{\k+2}}{\Gamma(\k+2)}$
due to the series expansion of $I_\k$ (see e.g. \cite{BorodinSalminem} p. 638), we get
$
 \frac{ \gamma^{{\k}}}{[\Gamma(\k+1)I_\k(2\sqrt{\gamma})]^2}
   \
   -\left(1+\frac{\gamma}{\kappa+1}\right)^{-2}
<
0
$ for every $\gamma>0$.
One consequence of this inequality is that $\Upsilon_0<0$, $\Upsilon_0$ being defined in \eqref{eqDefUpsilon0},
and another one is that
$\forall (u,y)\in(\R_+^*)^2,\ a_\e(u,y)< 0$.
This, and the fact that $x\mapsto E \left( e^{- \gamma F^-(x)}  \right)$ is nonincreasing for $\gamma\geq 0$, lead to
\begin{equation}\label{InegENcadrementIntegrabilitef}
\forall (u,y)\in(\R_+^*)^2,\quad
    a_\e(u,y)
\leq
    a_{t,\e}(u,y)
\leq
    \frac{u^\k e^{-u}}{\Gamma(\k)}
    \bigg[
        1-\left(1+\frac{ (1-\e)^2 y  }{\kappa+1}\right)^{-2}
    \bigg]
    y^{-\k-1}
=:
    h_\e(u,y).
\end{equation}
Hence, $|a_{t,\e}(u,y)|\leq |a_\e(u,y)|+|h_\e(u,y)|$ and $|b_{t,\e}(u,y)|\leq |b_\e(u,y)|$.
for every $(u,y)\in\R_+^2$
and $0<\e<1$. Moreover, since $0<\k<1$,  $|h_\e|$, $|b_\e|$ and $|a_\e|$ have finite integrals over $\R_+^2$
(notice for example that $e^u u^{-\k}a(u,y)=O(y^{-\k})$ as $y\to 0$,
since
$I_\k(x) =(x/2)^\k/\Gamma(\k+1)+(x/2)^{\k+2}/\Gamma(\k+2)+o(x^{\k+2})$ as $x\to 0$ by \eqref{eqDLBesselI} below,
and $-y^{-\k-1}u^\k e^{-u}/\Gamma(\k)\leq a_\e(u,y)  < 0$ as $y\to+\infty$).
Thus, by the dominated convergence theorem,
\begin{equation}\label{eqDefUpsilonEpsilon}
    \lim_{t \rightarrow + \infty}\Upsilon_{t,\epsilon}
=
    \int_0^\infty\int_0^\infty f_\e(u,y)\dd y \dd u
=:
    \Upsilon_{\epsilon},
\qquad
    0<\e<1.
\end{equation}
Hence applying \eqref{eqStep2},
\begin{equation}\label{eqStep4}
    \limsup_{t\to+\infty}\big(S_{4,t,\e}^{\pm}e^{\k \phi(t)}\lambda^{-\k} \big)
\leq
    8^{\k}\Upsilon_{\epsilon},
\qquad
0<\e<1.
\end{equation}

\noindent {\bf Step 5:} Let $0<\e<1$; we now prove that $S_{5,t,\e}^{\pm}$ is negligible.
First, we have
\begin{eqnarray}
S_{5,t,\e}^{\pm}
& = &
    \int_0^\infty e^{-u}
    E\left[\left(Z_{t,\e,\infty}^{2}(u)-e^{-\lambda_t^\pm u (1-\e)e^{h_t} \widehat A_\infty(\Im_1+\Im_2)}
    \right)\right.\nonumber\\
&&    \left.
    \un_{\{\Ip_1+\Ip_2< t e^{-\theta \phi(t)}\}}\un_{\{(1+\e)e^{h_t} \widehat A_\infty\geq t e^{-\theta\phi(t)}\}}
    \un_{\{(1-\e)e^{h_t} \widehat A_\infty\leq \Ip_1+\Ip_2\leq (1+\e)e^{h_t} \widehat A_\infty\}}
    \right] \dd u,\nonumber
\end{eqnarray}
where we used
$
    (\un_{E_1}-\un_{E_2})\un_{E_3}
=
    -\un_{\overline{E_1}}\un_{E_2}\un_{E_3}
+
    \un_{E_1}\un_{\overline{E_2}}\un_{E_3}
$,
for events $E_i$
with $\un_{E_1}\un_{\overline{E_2}}\un_{E_3}=0$ in our case.
Conditioning by $\sigma\big(\widehat A_\infty, \Ip_1,\Ip_2\big)$, which is independent of $\Im_1$ and $\Im_2$
gives since $\Im_1$ and $\Im_2$ are independent,
\begin{eqnarray}
S_{5,t,\e}^{\pm}
&= &
    \int_0^\infty e^{-u}
    E\left[\left(
                Z_{t,\e,\infty}^{2}(u)
                -E\left(e^{-\lambda_t^\pm u (1-\e)e^{h_t} \widehat A_\infty\Im_1}\Big| \widehat A_\infty\right)^2
            \right)
    \right.\nonumber\\
&&    \left.
    \un_{\{\Ip_1+\Ip_2< t e^{-\theta \phi(t)}\}}\un_{\{(1+\e)e^{h_t} \widehat A_\infty\geq t e^{-\theta\phi(t)}\}}
    \un_{\{(1-\e)e^{h_t} \widehat A_\infty\leq \Ip_1+\Ip_2\leq (1+\e)e^{h_t} \widehat A_\infty\}}
    \right] \dd u\nonumber
\\
& \leq &
    \int_0^\infty e^{-u}
    E\bigg[\bigg(
               \bigg[1+\frac{4 \lambda u (1-\e^2) e^{-\phi(t)}\widehat A_\infty}{\kappa+1}\bigg]^{-2}
               -E\bigg(e^{-2\lambda u (1-\e^2) e^{-\phi(t)} \widehat A_\infty\Im_1}\Big| \widehat A_\infty\bigg)^2
            \bigg)
\nonumber\\
&&
    \un_{\{\Ip_1+\Ip_2< t e^{-\theta \phi(t)}\}}\un_{\{(1+\e)e^{h_t} \widehat A_\infty\geq t e^{-\theta\phi(t)}\}}
    \un_{\{(1-\e)e^{h_t} \widehat A_\infty\leq \Ip_1+\Ip_2\leq (1+\e)e^{h_t} \widehat A_\infty\}}
    \bigg] \dd u
\label{InegS5Numero1}
\end{eqnarray}
for large $t$
where we used $2\lambda (1-\e)/t\leq \lambda_t^\pm\leq 2\lambda(1+\e)/t$.
Consequently,
\begin{eqnarray}
S_{5,t,\e}^{\pm}
& \leq &
    \int_0^\infty e^{-u}
    E\bigg[\bigg|
                \bigg[1+\frac{4 \lambda u (1-\e^2) e^{-\phi(t)}\widehat A_\infty}{\kappa+1}\bigg]^{-2}
                -
                E\bigg(e^{-2\lambda u (1-\e^2) e^{-\phi(t)} \widehat A_\infty\Im_1}\Big| \widehat A_\infty\bigg)^2
            \bigg|
\nonumber\\
&&    \left.
    \un_{\{(1-\e)e^{h_t} \widehat A_\infty\leq t e^{-\theta\phi(t)} \leq (1+\e)e^{h_t} \widehat A_\infty\}}
    \right]
    \dd u.
    \label{eqK6Modif1bis}
\end{eqnarray}
As in Step 2 after \eqref{IntegrabiliteFacile}, since $\gamma_\k=2/\widehat A_\infty$ has density $x^{\k-1} e^{-x}\un_{\R_+}(x)/\Gamma(\k)$,
and $\widehat A_\infty$ and then $\gamma_\k$ are independent of $\Im_1$, the right hand side of \eqref{eqK6Modif1bis} is equal to
\begin{equation}
    \int_0^\infty e^{-u}
    \int_{\frac{2(1-\e)}{e^{(1-\theta)\phi(t)}}}^{\frac{2(1+\e)}{e^{(1-\theta)\phi(t)}}}
    \bigg|
                \bigg[1+\frac{8 \lambda u (1-\e^2) }{(\kappa+1)e^{\phi(t)}x}\bigg]^{-2}
                -
                E\bigg(e^{-4\lambda u (1-\e^2) e^{-\phi(t)} x^{-1}\Im_1}
                \bigg)^2
    \bigg|
    \frac{x^{\k-1}e^{-x}}{\Gamma(\k)} \dd x
    \dd u.
\end{equation}
With the change of variables $y=8 u \lambda e^{-\phi(t)}x^{-1}$,
this is equal to
\begin{eqnarray}
&  &
    \frac{(8\lambda)^\k}{ e^{\k\phi(t)}}
    \int_0^\infty \frac{u^\k e^{-u}}{\Gamma(\k)}
    \int_{\frac{4 u \lambda }{(1+\e)e^{\theta\phi(t)}}}
    ^{\frac{4 u \lambda }{(1-\e)e^{\theta\phi(t)}}}
    \left|
                E\left(e^{-y(1-\e^2) \Im_1/2}
                \right)^2
                -\left[1+\frac{(1-\e^2 )y}{\kappa+1}\right]^{-2}
    \right|
    \frac{\dd y \dd u}{y^{\k+1}e^{8 u \lambda e^{-\phi(t)}y^{-1}}}.
\nonumber
\end{eqnarray}
So by \eqref{eqK6Modif1bis},
\begin{equation}\label{InegS5Leq}
    S_{5,t,\e}^{\pm}
\leq
    \frac{ (8\lambda )^\k}{ e^{\k \phi(t)}}\int_0^\infty
    \int_0^\infty
    \big|F_{t,(1-\e^2), (1-\e^2), \e}\big(u,y\big)\big|
    \dd y \dd u
\end{equation}
where for $z_1\geq 0$ and $z_2\geq 0$,
$$
    F_{t,z_1,z_2,\e}(u,y)
:=
    \un_{\{\frac{ 4\lambda u e^{-\theta\phi(t)}}{1+\e}\leq y\leq \frac{ 4\lambda u e^{-\theta\phi(t)}}{1-\e}\}}
    \frac{u^\k e^{-u}}{\Gamma(\k)}
    \bigg[
        E\left(e^{- z_1 y \Im_1/2}\right)^2
        -\left(1+\frac{ z_2 y  }{\kappa+1}\right)^{-2}
    \bigg]
    \frac{1}{y^{\k+1}}.
$$
Using the same method, we get for large $t$,
$$
    S_{5,t,\e}^{\pm}
\geq
    -\frac{ (8\lambda )^\k}{ e^{\k \phi(t)}}\int_0^\infty
    \int_0^\infty
    \big|F_{t, (1-\e)^2,(1+\e)^2,\e}(u,y)\big|
    \dd y \dd u.
$$
This together with \eqref{InegS5Leq} gives
\begin{equation}\label{InegS5ValeurAbsolue}
    |S_{5,t,\e}^{\pm}|
    \frac{ e^{\k \phi(t)}}{ (8\lambda )^\k}
\leq
    \int_0^\infty\int_0^\infty
        \big(\big|F_{t,(1-\e^2), (1-\e^2),\e}(u,y)\big|
        +
        \big|F_{t,(1-\e)^2,(1+\e)^2, \e}(u,y)\big|\big)
    \dd y \dd u.
\end{equation}
We have for every $z_1\geq 0$, $z_2\geq 0$,
\begin{equation}
    F_{t,z_1,z_2,\e} (u,y)
 \leq
    \frac{u^\k e^{-u}}{\Gamma(\k)}
    \bigg[
        1
        -\left(1+\frac{ z_2 y  }{\kappa+1}\right)^{-2}
    \bigg]
    \frac{1}{y^{\k+1}},
\qquad
    u\geq 0, \ y> 0,
\label{InegFt+integrable}
\end{equation}
which has a finite integral over $\R_+^2$ since $0<\k<1$.
We recall that
$
    \lim_{x\to+\infty} E\left(e^{-\gamma F^-(x)}\right)
=
    \frac{ (2\gamma)^{{\k/2}}}{\k\Gamma(\k)I_\k(2\sqrt{2\gamma})}
\in
    [0,1]
$
for $\gamma>0$ by \eqref{eqLimiteLaplaceF-},
and that $\Im_1\egloi F^-(h_t/2)$ by Proposition \ref{lemU1}.
Also, $x\mapsto E \left( e^{- \gamma F^-(x)}  \right)$ is nonincreasing for $\gamma> 0$.
So,
\begin{eqnarray}
    \nonumber
    F_{t,z_1,z_2,\e} (u,y)
& \geq &
    \un_{\{\frac{ 4\lambda u e^{-\theta\phi(t)}}{1+\e}\leq y\leq \frac{ 4\lambda u e^{-\theta\phi(t)}}{1-\e}\}}
    \frac{u^\k e^{-u}}{\Gamma(\k)}
    \bigg[
        \frac{ (z_1 y)^{\k}}{[\k\Gamma(\k)I_\k(2\sqrt{ z_1 y })]^2}
        -\left(1+\frac{z_2 y  }{\kappa+1}\right)^{-2}
    \bigg]
    \frac{1}{y^{\k+1}}
\\
& \geq &
    -\frac{u^\k e^{-u}}{\Gamma(\k)}
    \bigg|
        \frac{ (z_1 y)^{\k}}{[\k\Gamma(\k)I_\k(2\sqrt{ z_1 y })]^2}
        -\left(1+\frac{z_2 y  }{\kappa+1}\right)^{-2}
    \bigg|
    \frac{1}{y^{\k+1}}.
    \label{eqEHSIntegrable}
\end{eqnarray}
Notice that RHS of \eqref{eqEHSIntegrable}
has a finite integral over $\R_+^2$, since it is, for fixed $u$,
$O(y^{-\k})$ as $y\to 0$ and $O(y^{-\k-1})$ as $y\to+\infty$, where we used
$I_\k(x) =(x/2)^\k/\Gamma(\k+1)+(x/2)^{\k+2}/\Gamma(\k+2)+o(x^{\k+2})$ as $x\to 0$  by \eqref{eqDLBesselI} below,
and RHS of \eqref{eqEHSIntegrable} $\geq -2u^\k e^{-u}y^{-\k-1}/\Gamma(\k)$.
This and \eqref{InegFt+integrable} prove that
$|F_{t,z_1,z_2,\e}|$ is dominated by some function which does not depend on $t$ and has a finite integral over $\R_+^2$.
Moreover,
$
        \big|F_{t,(1-\e^2), (1-\e^2),\e}(u,y)\big|
        +
        \big|F_{t,(1-\e)^2,(1+\e)^2, \e}(u,y)\big|
\to_{t\to+\infty} 0
$
for every $(u,y)\in(\R_+^*)^2$ and $0<\e<1$.
So, by the dominated convergence theorem, RHS of \eqref{InegS5ValeurAbsolue}
$\to_{t\to+\infty} 0$.
So by  \eqref{InegS5ValeurAbsolue}, we get for every $0<\e<1$,
\begin{equation}\label{eqS5}
    S_{5,t,\e}^{\pm}
=
    o(e^{-\k \phi(t)}),
\qquad
    t\to+\infty.
\end{equation}

\noindent{\bf Step 6:}
We now prove that
$
    \Upsilon_{\epsilon}
\to_{\e\to 0}
    \Upsilon_0
$, and end the proof of Lemma \ref{LemmaPartB}.
We recall that
$
    \Upsilon_{\epsilon}
=
    \int_0^\infty\int_0^\infty f_\e(u,y)\dd y \dd u
$
as introduced in \eqref{eqDefUpsilonEpsilon}, where $f_\e$ is defined in \eqref{eqDeffepsilon}.
First, notice that  for every $(u,y)\in (0,\infty)^2$,
$$
    f_\e(u,y)
\to_{\e\to 0}
    \frac{u^\k e^{-u}}{\Gamma(\k)}
    \bigg[
        \frac{ y^{\k}}{[\Gamma(\k+1)I_\k(2\sqrt{ y})]^2}
        -\left(1+\frac{ y}{\kappa+1}\right)^{-2}
    \bigg]
    y^{-\k-1}
=:
    f_0(u,y).
$$

Using $|(1+c(1-x)^2)^{-2}-(1+c(1+x)^2)^{-2}|\leq 16 c x$ every for $x\in[0,1]$ and $c>0$, coming from Mean value theorem,
we have for $b_\e$ defined in \eqref{eqDefBEpsilon},
$$
    |b_\e(u,y)|
\leq
    u^\k e^{-u} y^{-\k-1}\Gamma(\k)^{-1}
    [
    16 y \un_{\{y< \k+1\}}/(\k+1)
    +
    2\un_{\{y\geq \k+1\}}
    ],
\qquad
u> 0, \ y>0,
$$
which has a finite integral over $(\R_+^*)^2$.
Moreover,
as in \eqref{eqEHSIntegrable},
$
    \Big| \frac{ \gamma^{{\k}}}{[\Gamma(\k+1)I_\k(2\sqrt{\gamma})]^2}
   -\left(1+\frac{\gamma}{\kappa+1}\right)^{-2}
    \Big|
$
is bounded by $2$ on $\R_+$, and
is $o(\gamma)$ as $\gamma\to0$.
and so is less than $|\gamma|$ for every $\gamma$ in some nonempty interval $[0,c)$.
So,
$
    |a_\e(u,y)|
\leq
    u^\k e^{-u}y^{-\k-1}\Gamma(\k)^{-1}
    [|y|\un_{\{y\leq c\}}+2 \un_{\{y>c\}}]
$
which also has a finite integral over $(\R_+^*)^2$.
Since
$
|f_\e(u,y)|
\leq
    |a_\e(u,y)|
    +|b_\e(u,y)|
$
for every $(u,y)\in(\R_+^*)^2$ (see the start of Step 4),
the dominated convergence theorem gives
\begin{equation}\label{eqConvergenceUpsilonEpsilon}
    \Upsilon_{\epsilon}
 \to_{\e\to0}
    \int_0^\infty
    \frac{u^\k e^{-u}}{\Gamma(\k)}
    \dd u
    \int_0^\infty
    \bigg[
        \frac{ y^{\k}}{[\Gamma(\k+1)I_\k(2\sqrt{ y})]^2}
        -\left(1+\frac{ y}{\kappa+1}\right)^{-2}
    \bigg]
    y^{-\k-1}
    \dd y
=
        \Upsilon_0,
\end{equation}
as defined in \eqref{eqDefUpsilon0}.

We proved that
$
    S_{5,t,\e}^{\pm}
=
    o(e^{-\k \phi(t)})
$
(in \eqref{eqS5} in Step 5)
and that
$
    \limsup_{t\to+\infty}[S_{4,t,\e}^{\pm}e^{\k \phi(t)}\lambda^{-\k} ]
\leq
    8^{\k}\Upsilon_{\epsilon}
$
(in \eqref{eqStep4} in Step 4).
So,
$
    \limsup_{t\to+\infty}[(S_{4,t,\e}^{\pm}+S_{5,t,\e}^{\pm})e^{\k \phi(t)}\lambda^{-\k} ]
\leq
    8^{\k}\Upsilon_{\epsilon}.
$
Then by \eqref{eqK1K4}, \eqref{eqConvergenceUpsilonEpsilon} and since $\k<1$,
\begin{equation}\label{InegLimsupS1S3}
    \limsup_{t\to+\infty} \left[(S_{1,t}^{\pm}-S_{3,t}^{\pm})e^{\k \phi(t)}\lambda^{-\k}\right]
\leq
    8^{\k}\Upsilon_{\epsilon}
\to_{\e\to 0}
    8^\k \Upsilon_0,
\end{equation}
so LHS of \eqref{InegLimsupS1S3} $\leq 8^\k \Upsilon_0$.
We prove with the similarly that
$
    \liminf_{t\to+\infty} \left[(S_{1,t}^{\pm}-S_{3,t}^{\pm})e^{\k \phi(t)}\lambda^{-\k}\right]
\geq
    8^{\k}\Upsilon_0
$.
This proves \eqref{eqLaplacePartB} and then Lemma \ref{LemmaPartB}.
\hfill$\Box$

\noindent {\bf Conclusion:} This and \eqref{eqConclusionPartA} gives
$
    S_{0,t}^{\pm}=S_{1,t}^{\pm}+S_{2,t}^{\pm}
=
    1+(\Upsilon_0-C_0 )8^{\k}\lambda^\k e^{-\k \phi(t)}+o(e^{-\k\phi(t)})
$.
This and \eqref{InegLaplaceUS0} give, since $\phi(t)=o(\log t)$,
$$
    1-\E\left(e^{-\frac{\lambda}{t} \bU}\right)
\sim_{t\to+\infty}
    8^\k(C_0- \Upsilon_0)\lambda^\k e^{-\k\phi(t)},
$$
which proves Proposition \ref{proplap}.
\hfill$\Box$





\mysection{Proof of the main results}\label{SectionPreuveLocEtAging}

\subsection{The renewal results : }

In this subsection we prove Propositions \ref{propRenewalJoint},  \ref{propTRT2} and Corollary \ref{propTRT}.
We start with some important intermediate result on the exit time $\bU$.
\noindent \\
In what follows, we mainly use the same ideas as in  Enriquez et al. (\cite{ESZ3} pages 441 to 443),
inspired by the book of Feller (\cite{Feller} pages 470 to 472). 
We start with the following lemma, for which we recall
that for all $i\geq 1$, $U_i=H(\tilde L_i)-H(\tilde m_i)$, and $C_\k$ is defined in Proposition \ref{proplap}.

\begin{lemma} \label{5.4} For $t>0$, let $\mu_t$ be the positive measure on $\R_+$ such that
$$
    \forall x\geq 0,
\qquad
    \mu_t([0,x])
:=
    e^{-\kappa \phi(t)}
    \sum_{j=1}^{n_t-2}\P\Bigg(\sum_{i=1}^j \frac{U_i}{t}\leq x\Bigg).
$$
Then, $(\mu_t)_t$ converges vaguely as $t\to+\infty$ to 
$\mu$ defined by
$\dd \mu(x):= (C_{\kappa} \Gamma(\kappa))^{-1}  x^{\kappa-1}\un_{(0,+\infty)}(x) \dd x$.
Moreover when $t\to+\infty$,
$x \mapsto e^{\kappa \phi(t)} \P(\bU/t \geq x)$
converges uniformly on every compact subset in $(0,+\infty)$ to
$x\mapsto C_{\kappa}x^{- \kappa} /\Gamma(1- \kappa)$.
\end{lemma}

\noindent{\bf Proof:}
First, let us prove that for all $\lambda>0$, we have as $t\to+\infty$,
\begin{eqnarray}
    \int_0^{+ \infty} e^{- \lambda x} \dd \mu_t(x)
& = &
    \int_0^{+ \infty} \frac{e^{-\lambda x}x^{\kappa-1}}{C_{\kappa}\Gamma(\kappa)}\textnormal{d}x +o(1),
    \label{cor2.1} \\
    \int_0^{+ \infty} e^{- \lambda x} e^{\kappa \phi(t)}\P(\bU/t \geq x) \textnormal{d}x
&  = &
    \int_0^{+ \infty} e^{- \lambda x} \frac{ C_{\kappa} } { \Gamma(1- \kappa) x^{\kappa} } \textnormal{d}x +o(1).
    \label{eq2}
\end{eqnarray}
Let $\lambda>0$. First, we have, by Proposition \ref{propind},
\begin{equation*}
\int_0^{+ \infty} e^{- \lambda x} \dd \mu_t(x)
=
     \frac{1}{e^{\kappa \phi(t)}}\sum_{j=1}^{n_t-2}\E\Big(e^{-\lambda\sum_{i=1}^{j} \frac{U_i}{t}}\Big)
=
    \sum_{j=1}^{n_t-2}
        \frac{1}{e^{\kappa \phi(t)}}
        \left[\E\left(e^{- \lambda \frac{\bU}{t}} \right) \right]^{ j}
    +O\left(\frac{(n_t)^2 e^{- \delta \kappa h_t}}{e^{\k \phi(t)}}\right).
\end{equation*}
We notice that, by Proposition  \ref{proplap},
$[\E\left(e^{- \lambda \bU/t}\right)]^{n_t-1}=o(1)$,
since $n_t e^{-\k\phi(t)}\to_{t\to+\infty}+\infty$ and $C_\k>0$.
Hence, 
we get as $t\to +\infty$,
again by Proposition  \ref{proplap}
and since $(n_t)^2 e^{- \delta \kappa h_t}e^{-\k \phi(t)}=o(1)$
because $\phi(t)=o(\log t)$,
\begin{equation*}
\int_0^{+ \infty} e^{- \lambda x} \dd \mu_t(x)
 =
    \frac{e^{-\kappa \phi(t)}(1
    +o(1))}{1-\E\left(e^{- \lambda \bU/t}\right)}
    +o(1)
 =
   \frac{1}{C_{\kappa} \lambda^{\kappa}}
   +o(1)
 =
    \int_0^{+ \infty} \frac{e^{-\lambda x}x^{\kappa-1}}{C_{\kappa}\Gamma(\kappa)} \textnormal{d}x+o(1),
\end{equation*}
which gives \eqref{cor2.1}. This implies the pointwise convergence of the Laplace transform of $(\mu_t)_t$ to that of $\mu$,
and then the vague convergence of $(\mu_t)_t$ to $\mu$
(e.g. by Feller \cite{Feller}, XIII.1 Th. 2c).
Now, we have as $t\to+\infty$ by Fubini and then Proposition \ref{proplap},
$$
    \lambda \int_{0}^{+ \infty}  e^{-\lambda x}\P\left(\bU/t \geq x\right)\textnormal{d}x
=
    \int_0^\infty\int_0^u\lambda e^{-\lambda x}\dd x\P(\bU/t\in \dd u)
=
    \E\left(1-e^{-\frac{\lambda}{t} \bU}\right)
=
    \frac{C_\k\lambda^\k+o(1)}{e^{\k\phi(t)}}.
$$
Since
$
    \lambda^{\kappa}
=
\lambda \int_{0}^{+ \infty} e^{-\lambda x}x^{- \kappa}\textnormal{d}x/ \Gamma(1- \kappa)
$,
we get \eqref{eq2}.  Again, this pointwise convergence of Laplace transforms gives the vague convergence of
$e^{\kappa \phi(t)}\P\left(\bU/t \geq x\right) \dd x$
to
$C_\k x^{-\k}\dd x /\Gamma(1-\k)$
as $t\to+\infty$.
Since $x\mapsto \P(\bU/t \geq x)$ is nonincreasing, we have for all $0<\e<x$,
\begin{equation*}
    \frac{1}{\e}\int_x^{x+\e}e^{\kappa \phi(t)}\P\left(\bU/t \geq u\right) \dd u
\leq
    e^{\kappa \phi(t)}\P\left(\bU/t \geq x\right)
\leq
    \frac{1}{\e}\int_{x-\e}^{x}e^{\kappa \phi(t)}\P\left(\bU/t \geq u\right) \dd u.
\end{equation*}
Using the previous vague convergence gives
\begin{equation*}
    \limsup_{t\to+\infty}e^{\kappa \phi(t)}\P\left(\bU/t \geq x\right)
\leq
    \frac{1}{\e}\lim_{t\to+\infty}\int_{x-\e}^{x}e^{\kappa \phi(t)}\P\left(\bU/t \geq u\right) \dd u
=
    \frac{1}{\e}\int_{x-\e}^{x} \frac{C_\k u^{-\k}}{\Gamma(1-\k)} \dd u.
\end{equation*}
Taking the limit as $\e\downarrow 0$ gives
$
    \limsup_{t\to+\infty}e^{\kappa \phi(t)}\P\left(\bU/t \geq x\right)
\leq
    C_\k x^{-\k}/\Gamma(1-\k)
$.
We prove similarly that
$
    \liminf_{t\to+\infty}e^{\kappa \phi(t)}\P\left(\bU/t \geq x\right)
\geq
    C_\k x^{-\k}/\Gamma(1-\k)
$.
This gives the pointwise convergence of
$x \mapsto e^{\kappa \phi(t)} \P(\bU/t \geq x)$
 to
$x\mapsto C_{\kappa}x^{- \kappa} /\Gamma(1- \kappa)$
on $(0,+\infty)$ as $t\to+\infty$.
Finally, since $x\mapsto  e^{\kappa \phi(t)} \P(\bU/t \geq x)$ is monotone and its pointwise limit
is continuous on $(0,\infty)$, Dini's theorem proves that this convergence is uniform
on every compact of $(0,\infty)$.
\hfill$\Box$

\medskip

We now introduce for $t>0$  the unique integer $\tilde N_t$ such that
$ H(\tilde m_{\tilde N_t}) \leq t < H(\tilde m_{\tilde N_t+1})$. We notice that
$\tilde N_t=N_t$ on $\mV_t$ due to Remark \ref{RemEgaliteAvecouSansTilde}.
We prove the following lemma:

\begin{lemma}\label{LemmaMajorationNt}
We have, as $t\to+\infty$,
$$
    \P\left( N_t \geq n_t\right)
=
    o(1),
\qquad\
    \P\big( N_t =0\big)
=
    o(1),
\qquad\
    \P\big( \tilde N_t \geq n_t\big)
=
    o(1),
\qquad\
    \P\big(\tilde N_t=0\big)
=
    o(1).
$$
\end{lemma}

\noindent{\bf Proof:}
First, by
equation \eqref{eqSommeHmk}
and the exponential Markov inequality,
$$
    \P\big( \tilde N_t \geq n_t\big)
\leq
    \P\bigg(\sum_{j=1}^{n_t-1} U_j  \leq  H(\tilde m_{n_t}) \leq t\bigg)
\leq
    e \E\Big(e^{- \sum_{i=1}^{n_t-1} U_i/t}\Big).
$$
This, together with Propositions \ref{propind} and \ref{proplap}, gives since $\phi(t)=o(\log t)$,
\begin{equation}
    \P\big( \tilde N_t \geq n_t\big)
\leq
    e \left(\E\big(e^{-\bU/t}\big) \right)^{n_t}
    +e C_2 n_t e^{- \delta \kappa h_t}
\leq
    C_+\exp\big(-c_-n_t e^{- \kappa \phi(t)}\big)
    + e C_2 n_t e^{- \delta \kappa h_t}
=
    o(1). \label{Ntnt}
\end{equation}
This proves the third inequality of the lemma.
Moreover,  $\P(\tilde N_t=0)=o(1)$ by Lemma \ref{lemtps}.
Finally, the first two inequalities follow from Lemma \ref{CVs}, since $\tilde N_t=N_t$ on $\mV_t$.
\hfill$\Box$

\medskip

\noindent \textbf{Proof of Proposition \ref{propRenewalJoint}:}
First, let $0<r <s<1$, and $a>0$.
Remark \ref{RemEgaliteAvecouSansTilde} and then
Lemmas \ref{CVs}, \ref{lemtps} and \ref{LemmaMajorationNt} give for $t>0$,
\begin{align}
&
    \P\big[
            1-s\leq H(m_{N_t})/t \leq 1-r,\
             H(m_{N_t+1})/t\geq 1+a\big]
    \label{TempsAtteinteDerniereVallee}
    \\
& \leq
    \sum_{j=1}^{n_t-1}
    \P\bigg(
            1-s\leq \frac{H(\tilde m_{\tilde N_t})}{t} \leq 1-r,
            \ \frac{H(\tilde m_{\tilde N_t+1})}{t}\geq 1+a,
            \ \tilde N_{t}=j, \ { \mathcal V}_t
      \bigg)
    \nonumber\\
&
\quad
    +\P\big(\tilde N_t\geq n_t\big)
    +\P\big(\tilde N_t=0\big)
    +P(\overline{ \mathcal V}_t)
    \nonumber\\
& \leq
    \sum_{j=1}^{n_t-1}
        \P\bigg( 1-s-\frac{2}{\log h_t} \leq \sum_{i=1}^{j-1} \frac{U_i}{t} \leq 1-r ,
                \ \sum_{i=1}^{j} \frac{U_i}{t} \geq 1+a-\frac{2}{\log h_t}
          \bigg)
    +o(1).
    \nonumber
\end{align}
We now use \eqref{DifferenceProbasDansLemmeLaplace} of Proposition \ref{propind}
and get for small $\e>0$, for large $t$,
\begin{equation}
    \eqref{TempsAtteinteDerniereVallee}
\leq
    \int_{1-s-\e}^{1-r} e^{\kappa \phi(t)}\P(\bU/t >1+a-\e-x)\dd \mu_t(x) +o(1).\label{PH}
\end{equation}
Using first the uniform convergence of
$u\mapsto e^{\kappa \phi(t)}\P(\bU/t >u)$ on the compact
$[a+r-\e,a+s]\subset (0,\infty)$
and then the vague convergence of $\mu_t$ (see Lemma \ref{5.4}),
we get
\begin{align*}
    \lim_{t \rightarrow +\infty}
    \int_{1-s-\e}^{1-r} e^{\kappa \phi(t)}\P(\bU/t >1+a-\e-x)\dd \mu_t(x)
=
    \int_{1-s-\e}^{1-r} \frac{x^{ \kappa-1} (1+a-\e-x)^{-\kappa}}{ \Gamma(\kappa) \Gamma(1-\kappa)}  \textnormal{d}x.
\end{align*}
Consequently, by letting $\e \to 0$, we obtain the first inequality of the following line:
$$
    \limsup_{t\to+\infty}     \eqref{TempsAtteinteDerniereVallee}
\leq
    \int_{1-s}^{1-r}
    \frac{ x^{\k-1}}{\Gamma(1-\k)\Gamma(\k)}\int_a^{\infty}\k(1+y-x)^{-\k-1}\dd y \dd x
\leq
    \liminf_{t\to+\infty}     \eqref{TempsAtteinteDerniereVallee}.
$$
We prove similarly the second inequality.
Consequently,
$
    \P[(H(m_{N_t})/t, H(m_{N_t+1})/t)\in \Delta]
$
converges to
$
    \int_\Delta
    \k x^{\k-1}(y-x)^{-\k-1}\un_{(0,1)}(x)\un_{(1,\infty)}(y)/(\Gamma(1-\k)\Gamma(\k))\dd x \dd y
$
as $t\to+\infty$ for $\Delta=[1-s,1-r]\times [1+a,+\infty)$
and so for every $\Delta=[a,b]\times [c,d]\subset (0,1)\times (1,\infty)$.
This gives the vague convergence on $(0,1)\times (1,+\infty)$,
and then the convergence in law of $(H(m_{N_t})/t, H(m_{N_t+1})/t)$
to
$\k x^{\k-1}(y-x)^{-\k-1}\un_{(0,1)}(x)\un_{(1,\infty)}(y)/(\Gamma(1-\k)\Gamma(\k))\dd x \dd y$,
which proves Proposition \ref{propRenewalJoint}
since $\Gamma(1-\k)\Gamma(\k)=\pi/\sin(\pi \k)$.
\hfill$\Box$

\medskip

\noindent{\bf Proof of  Corollary \ref{propTRT}:}
First by Proposition \ref{propRenewalJoint}, for $v\geq 0$ as $t\to +\infty$,
$$
    \lim_{t \rightarrow + \infty} \P\bigg(\frac{H(\mf_{N_t+1})}{t} \geq 1+ v\bigg)
=
    \frac{\sin(\pi \k)}{\pi}\int_v^\infty \dd u \int_0^1\k x^{\k-1}(1+u-x)^{-\k-1}\dd x.
$$
Using the change of variables $z=x/(1+u-x)$ in the second integral leads to  \eqref{Nt}.
\eqref{HNt} follows from Proposition \ref{propRenewalJoint} by straightforward computations.
Finally, for every continuous bounded function $\varphi$, Proposition \ref{propRenewalJoint}
and the change of variables $u=y-x$ give
$$
    \E\big[\varphi\big((H(m_{N_t+1})-H(m_{N_t}))/t\big)\big]
\to_{t\to +\infty}
    \int\varphi(u)\frac{\sin(\pi \k)}{\pi}u^{-\k-1}\int \k x^{\k-1}\un_{[0,1]}(x)\un_{\{x\geq 1-u\}}\dd x \dd u ,
$$
which gives the convergence in law of $(H(m_{N_t+1})-H(m_{N_t}))/t$ under $\P$ as $t\to+\infty$.
\hfill$\Box$

\medskip

\noindent{\bf Proof of Proposition \ref{propTRT2}:}
Let $\lambda>0$.
For $t>0$, we denote by $\nu_t$ the measure on $\R_+$ such that
$$
    \nu_t([0,x])
=
    e^{-\kappa \phi(t)}
    \sum_{j=1}^{n_t-1}
        \exp\bigg(- \frac{C_{\kappa} \lambda^{\kappa}  j} {e^{\kappa \phi(t)}}\bigg)
        \P\bigg(\sum_{i=1}^{j-1} \frac{U_i}{t}\leq x\bigg),
\qquad
    x\geq 0.
$$
In particular for $j=1$, $\P(\sum_{i=1}^{j-1} U_i/t\in .)=\P(0\in.)$ denotes the Dirac measure at $0$.
We first show that the Laplace transform of the measure $\nu_t$ converges when $t$ goes to infinity.
We consider $\alpha$ such that  $0<\lambda <\alpha$. We get as in the proof of Lemma \ref{5.4},
\begin{align*}
&
    \int_0^{+ \infty}
    e^{- \alpha u} \dd \nu_t(u)
=
    e^{-\k \phi(t)}
    \sum_{j=1}^{n_t-1}
        \exp\Big(- \frac{C_{\kappa} \lambda^{\kappa} j }{e^{\kappa \phi(t)}}\Big)
        \left[\E\left(e^{-\alpha\bU/t}\right)\right]^{j-1}
    +O(e^{-\k \phi(t)} n_t^2 e^{-\delta \k h_t})
\\
= &
    \frac
        {
            1- \Big[
                    \exp\Big(- \frac{C_{\kappa} \lambda^{\kappa}  } {e^{\kappa \phi(t)}}\Big)
                    \E\left(e^{-\alpha\bU/t}\right)
            \Big]^{n_t-1}
        }
        {
            \exp\Big(\kappa \phi(t)+\frac{C_{\kappa} \lambda^{\kappa}  } {e^{\kappa \phi(t)}}\Big)
            \Big[1-\exp\big(- \frac{C_{\kappa} \lambda^{\kappa}  } {e^{\kappa \phi(t)}}\big) \E\left(e^{-\alpha \bU/t}\right)\Big]
        }
    +o(1)
     =
      \frac{1}{C_{\kappa}(\alpha^{\kappa}+ \lambda^{\kappa})}
      +o(1),
\end{align*}
by Propositions  \ref{propind} and then \ref{proplap}.
We also notice that
\begin{align*}
    \frac{1}
         {\alpha^{\kappa}+\lambda^{\kappa}}
=
    \sum_{j=0}^{+ \infty}  \frac{(- \lambda^{\kappa})^j}{\alpha^{{\kappa(1+ j)}}}
=
    \sum_{j=0}^{+ \infty}  \frac{(- \lambda^{\kappa})^j}{\Gamma[\kappa(1+j)]} \int_{0}^{+ \infty}e^{- \alpha u} u^{\kappa(1+ j)-1}\textnormal{d}u.
\end{align*}
Moreover,
$
\int_0^\infty\sum_{j=0}^\infty
    \big|
        \frac{(-\lambda^\k)^{j}}{\Gamma[\kappa(1+j)]}
        \frac{e^{-\alpha u}u^{\kappa(1+ j)-1}}{C_{\kappa}}
    \big|
    \textnormal{d}u<\infty
$, since
$
    \sum_{j=0}^{+ \infty}  \frac{(\lambda u)^{\k j}}{\Gamma[\kappa(1+j)]}
=
    O(e^{(\lambda+\e) u})
$
as $u\to +\infty$ for any $\e>0$,
and $\k>0$.
So Fubini gives,
\begin{equation*}
    \forall \alpha>\lambda,
\qquad
    \lim_{t\to+\infty}
    \int_0^{+ \infty} e^{- \alpha u} \dd\nu_t(u)
=
    \int_{0}^{+ \infty} e^{- \alpha u} \dd\nu(u)+o(1),
\end{equation*}
where $\nu$ is the measure defined by
$
    \dd\nu(u)
=
    \frac{1}{C_{\kappa}} \sum_{j=0}^{+ \infty}
    \frac{(- \lambda^{\kappa})^ju^{\kappa(1+ j)-1}}{\Gamma[\kappa(1+j)]}\un_{\R_+}(u) \textnormal{d}u
$.
This pointwise convergence of the Laplace transform
of $\nu_t$ on $(\lambda,+\infty)$ leads to
the vague convergence of $\nu_t$ to $\nu$ as $t\to+\infty$
(e.g. by Feller \cite{Feller}, XIII.1 Th. 2c).

We have, with the arguments already used between \eqref{TempsAtteinteDerniereVallee} and \eqref{PH},
for any $a>0$ and $\e>0$,
\begin{eqnarray}
\label{EqLaplaceNbValleesVisiteesN}
&&
    \E\left[\exp\bigg(- \frac{C_{\kappa} \lambda^{\kappa}  N_t}{e^{\kappa \phi(t)}}\bigg),\
            \frac{H(m_{N_t+1})}{t}\geq 1+a
       \right]
\\
& \leq &
    \sum_{j=1}^{n_t-1}
            \exp\bigg(- \frac{C_{\kappa} \lambda^{\kappa}  j}{e^{\kappa \phi(t)}}\bigg)
    \P\left[
            \frac{H(\tilde m_{j})}{t}\leq 1,
            \frac{H(\tilde m_{j+1})}{t}\geq 1+a,
            \tilde N_t=j, \mV_t
    \right]+o(1)
\nonumber\\
& \leq &
    \sum_{j=1}^{n_t-1}\exp\bigg(- \frac{C_{\kappa} \lambda^{\kappa}  j} {e^{\kappa \phi(t)}}\bigg)
    \P\left[
             \sum_{i=1}^{j-1} \frac{U_i}{t}\leq 1,\ 1+a- \frac{2}{\log h_t} \leq \sum_{i=1}^{j} \frac{U_i}{t}
      \right]
     +o(1)
\nonumber\\
& \leq &
    \int_0^{1} e^{\kappa \phi(t)} \P\left(\frac{\bU}{t} > 1+a-\e-x\right)\dd \nu_t(x) +o(1),
\nonumber
\end{eqnarray}
since the term for $j=1$ in the third line is less than
$
\P[U_1/t\geq 1+a-\e]
\leq
    \P(\bU/t\geq 1+a-\e)
    +o(1)
=
    o(1)
$
by the case $n=1$ just after \eqref{DifferenceProbasDansLemmeLaplace}.

Using
the uniform convergence of
$x\mapsto e^{\kappa \phi(t)}\P\left({\bU}/{t} > 1+a-x\right)$ on $[0,1]$
provided by Lemma \ref{5.4} followed by the vague convergence of $\nu_t$ to $\nu$,  we get the first inequality of
\begin{align}\label{eqLimLaplaceNbValleesAveca}
    \limsup_{t\to+\infty} \eqref{EqLaplaceNbValleesVisiteesN}
 \leq
    \int_0^1\frac{C_\k}{\Gamma(1-\k)}(1+a-x)^{-\k}\dd \nu(x)
\leq
    \liminf_{t\to+\infty} \eqref{EqLaplaceNbValleesVisiteesN}.
\end{align}
We obtain the second one similarly. Since
$
    \lim_{a\downarrow 0}\lim_{t \to + \infty} \P(H(\mf_{N_t+1})/t< 1+a)=0
$
by equation \eqref{Nt} of Corrolary \ref{propTRT}, letting $a\downarrow 0$ in \eqref{eqLimLaplaceNbValleesAveca} gives
$$
    \lim_{t\to+\infty}
    \E\left[\exp\bigg(- \frac{C_{\kappa} \lambda^{\kappa}  N_t}{e^{\kappa \phi(t)}}\bigg)
       \right]
=
    \int_0^1\frac{C_\k(1-x)^{-\k}}{\Gamma(1-\k)}\dd \nu(x)
=
    \sum_{j=0}^{+ \infty}
    (- \lambda^{\kappa})^j
    \int_0^{1}
    \frac{(1-x)^{-\kappa}x^{\kappa(1+ j)-1}}
    {\Gamma[\kappa(1+j)]\Gamma(1- \kappa)}
    \textnormal{d}x.
$$
Since $\int_0^1 x^{a-1}(1-x)^{b-1}\dd x=\Gamma(a)\Gamma(b)/\Gamma(a+b)$ for every $a>0$ and  $b>0$,
changing $C_\k \lambda^\k$ into $u$ gives the pointwise convergence of
$\E[\exp(-u N_t/e^{\k\phi(t)})]$ to
$
    \sum_{j=0}^{+ \infty}
    \frac{1}{\Gamma(\kappa j+1) }
    \left(\frac{ -u}{C_{\kappa}}\right)^j
$,
which ends the proof of Proposition \ref{propTRT2}. \hfill$\Box$

\subsection{The localization : proof of Theorem \ref{ththm}}\label{SubSectProofLocalization}

We recall
the notation
$
    H_{x\rightarrow y }
=
    \inf\{s>H(x), \ X(s)=y\}-H(x)
$
for $(x,y)\in\R_+^2$,
which is equal to
$
    H(y)-H(x)
$
if $x<y$.
Let $ \phi^*(t):= \phi(t) / \g$, where $0<\g<1$ will be chosen later.
We define $t^*:= t- e^{ (1+2 \delta)\phi^*(t)}$,
\begin{equation*}
\begin{array}{llll}
&
\mathcal A_0:=\big\{1 \leq \tN_t < n_t \big\},
&&
\mathcal A_1:=\cap_{j=1}^{n_t-1} \big\{H_{\tL \rightarrow {\mt_{j+1} }}\leq 2t/ \log h_t   \big\},\\
&
\mathcal A_2:=\cap_{j=1}^{n_t-1} \big\{H_{\mt_j\rightarrow \mt_{j+1} }< H_{\mt_j\rightarrow \tilde L_{j}^- }  \big\},
&&
\mathcal A_3:= \big\{H\big(\mt_{\tN_t}\big)\leq t^{*} \big\}.
\end{array}
\end{equation*}
We also introduce
$\mI_j:=[\mt_{j}-  \phi^*(t)/ \g, \mt_{j}+ \phi^*(t)/ \g ]$, $j \in \N^*$.
Let $\e>0$. We have:
\begin{eqnarray}
    \P\left( X(t) \notin \mI_{\tN_t}\right)
& \leq &
    \P\left( X(t) \notin \mI_{\tN_t}, \tN_t=\tN_{t(1+\e)},\mathcal A_0,\mathcal A_1, \mathcal A_2, \mathcal A_3 \right)
    +\P\big(\tN_t\neq \tN_{t(1+\e)} \big)
    +\P\left(\overline{{\mathcal A}_0}\right)
\nonumber\\
&
    &+ \P\left( X(t) \notin \mI_{\tN_t}, \tN_t=\tN_{t(1+\e)},\mathcal A_0,\mathcal A_1, \mathcal A_2, \overline{{\mathcal A}_3} \right)
    +\P\left(\overline{{\mathcal A}_1}\right)
    +\P\left(\overline{{\mathcal A}_2}\right).
\label{eqLocalizationDecoupageEn3Parties}
\end{eqnarray}
We split the proof into three parts, in which we estimate these different probabilities. We start with:

\noindent
\textbf{Part 1:}  We prove that for large $t$,
\begin{equation}
\label{eqPart1}
    \P\left( X(t) \notin \mI_{\tN_t},\tN_t=\tN_{t(1+\e)},\mathcal A_0,\mathcal A_1, \mathcal A_2, \mathcal A_3 \right)
\leq
    n_t
    \left(
        C_+  h_t e^{- \k(1-\delta) \phi^*(t)/16}
        +
        e^{- (c_-) \delta  \phi^*(t)}
    \right).
\end{equation}
Let $\mathcal B_j:= \{\tN_t=\tN_{t(1+\e)}=j\} \cap\mathcal A_1 \cap \mathcal A_2$ for $j\in\N^*$.
We have
\begin{eqnarray}
    \P\left( X(t) \notin  \mI_{\tN_t},\ \tN_t=\tN_{t(1+\e)},\mathcal A_0,\mathcal A_1, \mathcal A_2, \mathcal A_3 \right)
=
    \sum_{j=1}^{n_t-1} \P\left( X(t) \notin \mI_{j}, \mathcal B_j \cap \mathcal A_3  \right).
\label{sum}
\end{eqnarray}
Let $1\leq j <n_t$.
Notice that on $\mathcal B_j$, $\tilde N_{t(1+\e)}=j$, so
$H(\tilde L_j)+H_{\tL \rightarrow {\mt_{j+1} }}=H(\tilde m_{j+1})> t(1+\e)$;
moreover for large $t$, $H_{\tL \rightarrow {\mt_{j+1} }}\leq \e t/2$,
then $H(\tilde L_j)>(1+\e/2)t$, and so for all $u\in[H(\tilde m_j), (1+\e/2)t]$,
$
    u
<
    H(\tilde L_j)
<
    H(\tilde m_{j+1})
=
    H(\tilde m_j)+H_{\tilde m_j\to \tilde m_{j+1}}
<
    H(\tilde m_j)+H_{\tilde m_j\to \tilde L_j^-}
$,
and then $X(u)\in(\tilde L_j^-,\tilde L_j)$.
That is, if $t$ is large enough, on $\mathcal B_j$, after first hitting $\tilde m_j$,
$X$ stays in $(\Lm, \tilde{L}_j)$ at least until time $(1+ \e/2)t$.
Therefore, conditioning on $H(\mt_j)$ and using the strong Markov property,
\begin{align}
    \P&\left( X(t) \notin \mI_{j}, \mathcal B_j \cap \mathcal A_3  \right)
\nonumber\\
&  \leq
    E\Big(
        \P^{\wk}\big[X(t) \notin \mI_{j},\ H(\mt_j)\leq t^{*}, \ \forall u\in[H(\mt_j), (1+\e/2)t],      X(u)\in\big(\Lm, \tL\big)\big]
     \Big)
\nonumber\\
& =
    E\bigg(\int_0^{t^{*}} \P^{\wk}( H(\mt_j) \in ds)
        \P^{\wk}_{\mt_j}\left[X(t-s) \notin \mI_{j},\ \forall u\in\big[0,(1+\e/2)t-s\big], X(u)\in\big(\Lm,\tL\big)  \right]
     \bigg)
\nonumber\\
& \leq
    E\bigg(
        \sup_{  0 \leq s \leq t^{*} }
        \P^{\wk}_{\mt_j}\left[ X(t-s) \notin \mI_{j},\forall u\in[0,(1+\e/2)t-s],X(u)\in\big(\Lm,\tL\big)   \right]
     \bigg).
\label{eqIntermediairePartie6}
 \end{align}
Now, as in Brox (\cite{Brox}, proof of Prop. 4.1) we introduce
a coupling between $X$ (under $\P^{\wk}_{\mt_j}$) and a reflected diffusion $Y_j$ defined below.
To this aim, let
$\big(Y_j^{(x)}(u),\ u\geq 0\big)$
be (in words) a diffusion in the potential $\wk$, starting from $x\in\big[\Lm,\tL\big]$ and reflected at $\Lm$ and $\tL$.
We denote its law by $\widehat{P}_{j,x}^{W_{\kappa}}$.
More precisely, this process $Y_j^{(x)}$ is defined as in Brox (\cite{Brox}, p. 1216) by
$$
    Y_j^{(x)}(u)
:=
    A^{-1}\Big(\widehat{B}_j^{(x)}\Big(\widehat{T}_{j,x}^{-1}(u)\Big)\Big),
\qquad
    u\geq 0,
    \ x \in [\Lm, \tL],
$$
where $(\widehat{B}^{(x)}_j(s),\ s\geq 0)$ is a one-dimensional Brownian motion independent from $\wk$,
starting from $A(x)$ and reflected at $ A\big(\Lm\big)$ and $ A\big( \tL\big)$,
and $\widehat{T}_{j,x}$ is defined  as $T$ (see \eqref{T}) replacing $B$ by $\hat B_j^{(x)}$.
This enables us to define $(Y_j(s),\ s\geq 0)$
by
$\widehat{P}_j^{W_{\kappa}}(Y_j\in .):=\int_{\Lm}^{ \tL} \widehat{P}_{j,x}^{W_{\kappa}}(.)\dd \tilde \mu_j(x)$, where
\begin{align}\label{def_mu}
    \dd \tilde \mu_j(x)
:=
    \exp(-\tilde V^{(j)}(x)) \un_{[\Lm,\tL]}(x)\dd x
    \bigg(\int_{\Lm}^{ \tL} \exp(-\tilde V^{(j)}(y))dy\bigg)^{-1}.
\end{align}
As is proved in (\cite{Brox}, proof of Prop. 4.1), $\tilde \mu_j$ is invariant for the semi-group of $Y_j$;
in particular $\widehat{P}^{W_{\kappa}}_j(Y_j(s) \in \Delta)=\tilde \mu_j(\Delta)$ for every $s\geq 0$ and $\Delta \subset [\Lm, \tL]$.

We can now, as in \cite{Brox}, build a coupling ${Q}^{W_{\kappa}}_{\mt_j}$ of $X$ and $Y_j$, such that
${Q}^{W_{\kappa}}_{\mt_j}(Y_j\in.)=\widehat{P}^{W_{\kappa}}_j(Y_j\in .)$,
and ${Q}^{W_{\kappa}}_{\mt_j}(X\in.)=\P^{\wk}_{\mt_j}(X\in.)$,
these two Markov processes $Y_j$ and $X$ move independently until
$$
    \widehat H_{\{X=Y_j\}}
:=
    \inf\{u\geq 0,\ X(u)=Y_j(u)\},
$$
which is the first collision,
then $X(u)=Y_j(u)$ until the next exit time
$$
    \widehat H_j^{{\textnormal{exit}}}
:=
    \inf\big\{u> \widehat H_{\{X=Y_j\}},\ X(u) \notin (\Lm, \tL) \big\},
$$
and then $X$ and $Y_j$  move independently again.

We introduce (see Figure \ref{fig2})
$$
    t_1:=t-t^{*}=e^{ (1+ 2 \delta)\phi^*(t)},
\qquad
    \widehat L_j^+
:=
    \tilde\tau_j(\phi^*(t)),
\qquad
    \widehat L_j^-
:=
    \tilde\tau_j^-(\phi^*(t)).
$$
We now prove that, with a large probability, under ${Q}^{W_{\kappa}}_{\mt_j}$,
$X$ and $Y_j$ first collide before time $t_1$,
that is, $\widehat H_{\{X=Y_j\}}\leq t_1$.
To this aim,
we first provide a result concerning only the hitting times $H$ of $X$:

\begin{figure}[h]
\begin{center}
\scalebox{1.5}{\input{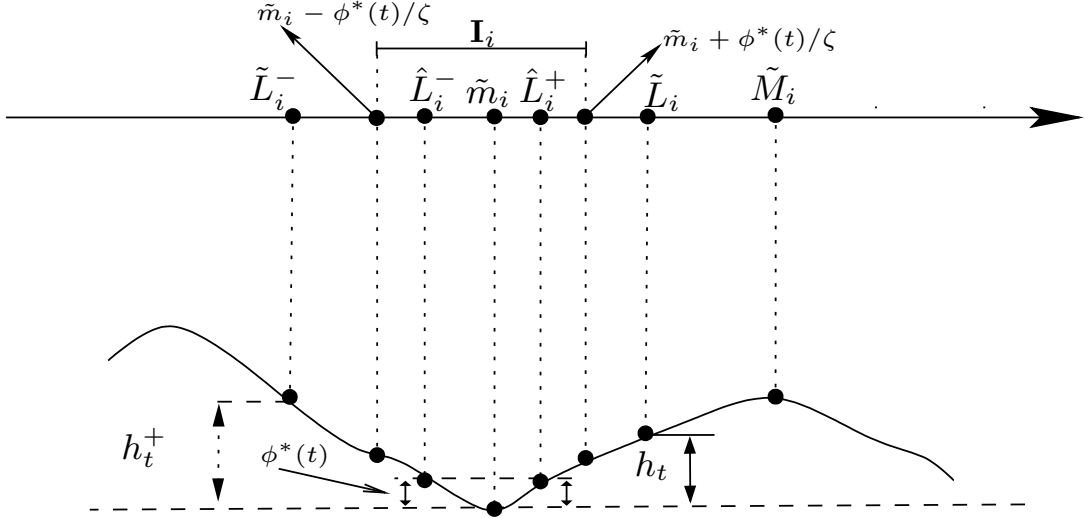}}
\caption{Schema of valley number $i$ and some notation around its minimum for small $\zeta$}\label{fig2}
\end{center}
\end{figure}

\begin{lemma}\label{LemmaProbaCouplingSortieXValleeJ}
For large $t$,
with a probability larger than $1-e^{- (c_-) \delta  \phi^*(t)}$,
\begin{equation}
    Q^{\wk}_{\mt_j}\left(H\big( \widehat L_j^+\big)  \vee  H\big(\widehat L_j^-\big) > e^{(1+\delta) \phi^*(t)} \right)
\leq
    C_+  \phi^*(t) e^{-\delta  \phi^*(t)/8},
\label{6.47}
\end{equation}
uniformly for $j\in[1, n_t)$ and with $x\vee y:=\max(x,y)$,
\end{lemma}
\noindent
The proof of this lemma is deferred to Subsection \ref{SubSectProofTechnicalLemmasProbaCouplingQW}.
We deduce from \eqref{6.47} and the definition of $x\vee y=\max(x,y)$
that with probability larger than $1-e^{- (c_-) \delta  \phi^*(t)}$,
\begin{align*}
&
{Q}^{\wk}_{\mt_j} \big( \widehat H_{\{X=Y_j\}} > t_1 \big)
\\
& \leq
   {Q}^{\wk}_{\mt_j} \big[\widehat H_{\{X=Y_j\}} > H( \widehat L_j^+)  \vee  H(\widehat L_j^-)   \big]
    +
    {Q}^{\wk}_{\mt_j} \big[H(\widehat L_j^+)\vee H(\widehat L_j^-)>t_1   \big]
    \\
& \leq
    {Q}^{\wk}_{\mt_j} \left( \widehat H_{\{X=Y_j\}} > H(\widehat L^-_j), Y_j(0) < \mt_j \right)
    +{Q}^{\wk}_{\mt_j} \left( \widehat H_{\{X=Y_j\}} > H(\widehat L_j^+), Y_j(0) \geq \mt_j \right)
    +  \frac{C_+\phi^*(t)}{e^{\delta   \phi^*(t)/8}}.
\end{align*}
On $\big\{\widehat H_{\{X=Y_j\}} > H(\widehat L^-_j), Y_j(0)<\mt_j\big\}$,
$Y_j(0)-X(0)=Y_j(0)-\tilde m_j<0$ under ${Q}^{\wk}_{\mt_j}$,
and by continuity $Y_j-X<0$ up to time $\widehat H_{\{X=Y_j\}}$
and in particular at time $H(\widehat L_j^-)$, so
\begin{eqnarray}
    {Q}^{\wk}_{\mt_j} \big( \widehat H_{\{X=Y_j\}} > t_1 \big)
& \leq &
    {Q}^{\wk}_{\mt_j} \left(\  Y_j[H(\widehat L_j^-)] \in [\Lm,\widehat L_j^-], \ \widehat H_{\{X=Y_j\}}> H(\widehat{L}^-_j) \right)
    \nonumber\\
&    & +
    {Q}^{\wk}_{\mt_j} \left(\  Y_j[H(\widehat L_j^+)] \in [\widehat L_j^+, \tL], \  \widehat H_{\{X=Y_j\}} > H(\widehat L_j^+)\right)
    + \frac{C_+  \phi^*(t) }{e^{\delta  \phi^*(t)/8}}
\nonumber\\
& \leq &
    \tilde \mu_j\big(\big[\Lm,\widehat L_j^-\big]\big)
    +\tilde \mu_j\big(\big[\widehat L_j^+,\tL\big]\big)+ C_+  \phi^*(t) e^{- \delta  \phi^*(t)/8},
    \label{P5B}
\end{eqnarray}
where the last line comes from the independence of $X$ and $Y_j$ until $\widehat H_{\{X=Y_j\}}$
and since
$
    {Q}^{W_{\kappa}}_{\mt_j}(Y_j(s)\in \Delta)
=
    \widehat{P}^{W_{\kappa}}_j(Y_j(s)\in \Delta)
=
    \tilde \mu_j(\Delta)
$,
for every $s\geq 0$ and $\Delta \subset [\Lm, \tL]$ as explained after \eqref{def_mu}.

We would like to bound \eqref{eqIntermediairePartie6}.
Let $s\in[0,t^{*}]$.
Using first $t_1 =t-t^*\leq t-s \leq (1+\e/2)t-s$,
and second
$X(u)=Y_j(u)$ for every $\widehat H_{\{X=Y_j\}} \leq u\leq \widehat H_j^{{\textnormal{exit}}}$
and
$
    {Q}^{W_{\kappa}}_{\mt_j}(Y_j(t-s)\in .)
=
    \tilde \mu_j(.)
$,
\begin{eqnarray}\label{P6BLine1}
&&
    Q^{\wk}_{\mt_j}\left( X(t-s) \notin \mI_{j},\  \forall u\in[0,(1+\e/2)t-s],\
    X(u)\in(\Lm,\tL)   \right)
\\
& \leq &
    Q^{\wk}_{\mt_j}\big(\widehat H_{\{X=Y_j\}}>t_1\big)
    +
    Q^{\wk}_{\mt_j}
        \left(
            \widehat H_{\{X=Y_j\}} \leq {t_1}\leq t-s\leq (1+\e/2)t-s \leq \widehat H_j^{{\textnormal{exit}}},
            X(t-s) \notin \mI_{j}
        \right)
\nonumber\\
& \leq &
    Q^{\wk}_{\mt_j}\big(\widehat H_{\{X=Y_j\}}>t_1\big)
    +
   Q^{\wk}_{\mt_j}\big(Y_j(t-s)\notin \mI_{j}\big)
\nonumber\\
& \leq &
    Q^{\wk}_{\mt_j}\big(\widehat H_{\{X=Y_j\}}>t_1\big)
    +
   \tilde  \mu_j\big(\big[\Lm, \tL\big] \smallsetminus \mI_{j}\big).
\label{P6B}
\end{eqnarray}
In the following lemma, we show that with high probability,
the invariant probability measure $\tilde \mu_j$ is highly concentrated on
the small neighborhood $\mI_{j}$ of $\tilde m_j$.
The proof of this lemma is deferred to Subsection \ref{SubSectProofTechnicalLemmasProbaCouplingQW}.
More precisely:

\begin{lemma} \label{lemma6.3}
If $\zeta<\k/8$,
for all $1\leq j < n_t$,
with a probability greater than $1-  e^{- (c_-)\delta  \phi^*(t)  } $,
\begin{equation}
    \tilde \mu_j\big(\big[\Lm, \tL\big] \smallsetminus \mI_{j}\big)
\leq
    \tilde  \mu_j\big(\big[\Lm,\widehat L_j^-\big]\big)+\tilde \mu_j\big(\big[\widehat L_j^+, \tL\big]\big)
\leq
    C_+ h_t e^{- (1- \delta) \kappa \phi^*(t) /16}.
\label{IM1}
\end{equation}
\end{lemma}
\noindent
Notice that for any $\g>0$,
the right hand side of \eqref{IM1}
go to $0$ as $t\to+\infty$
since $h_t\leq \log t$ and $\log \log t=o( \phi(t))=o( \phi^*(t))$, and is
$o(1/n_t)$ if $\zeta$ is chosen small enough.

Moreover, notice that we can replace $Q^{\wk}_{\mt_j}$ by $\P^{\wk}_{\mt_j}$ in the first line \eqref{P6BLine1}
since ${Q}^{W_{\kappa}}_{\mt_j}(X\in.)=\P^{\wk}_{\mt_j}(X\in.)$.
So, \eqref{P6B}, \eqref{P5B} and then Lemma \ref{lemma6.3}
give
\begin{align*}
&
    \sup_{0\leq s\leq t^*}
    \P^{\wk}_{\mt_j}
    \left( X(t-s) \notin \mI_{j},\  \forall u\in[0,(1+\e/2)t-s],\
    X(u)\in(\Lm,\tL)   \right)
\\
& \leq
        \tilde \mu_j\big(\big[\Lm,\widehat L_j^-\big]\big)
    +\tilde \mu_j\big(\big[\widehat L_j^+,\tL\big]\big)+  \frac{C_+ \phi^*(t)}{e^{ \delta  \phi^*(t)/8}}
    +
   \tilde  \mu_j\big(\big[\Lm, \tL\big] \smallsetminus \mI_{j}\big)
\leq
    \frac{C_+ h_t }{e^{ (1- \delta) \kappa \phi^*(t) /16}}
    + \frac{C_+  \phi^*(t)}{e^{ \delta  \phi^*(t)/8}}
\end{align*}
with probability at least $1-  e^{- (c_-)\delta  \phi^*(t)  } $.
Finally, integrating this and applying successively \eqref{sum} and \eqref{eqIntermediairePartie6}
leads to \eqref{eqPart1},  which ends the proof of this Part 1.


\bigskip\noindent
\textbf{Part 2 }: We prove that there exists $c_3>0$ such that
if $\zeta\leq \k/48$,
\begin{equation}\label{eqPart2}
    \P\left( X(t) \notin \mI_{\tN_t},\ \tN_t=\tN_{t(1+\e)},\mathcal A_0,\mathcal A_1, \mathcal A_2, \overline{{\mathcal A}_3} \right)
\leq
    n_t e^{- c_3 \phi^*(t)/\zeta}.
\end{equation}
First, we prove similarly as in Part 1 from \eqref{sum} to \eqref{eqIntermediairePartie6} that,
using $H(\tilde m_j)\leq t$ on $\mathcal B_j$,
\begin{align}
&
    \P\left( X(t) \notin \mI_{\tN_t},\ \tN_t=\tN_{t(1+\e)},\mathcal A_0,\mathcal A_1, \mathcal A_2, \overline{{\mathcal{A}}_3 }\right)
=
    \sum_{j=1}^{n_t-1} \P\left( X(t) \notin \mI_{j}, \mathcal B_j \cap \overline{\mathcal A_3}  \right)
\label{eqPart2Debut}\\[-2mm]
& \leq
    \sum_{j=1}^{n_t-1}
        E\left( \sup_{  t^{*} \leq s \leq t }
                \P^{\wk}_{\mt_j}
                \left(
                    X(t-s) \notin  \mI_{j},\
                    \forall u\in[0,(1+\e/2)t-s],\ X(u)\in\big(\Lm,\tL\big)
                \right)
        \right)
\nonumber\\
& \leq
\sum_{j=1}^{n_t-1}
    E\left( \sup_{  t^{*} \leq s \leq t }
            \P^{\wk}_{\mt_j}
            \big[
                H(\mt_{j}- \phi^*(t)/ \g) \wedge H(\mt_{j}+ \phi^*(t)/ \g)<t-s
            \big]
    \right),
\nonumber
\end{align}
by definition of $\mI_{j}=[\mt_{j}-  \phi^*(t)/ \g, \mt_{j}+ \phi^*(t)/ \g ]$,
and with $x\wedge y:=\inf(x,y)$.
Since
$t-t^*=e^{ (1+2 \delta)\phi^*(t)}$,
this gives
\begin{equation}\label{InegProbaPart2Suite}
    \eqref{eqPart2Debut}
\leq
     n_t \sup_{1\leq j\leq n_t-1}
    E\left(
            \P^{\wk}_{\mt_j}
            \big[
                H(\mt_{j}- \phi^*(t)/ \g) \wedge H(\mt_{j}+ \phi^*(t)/ \g)<e^{ (1+2 \delta)\phi^*(t)}
            \big]
    \right).
\end{equation}
We estimate this with the following lemma, the proof of which is deferred to Subsection \ref{SubSectProofTechnicalLemmasProbaCouplingQW}:
\begin{lemma}\label{LemmaProbaHitttingTimesHjPlusMoins}
Assume $\zeta\leq \k/48$.
For large $t$, for each $1\leq j<n_t$,  with probability at least
$1-e^{-(c_-)\phi^*(t)/\zeta} $,
$$
    \P^{\wk}_{\mt_j}\left(
                            H(\mt_{j} \pm \phi^*(t)/ \g)
                            \leq
                            e^{ (1+2 \delta)\phi^*(t)}
                    \right)
\leq
    e^{-  (c_-) \phi^*(t)/\zeta}.
$$
\end{lemma}
\noindent
This together with \eqref{InegProbaPart2Suite} leads to \eqref{eqPart2}, which ends this Part 2.


\medskip\noindent
\textbf{Part 3:} We prove that if $0<\delta<1/16$
and $(1+2\delta)\k<1$, as $t\to+\infty$,
\begin{equation}\label{eqPart3Localization}
\P\big( \tN_t \neq \tN_{t(1+\e)}\big)
+\P\left(\overline{{\mathcal A}_0}\right)
+\P\left(\overline{{\mathcal A}_1}\right)
+\P\left(\overline{{\mathcal A}_2}\right )
\leq
\e^{1- \kappa}/(1-\k)+o(1).
\end{equation}
First,
$
\P\big(\overline{{\mathcal A}_0}\big)=o(1)
$
by Lemma \ref{LemmaMajorationNt}.
The fact that
$\P\big(\overline{{\mathcal A}_1}\big)=o(1)$ for such $\delta$
follows from \eqref{InegProbaSommeTempsEntreValleesPourLemma35}
since $(\log h_t) n_t e^{h_t}/t =o(1)$ as stated in Lemma \ref{lemtps}.
Furthermore, by Lemma \ref{CVs} and then (\ref{Nt}),
$$
    \P\big(\tN_t \neq \tN_{t(1+\e)} \big)
\leq
    \P[H(m_{N_{t+1}})<t(1+\e)]
    +P\big(\overline{\mathcal V}_t\big)
\leq
    \e^{1- \kappa}/(1-\k)+o(1).
$$
Moreover by the strong Markov property, recalling that
$
    \tilde L_i^-
<
\tilde m_j
<
    \tilde L_i^*
<
    \tilde L_i
<
    \tilde m_{j+1}
$
by \eqref{eqDefLi*},
\begin{eqnarray*}
\P^{\wk}_{\mt_j} \left(H({\mt_{j+1} })< H({\tilde L_{j}^-}) \right)
&=&
    \P^{\wk}_{\mt_j} \left(H({\tL })< H({\tL^-}) \right)\times
    \P^{\wk}_{\tL} \left(H({\mt_{j+1} })< H({\tL^-}) \right)
\\
& \geq &
    \P^{\wk}_{\mt_j} \left(H({\tL })< H({\tL^-}) \right)
    \times \P^{\wk}_{\tL} \left(H({\mt_{j+1} })< H({\tL^*}) \right)
\\
& \geq &
    \big(1-e^{-\k h_t/2  }\big)
    \times
    \big(1-2n_t e^{-h_t/8}\big)
    =1-o(1/n_t),
\end{eqnarray*}
for all $1\leq j \leq n_t-1$
with probability $\geq 1-C_+n_t e^{-\k \delta h_t}$
by Lemmas \ref{LemmaProbaQuitterVallesaDroite} and \ref{LemmaProbaRetourEnmi}
since $\d<1/16$.
This proves that $\P\left(\overline{{\mathcal A}_2}\right )=o(1)$.
This leads to \eqref{eqPart3Localization}, which ends this Part 3.


\noindent{\bf Conclusion:}
We now choose $\zeta=\delta^2$.
We recall that $\phi(t)\leq h_t\leq \log t=\exp(o(\phi(t))$ since $\log \log t=o( \phi(t))$.
Combining \eqref{eqLocalizationDecoupageEn3Parties}, \eqref{eqPart1}, \eqref{eqPart2} with \eqref{eqPart3Localization},
and choosing $\delta$ small enough gives
$
    \limsup_{t\to+\infty}
    \P\big( X(t) \notin \mI_{\tN_t}\big)
\leq
    \e^{1- \kappa}/(1-\k)
$,
for every $\e>0$, and so is $0$.
So,
$$
    \lim_{t \rightarrow + \infty }
    \P\big( |X(t)- \mf_{N_{t}}| > \phi(t)/\zeta^2 \big)
=
    \lim_{t \rightarrow + \infty }
    \P\big( |X(t)- \tilde m_{\tilde N_{t}}| > \phi(t)/\zeta^2 \big)
=
    \lim_{t\to+\infty}
    \P\big( X(t) \notin \mI_{\tN_t}\big)
=
0,
$$
since $\tilde N_t=N_t$ and $\tilde m_j=m_j$ for  $1\leq j\leq n_t$
on ${\mathcal V}_t$ by Remark \ref{RemEgaliteAvecouSansTilde}, $P\big(\overline{\mathcal V}_t\big)=o(1)$
by Lemma \ref{CVs}
and $\P(\tilde N_t\notin [1, n_t])=o(1)$
by Lemma \ref{LemmaMajorationNt}.
This proves  Theorem \ref{ththm}.
$\hfill\Box$



\subsection{The aging : Proof of Proposition \ref{ThAging}}
We fix $\alpha>1$ and $\delta>0$. We recall that the r.v. $(m_i)_i$ depend on $h_t$.
In what follows, we apply Theorem \ref{ththm} first at time $t$ with function $\phi$,
and second at time $\alpha t$ with a function $\phi_\alpha$ defined by
$\log(\alpha t)-\phi_\alpha(\alpha t) = \log t-\phi(t)$,
so that the r.v. $m_i$ are the same in both cases.
By Theorem \ref{ththm}  and Lemma \ref{LemmaMajorationNt}, we get
\begin{align*}
 &
    \P(\left|X(\alpha t)-X(t)\right| \leq 3\mathcal{C}_1 \phi(t),\
        N_t<N_{\alpha t})
\\
& \leq
    \P\big[
        \left|m_{N_{\alpha t}}-m_{N_t}\right|
        \leq
        \mathcal{C}_1[ 4\phi(t)+\phi_\alpha(\alpha t)],
        N_t<N_{\alpha t}\leq n_{\alpha t},
        N_t\in[1, n_t)
      \big]
        +\P\big[N_t\notin[1,n_t)\big]
\\
&
        +\P\big[|X(t)-m_{N_t}|>\mathcal{C}_1\phi(t)\big]
        +\P\big[|X(\alpha t)-m_{N_{\alpha t}}|>\mathcal{C}_1\phi_\alpha(\alpha t)\big]
        +\P\big[N_{\alpha t}\notin[1,n_{\alpha t})\big]
\\
& \leq
    \sum_{i=1}^{n_t} \sum_{j=i+1}^{n_{\alpha t}}
        \P\Big(m_{j}-m_i \leq 4 \mathcal{C}_1 \phi(t)+\mathcal{C}_1\phi_\alpha(\alpha t) \Big)
    +o(1),
\end{align*}
as $t\to+\infty$. So, \eqref{2.0} and Lemma \ref{CVs} leads to
$
    \P(\left|X(\alpha t)-X(t)\right| \leq  \mathcal{C}_1 \phi(t), N_t<N_{\alpha t})
=
    o(1)
$
since $n_t=e^{o(h_t)}$ and $n_{\alpha t}=e^{o(h_t)}$.
Consequently,
\begin{eqnarray}\label{InegAging1}
    \P(\left|X(\alpha t)-X(t)\right|  \leq 3\mathcal{C}_1 \phi(t))
&    = &
    \P(
        \left|X(\alpha t)-X(t)\right| \leq 3\mathcal{C}_1 \phi(t),
        \  N_t=N_{\alpha t}
      )
    +o(1)
\nonumber\\
& \leq &
    \P(N_t=N_{\alpha t})+o(1).
\end{eqnarray}
Moreover, by Theorem \ref{ththm} applied at time $t$ with function $\phi$, we have for large $t$,
\begin{eqnarray*}
&&
    \P(\left|X(\alpha t)-X(t)\right| > 3\mathcal{C}_1 \phi(t),\ N_t=N_{\alpha t})\\
& \leq &
    \P( \left|X(\alpha t)-m_{N_t}\right| > 3\mathcal{C}_1 \phi(t) -\mathcal{C}_1 \phi(t),\ N_t=N_{\alpha t})
    +\P[|X(t)-m_{N_t}|>\mathcal{C}_1 \phi(t)]
\\
& \leq &
    \P[ \left|X(\alpha t)-m_{N_{\alpha t}}\right| > 2\mathcal{C}_1 \phi(t)]
    + o(1)
\\
& \leq &
    \P( \left|X(\alpha t)-m_{N_{\alpha t}}\right| > \mathcal{C}_1 \phi_\alpha(\alpha t))
    + o(1).
\end{eqnarray*}
This is $o(1)$ as $t\to+\infty$,
by Theorem \ref{ththm} applied at time $\alpha t$ with function $\phi_\alpha$.
Therefore,
\begin{eqnarray*}
    \P( N_t=N_{\alpha t})
& = &
    \P( \left|X(\alpha t)-X(t)\right| \leq 3\mathcal{C}_1 \phi(t), N_t=N_{\alpha t})
    +o(1)
\\
& \leq &
    \P(\left|X(\alpha t)-X(t)\right| \leq 3\mathcal{C}_1 \phi(t))
    +o(1).
\end{eqnarray*}
This together with \eqref{InegAging1} gives
$$
    \P(\left|X(\alpha t)-X(t)\right| \leq 3\mathcal{C}_1 \phi(t) )
=
    \P( N_t=N_{\alpha t})
    +o(1)
=
    \P[H(m_{N_t+1}) >\alpha t ]
    +o(1).
$$
This, combined with \eqref{Nt} and the change of variables
$u=1/(1+x)$ proves Proposition \ref{ThAging}, since $\phi$ is choosen up to a multiplicative constant.
\hfill $\Box$

\subsection{Proof of Lemmas \ref{LemmaProbaCouplingSortieXValleeJ}, \ref{lemma6.3} and  \ref{LemmaProbaHitttingTimesHjPlusMoins}}
\label{SubSectProofTechnicalLemmasProbaCouplingQW}~~

\noindent{{\bf Proof of Lemma \ref{LemmaProbaCouplingSortieXValleeJ}:}
Let $j \in[1, n_t)$.
First by \eqref{eqScaleFunctionX},
$
    \P^{\wk}_{\mt_j} \left(H(\widehat L_j^+) > H\big[\tilde\tau_j^-((1+\delta/2)\phi^*(t))\big]\right)
\leq
    \widehat{Q}_j/\widehat{D}_j
$,
where
$
    \widehat{Q}_j
:=
    \int_{\tilde m_j}^{\widehat L_j^+}
    e^{\tilde V^{(j)}(x)}\dd x
$
and
$
    \widehat{D}_j
:=
    \int_{\tilde\tau_j^-((1+\delta/2)\phi^*(t))}^{\widehat m_j}
    e^{\tilde V^{(j)}(x)}\dd x
$.
Recall that
$
    \widehat L_j^+
=
    \tilde\tau_j(\phi^*(t))
$.
With a method similar as for \eqref{eq3.13}, we get
with a probability larger than $1-e^{-(c_-) \delta  \phi^*(t)}$,
\begin{equation*}
    \widehat{Q}_j
\leq
    (\widehat L_j^+-\tilde m_j)\exp\big[\max\nolimits_{[\tilde m_j, \widehat L_j^+]}\tilde V^{(j)}\big]
\leq
    8\k^{-1}\phi^*(t)\exp(\phi^*(t)),
\qquad
    \widehat{D}_j
\geq
    e^{(1+\delta/8)\phi^*(t)}
\end{equation*}
by \eqref{bessel} for the first inequality and by
Fact \ref{Fact_Williams} {\bf (ii)} and
\eqref{MinorationAVallee} with $h=(1+\delta/2)\phi^*(t)$ and $\alpha=1-\frac{1+\delta/8}{1+\delta/2}$
for the second one, together with Lemma \ref{CVs} since $\phi(t)=o(\log t)$.
Hence, with such a probability,
\begin{equation}\label{InegProbaSortirPetitVoisinageGauche}
    \P^{\wk}_{\mt_j} \big[H(\widehat L_j^+) > H[\tilde\tau_j^-((1+\delta/2)\phi^*(t))]\big]
\leq
    C_+  \phi^*(t)e^{-\delta\phi^*(t)/8}.
\end{equation}
Moreover, as in \eqref{4.38} and by scaling, we have under $\P^{\wk}_{\mt_j}$ on
$\big\{ H(\widehat L_j^+) <  H[\tilde\tau_j^-((1+\delta/2)\phi^*(t))]  \big\}$,
\begin{equation}\label{HittingTimeLjHat+}
    H\big(\widehat L_j^+\big)
\egloi
    \tilde{A_j}\big(\widehat L_j^+\big)
    \int_{\tilde\tau_j^-[(1+\delta/2)\phi^*(t)] }^{\widehat L_j^+}
    e^{-\tilde V^{(j)}(u)} \lo_{ B}\big[\tau_{ B}(1),\tilde{A_j}(u)/\tilde{A_j}(\widehat L_j^+)\big]
    \textnormal{d}u,
\end{equation}
where
$
    \tilde{A}_j(u)
=
    \int_{\tilde m_j}^{u} e^{\tilde V^{(j)}(x)}\textnormal{d}x
$
for $u\in\R$ as before.
Since
$
    \tilde{A_j}\big(\widehat L_j^+\big)
\leq
    [\tilde \tau_j(\phi^*(t))-\tilde m_j]
    e^{\phi^*(t)}
$,
we get
$$
    \textnormal{RHS of }\eqref{HittingTimeLjHat+}
\leq
    e^{\phi^*(t)}
    \big[\tilde \tau_j(\phi^*(t))-\tilde\tau_j^-[(1+\delta/2)\phi^*(t)] \big]^2
    \sup\nolimits_{x\in\R}\lo_{ B}\big[\tau_{ B}(1),x\big].
$$
We know that
$
    P\big[\tilde \tau_j(\phi^*(t))-\tilde m_j\geq e^{\delta \phi^*(t)/3}
    \big]
\leq
    C_+ e^{- \kappa \phi^*(t)/(2\sqrt{2})}
$
and
$
    P\big[\tilde m_j - \tilde \tau_j^-[(1+\delta/2)\phi^*(t)]\geq e^{\delta \phi^*(t)/3}
    \big]
\leq
    C_+ e^{- \kappa \phi^*(t)/(2\sqrt{2})}
$
by \eqref{bessel} and \eqref{eqLongueurTempsAtteinte-}.
So,
we get with probability at least $1-C_+ e^{- \kappa \phi^*(t)/(2\sqrt{2})}$,
$
    \tilde \tau_j(\phi^*(t))-\tilde\tau_j^-[(1+\delta/2)\phi^*(t)
\leq 2 e^{\delta \phi^*(t)/3}
$ and then
$$
    P^{\wk}_{\mt_j}\big(\textnormal{RHS of }\eqref{HittingTimeLjHat+} \geq e^{ (1+ \delta)\phi^*(t)}\big)
\leq
    \P^{\wk}_{\mt_j}\big(\sup\nolimits_{x\in\R} \mathcal{L}_B(\tau^B(1), x) \geq e^{\delta \phi^*(t)/3}/4 \big)
\leq
    C_+ e^{-\delta  \phi^*(t)/3},
$$
where the last inequality comes from
\eqref{Diel}, \eqref{Diel2} and the independence of $B$ and $\wk$.
This, \eqref{InegProbaSortirPetitVoisinageGauche}, \eqref{HittingTimeLjHat+} and
${Q}^{W_{\kappa}}_{\mt_j}(X\in.)=\P^{\wk}_{\mt_j}(X\in.)$
give with a probability larger than $1-e^{-c_- \delta  \phi^*(t)}$,
\begin{equation}\label{HittingTimeL+jHat}
    Q^{\wk}_{\mt_j} \left(H(\widehat L_j^+)> e^{ (1+ \delta)\phi^*(t)}  \right)
=
    P^{\wk}_{\mt_j} \left(H(\widehat L_j^+)> e^{ (1+ \delta)\phi^*(t)}  \right)
\leq
    C_+  \phi^*(t) e^{- \delta  \phi^*(t)/8}.
\end{equation}
We get the same result for $H(\widehat L_j^-)$, since the law of $V^{(j)}$ restricted to $[\tau_j^-(h_t), \tau_j(h_t)]$
is symmetric with respect to $m_j$ for $j\geq 2$
by Fact \ref{Fact_Williams};  the result from $j=1$ follows from the fact that the valleys
for $j=1$ and $j=2$
have the same law by Lemma \ref{CVs}.
This together with \eqref{HittingTimeL+jHat} gives \eqref{6.47}.
\hfill$\Box$


\noindent{\bf Proof of Lemma \ref{lemma6.3}:}
Let $1\leq j\leq n_t-1$ and assume $0<\zeta<\k/8$.
Recall that
$
    \widehat L_j^-
=
    \tilde\tau_j^-(\phi^*(t) )
$
and
$
    \widehat L_j^+
=
    \tilde\tau_j(\phi^*(t) )
$.

First by \eqref{bessel}
since $\phi^*(t)=o(\log t)$ and $\zeta<\k/8$, we have for large $t$,
$$
    P\big[\widehat L_j^+\geq \tilde m_j+\phi^*(t)/\zeta\big]
\leq
    P\big[\tilde\tau_j(\phi^*(t)) -\tilde m_j>8\phi^*(t)/\k
    \big]
\leq
    C_+e^{-\k \phi^*(t)/(2\sqrt{2})}.
$$
Similarly by \eqref{eqLongueurTempsAtteinte-},
$
    P\big[\widehat L_j^-\leq \tilde m_j-\phi^*(t)/\zeta\big]
\leq
    C_+e^{-\k \phi^*(t)/(2\sqrt{2})}
$.
So with probability at least $1-C_+e^{-\k \phi^*(t)/(2\sqrt{2})}$,
$
    \big[\widehat L_j^-, \widehat L_j^+\big]
\subset
    [\mt_{j}-  \phi^*(t)/ \g, \mt_{j}+ \phi^*(t)/ \g ]
=
    \mI_j
$
as in Figure \ref{fig2},
and then
$
    \tilde \mu_j\big(\big[\Lm, \tL\big] \smallsetminus \mI_{j}\big)
\leq
    \tilde  \mu_j\big(\big[\Lm,\widehat L_j^-\big]\big)+\tilde \mu_j\big(\big[\widehat L_j^+, \tL\big]\big)
$.
This gives the first inequality of \eqref{IM1} with such probability.

We now prove the second inequality of \eqref{IM1}.
First, we observe that for large $t$,
$$
    \tilde{Z}
:=
    \int_{\Lm}^{ \tL} \frac{\dd y}{e^{\tilde V^{(j)}(y)}}
\geq
    \int_{\tau_j[\alpha \phi^*(t)/2]}^{\tau_j[\alpha \phi^*(t)]} \frac{\dd y}{e^{V^{(j)}(y)}}
\geq
    \big[\tau_j[\alpha \phi^*(t)]-\tau_j[\alpha \phi^*(t)/2]\big]e^{-\alpha \phi^*(t)}
\geq
    e^{-\alpha \phi^*(t)}
$$
on
$    \B_1^{\ref{lemma6.3}}\cap \mV_t$,
where
$
    \B_1^{\ref{lemma6.3}}
:=
    \{\tau_j[\alpha \phi^*(t)]-\tau_j[\alpha \phi^*(t)/2]  \geq 1 \}
$
and
$\alpha:= (1- \delta) \kappa/16 $.
By \eqref{3.10b} and Fact \ref{Fact_Williams} {\bf (ii)},
$P( \overline{\B_1^{\ref{lemma6.3}}}) \leq  4e^{-\alpha^2  \phi^*(t)^2/12}$.
Moreover, due to the definition \eqref{def_mu} of $\tilde \mu_j$, we have
\begin{equation}\label{InegMujDecomposition}
    \tilde  \mu_j\big(\big[\Lm,\widehat L_j^-\big]\big)
    +\tilde \mu_j\big(\big[\widehat L_j^+, \tL\big]\big)
=
    [{\mathcal{J}}_3+{\mathcal{J}}_4+{\mathcal{J}}_5+{\mathcal{J}}_6]/\tilde{Z},
\end{equation}
where
$$
    {\mathcal{J}}_3
:=
    \int_{\Lm}^{\tilde \tau_j^-(h_t)} \frac{\dd y}{e^{\tilde V^{(j)}(y)}},
\
    {\mathcal{J}}_4
:=
    \int_{\tilde \tau_j^-(h_t)}^{\tilde \tau_j^-[\phi^*(t)]} \frac{\dd y}{e^{\tilde V^{(j)}(y)}},
\
    {\mathcal{J}}_5
:=
    \int_{\tilde \tau_j[\phi^*(t)]}^{\tilde \tau_j(h_t)} \frac{\dd y}{e^{\tilde V^{(j)}(y)}},
\
    {\mathcal{J}}_6
:=
    \int_{\tilde \tau_j(h_t)}^{ \tL} \frac{\dd y}{e^{\tilde V^{(j)}(y)}}.
$$
Recalling that $\tilde V^{(j)}=V^{(j)}$ on $\mV_t$
by Remark \ref{RemEgaliteAvecouSansTilde}, we introduce $\gamma=(1-\delta)\k/8 $,
$$
    \B_2^{\ref{lemma6.3}}
 :=
    \big\{\inf\{ V^{(j)}(s),\ \tau_j(\phi^*(t))\leq s\leq \tau_j(h_t)\}
    >\gamma \phi^*(t)\big\},
$$
$$
    \B_3^{\ref{lemma6.3}}
:=
    \big\{\tilde \tau_j(h_t)-\tilde m_j \leq 8 h_t/ \kappa \big\}.
\qquad
    \B_4^{\ref{lemma6.3}}
:=
    \big\{\tilde L_j -\tilde \tau_j(h_t) \leq 2 h_t /\kappa\big\}.
$$
Equation \eqref{3.14b} with $h=h_t$, $\a=1/2$ and $\o=1$ gives
for large $t$ since $\tilde \tau_j(h_t)$ is a stopping time,
$$
    P\big( \overline{\B_4^{\ref{lemma6.3}}}\big)
=
    P\big(\tau^{W_{\kappa} } (-h_t/2) > 2 h_t /\kappa \big)
\leq
    e^{-\kappa h_t/16}.
$$
We have
$
    P\big(\overline{\B_3^{\ref{lemma6.3}}}
    \big)
\leq
        C_+ e^{- \kappa h_t /(2\sqrt{2})}
$
by \eqref{bessel}.
Moreover using Fact \ref{Fact_Williams} {\bf (ii)} if $j\geq 2$ and
taking the limit as $t\to+\infty$ in \eqref{preuve3.2} applied with $h=\phi^*(t)$,
$\alpha =1$, $\gamma=(1-\delta)\k/8$ and $\omega=h_t/\phi^*(t)$ gives
$
    P(\overline{\B_2^{\ref{lemma6.3}}})
\leq
    2e^{\k(\gamma-1)\phi^*(t)}
\leq
    2e^{-\delta \k^2\phi^*(t)/8}
$.
We have on $\B_2^{\ref{lemma6.3}}\cap \B_3^{\ref{lemma6.3}}\cap \B_4^{\ref{lemma6.3}}\cap \mV_t$,
$$
    {\mathcal{J}}_5
\leq
    \big[\tilde \tau_j(h_t)-\tilde m_j\big]e^{-\gamma \phi^*(t)}
\leq
    8\k^{-1} h_t e^{-\gamma \phi^*(t)},
\qquad
    {\mathcal{J}}_6
\leq
    \big[\tilde L_j -\tilde \tau_j(h_t)\big]e^{-h_t/2}
\leq
    2 \k^{-1} h_t e^{-h_t/2}.
$$
We prove similarly that there exists an event $\B_5^{\ref{lemma6.3}}$
such that
$
    P\big(\overline{\B_5^{\ref{lemma6.3}}}\cap \mV_t\big)
\leq
    2e^{\k(\gamma-1)\phi^*(t)}+C_+ e^{- \kappa h_t /(2\sqrt{2})}
$
and
$
    {\mathcal{J}}_4
\leq
    8\k^{-1} h_t e^{-\gamma \phi^*(t)}
$
on $\B_5^{\ref{lemma6.3}}\cap \mV_t$.
Furthermore, by Lemma \ref{LemmaProbasSansFaggionato} equations \eqref{eqLemmaSansFaggio1} and \eqref{eqLemmaSansFaggio3},
on some event $\B_6^{\ref{lemma6.3}}$ which has probability
at least $1-e^{-\k h_t/8}\geq 1-C_+ e^{-(c_-)\delta \phi^*(t)}$,
$$
    {\mathcal{J}}_3
\leq
    (\tilde\tau_j^-(h_t)-\tilde L_j^-)
    e^{-h_t/2}
\leq
    40 \k^{-1} h_t^+
    e^{-h_t/2}
\leq
    8\k^{-1} h_t e^{-\gamma \phi^*(t)}
$$
for large $t$.
These inequalities combined with \eqref{InegMujDecomposition} and Lemma \ref{CVs} give
on $\cap_{i=1}^6 \B_i^{\ref{lemma6.3}}\cap \mV_t$,
$$
    \tilde  \mu_j\big(\big[\Lm,\widehat L_j^-\big]\big)
    +\tilde \mu_j\big(\big[\widehat L_j^+, \tL\big]\big)
\leq
    8\k^{-1} h_t
    [3e^{-\gamma \phi^*(t)}+e^{-h_t/2}]
    e^{\alpha \phi^*(t)}
\leq
    C_+ h_t e^{-(1-\delta)\kappa\phi^*(t)/16},
$$
since $\phi^*(t)=o(\log t)$.
Since
$
    P(\cap_{i=1}^{6} \B_i^{\ref{lemma6.3}}\cap \mV_t)
\geq
    1-C_+ e^{-(c_-)\delta  \phi^*(t)}
$,
due to the previous inequalities and to Lemma \ref{CVs}, this proves the second inequality of \eqref{IM1}
for $2\leq j\leq n_t-1$.
This is also true if $j=1$ since the first valley has the same law by Lemma \ref{CVs}.
\hfill $\Box$

\medskip



\noindent{\bf Proof of Lemma \ref{LemmaProbaHitttingTimesHjPlusMoins}:}
Let $j\in[1, n_t)$.
By \eqref{eqHittingTimePartantdemi} applied with
$i=j$ and $r=\tilde \tau_j[\k \phi^*(t)/(8\zeta)]-\tilde m_j$,
there exists a Brownian motion
$\tilde B$, independent of $\tilde V^{(j)}$, such that under $\P_{\mt_j}^{\wk}$,
$$
    H\Big(
        \tilde \tau_j[\k \phi^*(t)/(8\zeta)]
     \Big)
\geq
    \int^{\tilde \tau_j[\k \phi^*(t)/(8\zeta)]}_{ \mt_j} e^{-V^{(j)}(u)}
    \lo_{\tilde B}\Big[\tau^{\tilde B}\Big(\tilde{A}_j \Big(\tilde \tau_j\Big[ \frac{\k\phi^*(t)}{8\zeta}\Big]\Big)\Big),\tilde{A}_j(u)\Big]
    \textnormal{d}u
=:
    H^+_j,
$$
where for all  $z \in \R$, $\tilde{A}_j(z)=\int_{ \mt_j}^{z} e^{\tilde V^{(j)}(x)}\textnormal{d}x$.
Now, let
$
    \B_7^{\ref{lemma6.3}}
:=
    \{\tilde \tau_j[\k \phi^*(t)/(8\zeta)]-\mt_j\leq \phi^*(t)/\zeta\}
$.
By \eqref{bessel},
$
    P\big(\B_7^{\ref{lemma6.3}}\big)
\geq
    1-C_+e^{-\k^2\phi^*(t)/(16\zeta\sqrt{2})}
$.
We have on $\B_7^{\ref{lemma6.3}}$ under $\P_{\mt_j}^{\wk}$,
$$
    H\big(\mt_{j}+\phi^*(t)/\g\big)
\geq
    H\big(\tilde \tau_j[\k \phi^*(t)/(8\zeta)]\big)
\geq
    H^+_j.
$$
Assume  $\zeta\leq \k/48$.
In order to estimate $H_j^+$, we introduce
$$
    \delta_t^*
:=
    e^{-\k\phi^*(t)/(48\zeta)}\in(0,1),
\qquad
    \mathcal{J}_7
:=
   \tilde  A_j\Big(\tilde \tau_j\Big[\frac{\k \phi^*(t)}{8\zeta}\Big]\Big),
\qquad
    \mathcal{J}_8
:=
    \int_{\tilde m_j}^{\tilde \tau_j[\frac{\k \phi^*(t)}{16\zeta}]}
            e^{-\tilde V^{(j)}(u)}
            \dd u.
$$
By scaling,
there exists some Brownian motion $B'$ independent of $\tilde V^{(j)}$ such that
$$
    H_j^+
\geq
    \int_{\tilde m_j}^{\tilde \tau_j[\frac{\k \phi^*(t)}{16\zeta}]} e^{-\tilde V^{(j)}(u)}
            \mathcal{J}_7
            \lo_{B'}\big(\tau^{B'}(1), \tilde A_j(u)/\mathcal{J}_7\big)
            \dd u
\geq
            \mathcal{J}_7
            \mathcal{J}_8
            \mathcal{J}_9
$$
on
$
    \B_8^{\ref{lemma6.3}}
:=
    \{
           \tilde  A_j\big(\tilde \tau_j\big[\frac{\k \phi^*(t)}{16\zeta}\big]\big)
        \leq
            \delta_t^* \mathcal{J}_7
    \}
$
with
$
    \mathcal{J}_9
:=
    \inf_{x\in[0,\delta_t^*]}
    \lo_{B'}\big(\tau^{B'}(1), x\big)
$.
We have
$$
    \mathcal{J}_8
\geq
    \big[
        \tilde \tau_j[\k \phi^*(t)/(48\zeta)]
        -
        \tilde \tau_j[\k \phi^*(t)/(96\zeta)]
    \big]
    \exp[-\k \phi^*(t)/(48\zeta)]
\geq
    \exp[-\k \phi^*(t)/(48\zeta)]
$$
with probability $\geq 1-e^{-(c_-)(\phi^*(t))^2/\zeta^2}-e^{-h_t/4}$ for large $t$
by \eqref{3.10b}, Fact \ref{Fact_Williams} {\bf (ii)} and Lemma \ref{CVs}.
Moreover,
$$
    \mathcal{J}_9
\geq
    \big[1-(\delta_t^*)^{1/3}\big]
    \lo_{B'}\big(\tau^{B'}(1), 0\big)
\geq
    (1/2)
    \lo_{B'}\big(\tau^{B'}(1), 0\big)
\geq
    e^{-\k\phi^*(t)/(48\zeta)}
$$
for large $t$
with probability $\geq 1-c_2(\delta_t^*)^{1/30}$ by \eqref{Dev}
for the first inequality,
and with probability $\geq 1-e^{-\k\phi^*(t)/(48\zeta)}$
for the last one since
$    \lo_{B'}\big(\tau^{B'}(1), 0\big)$ is exponentially distributed with mean $2$
as before by the first Ray-Knight theorem.
Furthermore,  by Lemma \ref{CVs}, Fact \ref{Fact_Williams} and \eqref{MinorationAVallee},
$$
    P\Big(
        \mathcal{J}_7\geq e^{\frac{5\k\phi^*(t)}{48\zeta}}
     \Big)
\geq
    P\Big(
        \int_0^{\tau^R(\k\phi^*(t)/(8\zeta))}e^{R(u)}\dd u
        \geq
        e^{(1-1/6)\frac{\k\phi^*(t)}{8\zeta}}
     \Big)
    -
    P\big(\overline{\mV_t}\big)
\geq
    1-4e^{-\frac{(c_-)\phi^*(t)}{\zeta}}.
$$
Finally by \eqref{bessel},
with probability at least $1-C_+ \exp(-\frac{\k^2\phi^*(t)}{32\zeta\sqrt{2}})$,
$$
    \tilde A_j\Big( \tilde \tau_j\Big[\frac{\k\phi^*(t)}{16\zeta} \Big]\Big)
\leq
    \Big[\tilde \tau_j\Big[\frac{\k\phi^*(t)}{16\zeta} \Big]-\tilde m_j\Big]
    \exp\Big(\frac{\k\phi^*(t)}{16\zeta} \Big)
\leq
    \frac{\phi^*(t)}{2\zeta}
    \exp\Big(\frac{\k\phi^*(t)}{16\zeta} \Big).
$$
The last two inequalities give
$
    P[    \B_8^{\ref{lemma6.3}}]
\geq
    1-e^{-(c_-)\phi^*(t)/\zeta}
$.
As a consequence, we have
$$
\P_{\mt_j}^{\wk}
\Big[
    H\big(\mt_{j}+\phi^*(t)/\g\big)
\geq
    H^+_j
\geq
            \mathcal{J}_7
            \mathcal{J}_8
            \mathcal{J}_9
\geq
    e^{\k \phi^*(t)/(16\zeta)}
\geq
    e^{(1+2\delta)\phi^*(t)}
\Big]
\geq
    1-e^{-(c_-)\phi^*(t)/\zeta}
$$
with probability at least $1-e^{-(c_-)\phi^*(t)/\zeta}$
since $\zeta\leq \k/48$ and $\delta < 1$.
We obtain the same result for $H(\tilde m_j-\phi^*(t)/\zeta)$ by symmetry of the law of $V^{(j)}$
for $j\geq 2$ by Fact \ref{Fact_Williams} {\bf (ii)}, and then for $j=1$ by Lemma \ref{CVs} as before.
\hfill$\Box$




\mysection{Proofs of some technical estimates related to the environment}\label{SectionAnnexe}

\subsection{Proof of Lemma \ref{lem4.5}}

We denote by $I_{\kappa}$ and $K_{\kappa}$ the modified Bessel functions, respectively of the first and second kind.
We remind that as $x\downarrow 0$ (see e.g. \cite{BorodinSalminem} p. 638), since $0<\k<1$,
\begin{align}
    I_\k(x)
& =  \frac{1}{\Gamma(\k+1)}\bigg(\frac{x}{2}\bigg)^\k
      +
      \frac{1}{\Gamma(\k+2)}\bigg(\frac{x}{2}\bigg)^{\k+2}
      +
      O\big(x^{\k+4}\big),
\label{eqDLBesselI}
\\
    K_\k(x)
& =
    \frac{\pi\big[I_{-\k}(x)-I_{\k}(x)\big]}{2\sin(\pi \k )}
 =
    \frac{\pi}{2\sin(\pi \k )}
    \bigg[
          \frac{(x/2)^{-\k}}{\Gamma(1-\k)}
    -
          \frac{(x/2)^\k}{\Gamma(\k+1)}
    +
          \frac{(x/2)^{2-\k}}{\Gamma(2-\k)}
    \bigg]
    +O(x^{\k+2}).
\label{eqDLBesselK}
\end{align}
Moreover, we remind that (see e.g. \cite{BorodinSalminem} p. 638),
\begin{equation}\label{eqSommeDiffBessel}
    I_\k'(u)K_\k(u)-I_\k(u)K_\k'(u)
=
    1/u,
\qquad
    u>0.
\end{equation}

Let $y>0$. First,  (\cite{BorodinSalminem}, 2.10.3 page 302) with $\alpha=0$,  $x=y$,
$z=y/2<x $, $\beta=1/2$,  and $\mu=-\k/2$
gives for $G^+$, which is defined in \eqref{eqDefF+G+},
$$
    E \left( e^{- \gamma G^+(y/2,y )}  \right)
=
    E\bigg[
        \exp\bigg(
                -\gamma \int_0^{\tau^{ W_\k^y}(y/2)}e^{ W_\k^y(s)}\dd s
            \bigg)
        \bigg]
=
    e^{\k y/4 }
    \frac
        {K_\k(2\sqrt{2\gamma}e^{y/2})}
        {K_\k(2\sqrt{2\gamma}e^{y/4})},
\qquad
    \gamma >0.
$$
So,
$
    E \left( e^{- \gamma G^+(y/2,y )/e^y } \right)
=
    e^{\k y/4 }
    \frac
        {K_\k(2\sqrt{2\gamma })}
        {K_\k(2\sqrt{2\gamma }e^{-y/4})}
=
    g_+\big(2\sqrt{2\gamma},[2\sqrt{2\gamma}e^{-y/4}]^\k\big)
$.
where
$$
    g_+(u,v)
:=
    \frac
        {u^\k K_\k(u)}
        {v K_\k(v^{1/\k})},
\qquad
    u>0,\ v>0.
$$
We have, as $\max(u,v)\downarrow 0$, by \eqref{eqDLBesselI} and \eqref{eqDLBesselK},
\begin{eqnarray*}
    g_+(u,v)
& = &
    \frac
        {u^\k
            \frac{\pi}{2\sin(\pi \k )}
            \big[
                \frac{1}{\Gamma(1-\k)}(u/2)^{-\k}
                -
                \frac{1}{\Gamma(\k+1)}(u/2)^{\k}
                +
                \frac{1}{\Gamma(2-\k)}(u/2)^{2-\k}
                +o(u^{2-\k})
            \big]
        }
        {v
            \frac{\pi}{2\sin(\pi \k )}
                \big[
                \frac{1}{\Gamma(1-\k)}\frac{v^{-1}}{2^{-\k}}
                -
                \frac{1}{\Gamma(\k+1)}\frac{v}{2^\k}
                +o(v)
                \big]
        }
\\
& =  &
                1
                -
                \frac{\Gamma(1-\k)}{\Gamma(\k+1)}\frac{u^{2\k}}{4^\k}
            +O\big([\max(u,v)]^2\big).
\end{eqnarray*}
This gives, with $u=2\sqrt{2\gamma}$ and $v=[2\sqrt{2\gamma}e^{-y/4}]^\k$
as $\gamma\downarrow 0$ and $y\to+\infty$,
\begin{equation}\label{eqDLg+}
    E \left( e^{- \gamma G^+(y/2,y )/e^y } \right)
=
               1
                -
                \frac{\Gamma(1-\k)}{\Gamma(\k+1)}(2\gamma)^{\k}
            +O\big(\max(\gamma,\gamma^\k e^{-\k y/2})\big).
\end{equation}

We now turn to $F^\pm(y)$, defined in \eqref{eqDefF+G+}.
We have for $\gamma>0$,  by \eqref{eqLienEntreRetWzNul},
\begin{equation*}
    E \left( e^{- \gamma F^\pm(y)}  \right)
 =
    \lim_{x\downarrow 0}
    E\bigg[
            \exp\bigg(-\gamma\int_0^{\tau^{W_{-\k}^x}(y)}e^{\pm W_{-\k}^x(s)}\dd s \bigg)
            \frac{\un_{\{\tau^{W_{-\k}^x}(0)=\infty\}}}{P\big[\tau^{W_{-\k}^x}(0)=\infty\big]}
    \bigg].
\end{equation*}
The expectation in the right hand side of this equality is equal to,
first by the strong Markov property,
and second by
(\cite{BorodinSalminem}, 3.10.7(b) page 317) with $\alpha=0$, $a=0$, $b=y$, $\beta=\pm 1/2$,
$\mu=\k/2$ and $x>0$,
and since $P\big(\tau^{W_{-\k}^x}(0)=\infty\big)=1-e^{-\k x}$ due to the scale function \eqref{eqScaleFunctionWk},
\begin{eqnarray}
&&
    E\bigg[
            \exp\bigg(-\gamma\int_0^{\tau^{W_{-\k}^x}(y)}e^{\pm W_{-\k}^x(s)}\dd s \bigg)
            \un_{\{\tau^{W_{-\k}^x}(y)<\tau^{W_{-\k}^x}(0)\}}
    \bigg]
    \frac{P\big[\tau^{W_{-\k}^y}(0)=\infty\big]}
         {P\big[\tau^{W_{-\k}^x}(0)=\infty\big]}
\nonumber\\
& = &
    \frac{e^{(\mu-|\mu|(\pm 1))y}
            S_{\k}\big(2\sqrt{2\gamma}e^{\pm  x/2},2\sqrt{2\gamma}\big)}
         {e^{(\mu-|\mu|(\pm 1))x}
            S_{\k}\big(2\sqrt{2\gamma}e^{\pm y/2},2\sqrt{2\gamma}\big)}
    \frac{1-e^{-\k y}}{1-e^{-\k x}},
\nonumber
\end{eqnarray}
where
$
    S_\k(u,v)
=
    (u v)^{-\k}[I_\k(u)K_\k(v)-K_\k(u)I_\k(v)]
    $ as defined in (\cite{BorodinSalminem} p. 645).
So
$E \big( e^{- \gamma F^\pm(y)}  \big)$ is the limit, as $x\downarrow 0$, of
\begin{eqnarray*}
    \frac{
            \big[
                I_\k(2\sqrt{2\gamma}e^{\pm x/2})K_\k(2\sqrt{2\gamma})
                -K_\k(2\sqrt{2\gamma}e^{\pm x/2})I_\k(2\sqrt{2\gamma})
            \big]}
         {
            \big[
                I_\k(2\sqrt{2\gamma}e^{\pm y/2})K_\k(2\sqrt{2\gamma})
                -K_\k(2\sqrt{2\gamma}e^{\pm y/2})I_\k(2\sqrt{2\gamma})
            \big]
         }
    \frac{\sinh(\k y/2)}
         {\sinh(\k x/2)}
=:
    \frac{N_\pm(x,y)}{D_\pm(x,y)}.
\end{eqnarray*}
Now, notice that as $x\downarrow 0$,
\begin{align*}
    I_\k\big(2\sqrt{2\gamma}e^{\pm x/2}\big)
& =
    I_\k\big[2\sqrt{2\gamma}(1\pm x/2+o(x))\big]
=
    I_\k\big[2\sqrt{2\gamma}\big]
    \pm
    \sqrt{2\gamma}I_\k'[2\sqrt{2\gamma}] x
    +
    o(x),
\\
    K_\k\big(2\sqrt{2\gamma}e^{\pm x/2}\big)
&=
    K_\k\big[2\sqrt{2\gamma}(1\pm x/2+o(x))\big]
=
    K_\k\big[2\sqrt{2\gamma}\big]
    \pm
    \sqrt{2\gamma}K_\k'[2\sqrt{2\gamma}] x
    +
    o(x).
\end{align*}
So by \eqref{eqSommeDiffBessel},
\begin{eqnarray*}
    N_\pm(x,y)
& \sim_{x\downarrow 0}  &
    \pm\sqrt{2\gamma}
        \big[
            I_\k'[2\sqrt{2\gamma}]     K_\k(2\sqrt{2\gamma})
            -
            K_\k'[2\sqrt{2\gamma}]     I_\k(2\sqrt{2\gamma})\big]
    \sinh(\k y/2)x
\\
& \sim_{x\downarrow 0} &
    \pm\sinh(\k y/2) x/2.
\end{eqnarray*}
Moreover, $\sinh(\k x/2)\sim_{\downarrow 0}\k x/2$, and then
\begin{equation}\label{eqLaplaceFpm}
    E \left( e^{- \gamma F^\pm(y)}  \right)
=
    \lim_{x\downarrow 0}
    \frac{N_\pm(x,y)}{D_\pm(x,y)}
=
    \frac{\pm \k^{-1}\sinh(\k y/2) }
         {
                I_\k(2\sqrt{2\gamma}e^{\pm y/2})K_\k(2\sqrt{2\gamma})
                -K_\k(2\sqrt{2\gamma}e^{\pm y/2})I_\k(2\sqrt{2\gamma})
         }.
\end{equation}
We now consider
\begin{equation}\label{JoliDL-}
    f(u,v)
:=
    \frac{-(v^{-\k}-v^\k)}{2\k[I_\k(u v)K_\k(u)-K_\k(u v)I_\k(u)]}
\end{equation}
so that
$
    E \big( e^{- \gamma F^-(y)}  \big)
=
    E \big( e^{- \gamma F^+(y)/e^y}  \big)
=
    f(2\sqrt{2\gamma},e^{-y/2})
$,
and so $F^-(y) \egloi F^+(y)/e^y$.
We get successively, as $\max(v^{2\k},u^3)\to 0$, by \eqref{eqDLBesselI} and \eqref{eqDLBesselK},
using $\Gamma(1-\k)\Gamma(\k)=\pi/\sin(\pi \k)$,
\begin{eqnarray*}
    K_\k(uv)I_\k(u)
& = &
    \frac{1}{2\k}v ^{-\k}
    +
    \frac{1}{8(\k+1)\k}v^{-\k}u^{2}
    +
    v^{-\k}o(\max(v^{2\k},u^3)),
\\
    I_\k(u v)K_\k(u)
& = &
    v^\k/(2\k)
    +
    v^{-\k}o(\max(v^{2\k},u^3)),
\\
    2\k\big[
        I_\k(u v)K_\k(u)
        -
        K_\k(uv)I_\k(u)
    \big]
& = &
    v^\k-v^{-\k}-\frac{1}{4(\k+1)}v^{-\k}u^2
+
    v^{-\k}o(\max(v^{2\k},u^3)).
\end{eqnarray*}
This yields
\begin{equation*}
    f(u,v)
 =
    \frac{-(v^{-\k}-v^\k)}
    {
        v^\k-v^{-\k}-\frac{1}{4(\k+1)}v^{-\k}u^2
        +
        v^{-\k}o(\max(v^\k,u^3))
    }
 =
        1-\frac{1}{4(\k+1)}u^2
        +
        o(\max(v^{2\k},u^3)).
\end{equation*}
Consequently, as $\max(e^{-\k y}, \gamma^{3/2})\to 0$,
\begin{equation*}
    f(2\sqrt{2\gamma},e^{-y/2})
=
    1-\frac{2\gamma}{\k+1}
    +
    O(\max(e^{-\k y},\gamma^{3/2}))
=
    \bigg(1+\frac{2\gamma}{\k+1}\bigg)^{-1}
    +
    O(\max(e^{-\k y },\gamma^{3/2})).
\end{equation*}
Since
$
    E \big( e^{- \gamma F^+(y)/e^y}  \big)
 =
    E \big( e^{- \gamma F^-(y)}  \big)
=
    f(2\sqrt{2\gamma},e^{-y/2})
$,
this and \eqref{eqDLg+} proves the existence of
$C_4>0$, $M>0$ and $\eta_1\in(0,1)$
such that
\eqref{Fmoins}, \eqref{Ff2} and \eqref{Gpf2} are satisfied simultaneously
for all $y>M$ and $\gamma\in(0, \eta_1]$.

We now prove $E[{F^+(y)}/{e^y}]$ exists and is bounded,
by computing
$\frac{\partial}{\partial \gamma}E \big( e^{- \gamma F^+(y)/e^y}  \big)$
at $\gamma=0$. To this aim, we fix $y>0$ and observe that as $\gamma\downarrow 0$,
once more by \eqref{eqDLBesselI} and \eqref{eqDLBesselK},
\begin{equation*}
    E \left( e^{- \gamma F^+(y)/e^y}  \right)
=
    f(2\sqrt{2\gamma},e^{-y/2})
=
            1
            -
            \frac{[e^{\k y /2}-e^{-\k y/2-y}]}{(\k+1)\sinh(\k y/2)}
          \gamma
            -\frac{(e^{-y+\k y/2}-e^{-\k y/2})}{(1-\k)\sinh(\k y/2)}\gamma
            +
            o(\gamma).
\end{equation*}
Hence,
$$
    E(F^+(y)/e^y)
=
    -\left(\frac{\partial}{\partial \gamma}\right)_{\gamma=0}E \left( e^{- \gamma F^+(y)/e^y}  \right)
=
            \frac{[e^{\k y /2}-e^{-\k y/2-y}]}{(\k+1)\sinh(\k y/2)}
            +\frac{(e^{-y+\k y/2}-e^{-\k y/2})}{(1-\k)\sinh(\k y/2)},
$$
which is a bounded function of $y$ on $\R_+$.

Finally, taking the limit of $E \big( e^{- \gamma F^-(y)}  \big)$
in \eqref{eqLaplaceFpm} as $y\to +\infty$ with the help of \eqref{eqDLBesselI} and
\eqref{eqDLBesselK} proves \eqref{eqLimiteLaplaceF-}.
\hfill$\Box$

\subsection{Proof of Lemma \ref{LemmaProbasSansFaggionato}}\label{SubScetPreuveLemmasTechniquesht+}

We consider
$
    \B_1^{\ref{SubScetPreuveLemmasTechniquesht+}}
:=
    \cap_{i=1}^{n_t}\{\tilde L_i^+-\tilde \tau_i^-(h_t^+) \leq  40 h_t^+/\k\}
$.
We also introduce
$
    \B_2^{\ref{SubScetPreuveLemmasTechniquesht+}}
:=
    \{\wk(\tilde m_{n_t}) \geq -n_t e^{5\k h_t/4}\}
$.
We have,
\begin{equation*}
    P\big[\overline{\B_2^{\ref{SubScetPreuveLemmasTechniquesht+}}}\big]
=
    P\big[\wk(\tilde m_{n_t})< -n_t e^{5\k h_t/4}\big]
=
    P\bigg(\sum_{i=1}^{n_t}\big(2h_t^++\wk(\tilde L_i^\sharp)-\wk(\tilde m_i)\big)> n_t e^{5\k h_t/4}\bigg).
\end{equation*}
Recall that the r.v. $\wk(\tilde L_i^\sharp)-\wk(\tilde m_i)$, $i\geq 1$ are
equal in law to
$-\inf_{[0,\tau_1^*(h_t)]}\wk$, which
expectancy is
$\frac{2}{\k} \sinh(\k h_t/2)e^{\k h_t/2}$ (as recalled before \eqref{InegBeta*}).
By Markov inequality, for large $t$,
\begin{equation}\label{eqProbaB4SansFaggio}
    P\big[\overline{\B_2^{\ref{SubScetPreuveLemmasTechniquesht+}}}\big]
\leq
    \frac{n_t}{n_t e^{5\k h_t/4}}\Big(2h_t^++\frac{2}{\k} \sinh(\k h_t/2)e^{\k h_t/2}\Big)
\leq
    C_+ e^{-\k h_t/4}.
\end{equation}


On $\overline{\B_1^{\ref{SubScetPreuveLemmasTechniquesht+}}}$ there exists $1\leq i\leq n_t$ such that
$
    \tilde L_i^+-\tilde \tau_i^-(h_t^+)
>
    40 h_t^+/\k
$.
There exists an integer $j\in\Z$ such that $-j\geq \wk[\tilde \tau_i^-(h_t^+)]>-j-1$.
So on $\B_2^{\ref{SubScetPreuveLemmasTechniquesht+}}$ by \eqref{InegLiPremieresDescentes},
$-j>\wk(\tilde m_i)\geq \wk(\tilde m_{n_t})\geq -n_t e^{5\k h_t/4}$,
thus $j< n_t e^{5\k h_t/4}$.
Moreover, we have by \eqref{InegLiPremieresDescentes},
\begin{equation}\label{eqPourDeterminerValeursdej}
    \wk[\tilde \tau_i^-(h_t^+)]
=
    \wk[\tilde m_i]+h_t^+
\leq
    \wk(\tilde L_i^\sharp)+h_t^+
=
    \wk(\tilde L_{i-1}^+)
\leq
    0,
\end{equation}
so $j\geq 0$.
Since
$0\geq -j\geq \wk[\tilde \tau_i^-(h_t^+)]$, we have
$\tau^{\wk}(-j)\leq \tilde \tau_i^-(h_t^+)$.
On the other hand,
$
    -j-1
<
    \wk[\tilde \tau_i^-(h_t^+)]
=
    \wk(\tilde m_i)+ h_t^+
$, so by \eqref{InegLiPremieresDescentes},
$$
    \wk(\tilde L_i^+)
=
    \inf_{[0, \tilde L_i^+]} \wk
=
    \wk(\tilde m_i)-h_t^+
>
    -j-1-2h_t^+,
$$
and then
$
    \tau^{\wk}[-j-1-2h_t^+]
\geq
    \tilde L_i^+
$.
Hence,
$
    \tau^{\wk}[-j-1-2h_t^+]-\tau^{\wk}(-j)
\geq
    \tilde L_i^+-\tilde \tau_i^-(h_t^+)
>
    40 h_t^+/\k
$.
So for large $t$,
\begin{eqnarray*}
     P\big(\overline{\B_1^{\ref{SubScetPreuveLemmasTechniquesht+}}}\cap \B_2^{\ref{SubScetPreuveLemmasTechniquesht+}}\big)
& \leq &
    P\Big(\cup_{j=0}^{\lfloor n_t e^{5\k h_t/4}\rfloor}
            \big\{
                \tau^{\wk}[-j-1-2h_t^+]-\tau^{\wk}(-j)
                \geq
                40h_t^+/\k
            \big\}
     \Big)
\\
& \leq &
    (n_t e^{5\k h_t/4}+1)
    P\Big(
         \tau^{\wk}[-1-2h_t^+]
         \geq
         40 h_t^+/\k
     \Big)
\\
& \leq &
    2n_t e^{5\k h_t/4}
    P\Big(
         \tau^{\wk}[-3 h_t^+]
         \geq
         40 h_t^+/\k
     \Big)
 \leq
    2 n_t e^{5\k h_t/4}
    \exp\big[-3 \k h_t^+\big]
 \leq
    e^{-\k h_t},
\end{eqnarray*}
by Lemma \ref{lemmaProbaAtteintewk} with $h=h_t^+$, $\alpha=3$ and $\omega=20$
and since $n_t=e^{o(h_t)}$ because $\phi(t)=o(\log t)$.

We also consider
$$
    \B_3^{\ref{SubScetPreuveLemmasTechniquesht+}}
:=
    \cap_{i=1}^{n_t}
            \big\{
                \inf\nolimits_{[\tilde \tau_i^-(h_t^+), \tilde \tau_i^-(h_t^+)+1]} \tilde V^{(i)}
                \geq
                h_t^+-\k h_t/2
            \big\},
\qquad
    \B_4^{\ref{SubScetPreuveLemmasTechniquesht+}}
:=
    \{\tilde m_{n_t}\leq e^{2\k h_t}\}.
$$
We notice that since $E[\tau^{\wk}(-1)]<\infty$,
$$
    E\big[\tilde L_1^+-\tilde \tau_1(h_t)\big]
=
    E\big[\tau^{\wk}(-h_t^+-h_t)\big]
\leq
    (\lfloor h_t\rfloor +1)E\big[\tau^{\wk}(-1)\big]
\leq
    C_+ h_t.
$$
Similarly,
$
    E\big[\tilde L_1^\sharp\big]
=
    E\big[\tau^{\wk}(-h_t^+)\big]
\leq
    C_+ h_t
$.
Moreover,
$
    E\big[\tilde\tau_1(h_t)-\tilde L_1^\sharp\big]
=
    E\big[\tau_1^*(h_t)\big]
\leq
    e^{\k h_t}
$
by \eqref{InegProbaTau1Etolie}. Combining these inequalities gives
$
    E\big[\tilde L_1^+\big]
\leq
    C_+ e^{\k h_t}
$, and then
$
    E[\tilde m_{n_t}]
\leq
    E\big[\tilde L_{n_t}^+\big]
=
    n_t E\big[\tilde L_{1}^+\big]
\leq
    C_+ n_t e^{\k h_t}
$ since $\big(\tilde L_{i+1}^+-\tilde L_i^+\big)$, $i\geq 0$, are i.i.d.
Consequently for large $t$,
$$
    P\big[\overline{\B_4^{\ref{SubScetPreuveLemmasTechniquesht+}}}\big]
\leq
    E[\tilde m_{n_t}]/e^{2\k h_t}
\leq
    e^{-\k h_t/2}.
$$

On $\overline{\B_3^{\ref{SubScetPreuveLemmasTechniquesht+}}}$, there exists $1\leq i\leq n_t$
such that
$
    \inf_{[\tilde \tau_i^-(h_t^+), \tilde \tau_i^-(h_t^+)+1]}\tilde V^{(i)}
<
     h_t^+-\k h_t/2
$.
Since $\tilde V^{(i)}[\tilde \tau_i^-(h_t^+)]=h_t^+$,
for $k=\lfloor \tilde \tau_i^-(h_t^+)\rfloor$, we have
$
    \sup_{[k,k+2]} \wk - \inf_{[k,k+2]} \wk
\geq
    \k h_t/2
$,
and $0\leq k\leq \tilde m_{n_t}\leq e^{2\k h_t}$ on $\B_4^{\ref{SubScetPreuveLemmasTechniquesht+}}$.
Moreover we have for every $k\in\N$,
\begin{eqnarray*}
    P\Big[
        \sup_{[k,k+2]} \wk - \inf_{[k,k+2]} \wk
        \geq
        \k h_t/2
     \Big]
& \leq &
    P\Big[
        \sup_{[0,2]} W
        \geq
        \frac{\k h_t}{4}-\k
     \Big]
    +
    P\Big[
        -\inf_{[0,2]} W
        \geq
        \frac{\k h_t}{4}-\k
     \Big]
\\
& = &
    2P\big[|W(2)|\geq \k h_t/4-\k\big]
 \leq
    4\exp\big[-(\k h_t/4-\kappa)^2/4\big]
\end{eqnarray*}
for large $t$ since $\sup_{[0,2]} W\egloi |W(2)|$ and $P(W(1)\geq x)\leq e^{-x^2/2}$ for large $x$. Consequently,
$$
    P\big(\overline{\B_3^{\ref{SubScetPreuveLemmasTechniquesht+}}}\big)
\leq
    P\big(\overline{\B_3^{\ref{SubScetPreuveLemmasTechniquesht+}}}
    \cap \B_4^{\ref{SubScetPreuveLemmasTechniquesht+}}\big)
    +
    P\big(\overline{\B_4^{\ref{SubScetPreuveLemmasTechniquesht+}}}\big)
\leq
    \sum_{k=0}^{\lfloor e^{2\k h_t} \rfloor}
    P\Big[
        \sup_{[k,k+2]} \wk - \inf_{[k,k+2]} \wk
        \geq
        \frac{\k h_t}{2}
     \Big]
    +
    P\big(\overline{\B_4^{\ref{SubScetPreuveLemmasTechniquesht+}}}\big)
\leq
    e^{-\k h_t/4}
$$
for large $t$.
Notice in particular that on $\B_3^{\ref{SubScetPreuveLemmasTechniquesht+}}$,
we have for every $1\leq i \leq n_t$,
$\tilde \tau_i^-(h_t^+)+1<\tilde m_i$ since $\tilde V^{(i)}(\tilde  m_i)=0<h_t^+-\k h_t/2$, and
\begin{equation}\label{eqIntegraleSansFaggio}
    \int_{\tilde L_i^-}^{\tilde m_i} e^{\tilde V^{(i)}(u)}\dd u
\geq
    \int_{\tilde \tau_i^-(h_t^+)}^{\tilde \tau_i^-(h_t^+)+1} e^{\tilde V^{(i)}(u)}\dd u
\geq
    e^{h_t^+-\k h_t/2}.
\end{equation}

Finally, let
$
    \B_5^{\ref{SubScetPreuveLemmasTechniquesht+}}
:=
    \cap_{i=1}^{n_t}
    \big\{
        \inf_{[\tilde \tau_i^-(h_t^+), \tilde \tau_i^-(h_t)]}\tilde V^{(i)}
        \geq
        h_t/2
    \big\}
$.
On $\overline{\B_5^{\ref{SubScetPreuveLemmasTechniquesht+}}}$, we consider $1\leq i \leq n_t$ such that
$
    \inf_{[\tilde \tau_i^-(h_t^+), \tilde \tau_i^-(h_t)]}\tilde V^{(i)}
    <
    h_t/2
$.
There exists an integer $j\in\Z$ such that $-j\geq \wk[\tilde \tau_i^-(h_t^+)]>-j-1$.
As before, just before and after \eqref{eqPourDeterminerValeursdej},
this yields $0\leq j<n_t e^{5\k h_t/4} $ on $\B_2^{\ref{SubScetPreuveLemmasTechniquesht+}}$.

Consider now
$
    y_j
:=
    \inf\big\{x\geq \tau^{\wk}(-j-1),\ \wk(x)-\inf_{[\tau^{\wk}(-j-1), x]}\wk\geq h_t/2\big\}
$.
We have
$
    y_j
>
    \tau^{\wk}(-j-1)
\geq
    0
$.
Moreover, due to the definition  of our $i$,
and since for all $0\leq x\leq \tau^{\wk}(-j-1)$,
$
    \wk[\tau^{\wk}(x)]
\geq
    -j-1
>
    \wk(\tilde m_i)+h_t^+-1
$
and so, for large $t$,
$
    \tilde V^{(i)}(x)
>
    h_t/2
\geq
        \inf_{[\tilde \tau_i^-(h_t^+), \tilde \tau_i^-(h_t)]}\tilde V^{(i)}
$,
we have
$$
    \tilde V^{(i)}[\tilde \tau_i^-(h_t)]
    -\inf_{[\tau^{\wk}(-j-1), \tilde \tau_i^-(h_t)]}\tilde V^{(i)}
\geq
    \tilde V^{(i)}[\tilde \tau_i^-(h_t)]
    -\inf_{[\tilde \tau_i^-(h_t^+), \tilde \tau_i^-(h_t)]}\tilde V^{(i)}
\geq
    h_t/2.
$$
Thus, $y_j\leq \tilde \tau_i^-(h_t)$.

Hence,
$
    \inf_{[\tau^{\wk}(-j-1), y_j]} \wk
\geq
    \inf_{[0, \tilde \tau_i(h_t)]} \wk
=
    \wk(\tilde m_i)
$
by \eqref{InegLiPremieresDescentes}.
Thus for these  $j$ and $i$,
\begin{equation}\label{eqFinaleLemmmeSansFaggio1}
    \inf_{[\tau^{\wk}(-j-1), y_j]} \wk-(-j-1)
\geq
    \wk(\tilde m_i)-(-j-1)
\geq
    \wk(\tilde m_i)-\wk[\tilde\tau_i^-(h_t^+)]
=
    -h_t^+,
\end{equation}
\begin{equation}\label{eqFinaleLemmmeSansFaggio2}
    \inf_{[y_j, \tilde\tau_i(h_t)]}\wk-\wk(y_j)
\geq
    \wk(\tilde m_i)-[(-j-1)+h_t/2]
\geq
    -h_t^+-h_t/2.
\end{equation}
So on
$
\overline{\B_5^{\ref{SubScetPreuveLemmasTechniquesht+}}}
\cap \B_2^{\ref{SubScetPreuveLemmasTechniquesht+}}
$
there exists $0\leq j\leq \lfloor n_t e^{5\k h_t/4}\rfloor$ and some $1\leq i\leq n_t$ such that
such that
\eqref{eqFinaleLemmmeSansFaggio1} and \eqref{eqFinaleLemmmeSansFaggio2} are satisfied.
More over for this $i$, we showed that
$y_j\leq \tilde \tau_i^-(h_t)$,
one consequence of this being that $\tilde \tau_i(h_t)-y_j$ plays the role of
$\tau_1^*(h_t)$ for the process $\wk(.+y_j)-\wk(y_j)$.
Moreover $y_j-\tau^{\wk}(-j-1)$ plays the role of $\tau_1^*(h_t/2)$ for $\wk[.+\tau^{\wk}(-j-1)]-(-j-1)$.
Applying the strong Markov property for every $0\leq j\leq \lfloor  n_t e^{5\k h_t/4}\rfloor$
 at stopping times $\tau^{\wk}(-j-1)$ and $y_j$, we get for large $t$,
\begin{eqnarray*}
    P\big[
        \overline{\B_5^{\ref{SubScetPreuveLemmasTechniquesht+}}}
        \cap
        \B_2^{\ref{SubScetPreuveLemmasTechniquesht+}}
    \big]
& \leq &
    \sum_{j=0}^{\lfloor n_t e^{5\k h_t/4}\rfloor}
    P\Big[\inf_{[0,\tau_1^*(h_t/2)]}\wk \geq -h_t^+\Big]
    P\Big[\inf_{[0,\tau_1^*(h_t)]}\wk \geq -h_t^+-h_t/2\Big]
\\
& \leq &
    2n_t e^{5\k h_t/4}
    \big[h_t^+ e^{-\k h_t/2 }\big]
    \big[(h_t^+ + h_t/2) e^{-\k h_t }\big]
\leq
e^{-h_t/8}/10
\end{eqnarray*}
where we applied \eqref{InegBeta*}
and $n_t =e^{o(h_t)}$ since $\phi(t)=o(\log t)$.
This together with \eqref{eqProbaB4SansFaggio} gives
$
    P\big[
        \overline{\B_5^{\ref{SubScetPreuveLemmasTechniquesht+}}}
    \big]
\leq
    e^{-\k h_t/8}/2.
$



Since \eqref{eqLemmaSansFaggio1}, \eqref{eqLemmaSansFaggio2} and \eqref{eqLemmaSansFaggio3}
are true for $1\leq i \leq n_t$
on
$
\B_1^{\ref{SubScetPreuveLemmasTechniquesht+}}
\cap \B_3^{\ref{SubScetPreuveLemmasTechniquesht+}}
\cap \B_5^{\ref{SubScetPreuveLemmasTechniquesht+}}$,
due to \eqref{eqIntegraleSansFaggio} for the second one,
the lemma is proved
\hfill $\Box$

\medskip
\noindent\textbf{Acknowledgment : } We would like to thank Dominique Lepingle for
pointing out reference \cite{RogPit}.
We are grateful to an anonymous referee for  comments that were useful to
        help improve the presentation of the paper.

\end{document}